\documentstyle[11pt]{article}
\newcommand{\mylabel}[1]{\label{#1}}

\newcommand{\numero}[1]{
\addtocounter{section}{1}
\begin{center}{\bf \thesection .\
#1\vspace{-.1in}}\end{center}
\setcounter{subsection}{0}
\setcounter{lemma}{0}\indent}

\newcommand{\subnumero}[1]{
\pagebreak[1]\begin{center}{\em #1}\nopagebreak\end{center}
}

\newcommand{\eop}{\hfill $\Box$\vspace{.1in}}

\newtheorem{lemma}{Lemma}[section]
\newtheorem{theorem}[lemma]{Theorem}

\newtheorem{corollary}[lemma]{Corollary}

\newtheorem{conjecture}{Conjecture}
\newtheorem{proposition}[lemma]{Proposition}

\newcommand{\cc}{{\bf C}}

\newcommand{\zz}{{\bf Z}}
\newcommand{\nn}{{\bf N}}

\newcommand{\Ee}{{\cal E}}
\newcommand{\Ff}{{\cal F}}
\newcommand{\Gg}{{\cal G}}

\newcommand{\Tt}{{\cal T}}
\newcommand{\Oo}{{\cal O}}

\newcommand{\Aa}{{\cal A}}
\newcommand{\Bb}{{\cal B}}
\newcommand{\Mm}{{\cal M}}

\newcommand{\Kk}{{\cal K}}
\newcommand{\Vv}{{\cal V}}

\newcommand{\Ll}{{\cal L}}
\newcommand{\Hh}{{\cal H}}

\newcommand{\Xx}{{\cal X}}

\newcommand{\Zz}{{\cal Z}}

\begin{document}

\section*{A relative notion of algebraic Lie group and applications to
$n$-stacks}

Carlos Simpson\newline
{\small Laboratoire Emile Picard
(UMR 5580 CNRS) \newline
Universit\'e Paul Sabatier\newline
31062 Toulouse CEDEX,  France}

\bigskip

Let $\Xx$ be the big etale site of schemes over $k=\cc$.
If $S$ is a scheme of finite type over $k$, let $\Xx /S$ denote the big
etale site of schemes over $S$.
The goal of this paper is to introduce a full subcategory of the category
of sheaves of groups on $\Xx /S$, which we will call  {\em the category
of presentable group sheaves} (\S 2),  with the following properties.
\newline
1. \, The category of presentable group sheaves contains those group sheaves
which are representable by group schemes of finite type over $S$ (Corollary
\ref{uvw}).
\newline
2. \, The category of presentable group sheaves is closed under kernel,
quotient (by a normal subgroup sheaf which is presentable), and extension
(Theorem \ref{I.1.e}).
\newline
3. \, If $S'\rightarrow S$ is a morphism then pullback takes presentable group
sheaves on $S$ to presentable group sheaves on $S'$ (Lemma \ref{I.1.h}).
\newline
4. \, If $S'\rightarrow S$ is a finite morphism then direct image takes
presentable group sheaves on $S'$ to presentable group sheaves on $S$ (Lemma
\ref{I.1.i}).
\newline
5. \, If $S=Spec (k)$ then presentable group sheaves are just group schemes of
finite type over $Spec (k)$ (Theorem \ref{I.1.m}).  In particular if $\Gg$ is
a
presentable group sheaf over any $S$ then the pullback to each point $Spec (k
)\rightarrow S$ is an algebraic group.
\newline
6. \, There is a notion of connectedness extending the usual
notion over $Spec(k )$ and compatible with quotients, extensions, pullbacks
and
finite direct images; and a presentable group sheaf $\Gg$ has a
largest connected presentable subsheaf $\Gg ^0\subset \Gg$ which we call the
{\em connected component} (Theorem \ref{I.1.o}).   \newline
7.  \, A presentable group sheaf $\Gg$ has a Lie algebra
object $Lie(\Gg )$ (Theorem \ref{lmn}) which is a vector sheaf with bracket
operation (see below  for a
discussion of the notion of vector sheaf---in the case $S=Spec (k)$ it is
the same thing as a finite dimensional $k$-vector space).
\newline
8. \, If $\Gg$  is a connected presentable group sheaf then $\Gg /Z(\Gg )$ is
determined up to isomorphism by the Lie algebra sheaf $Lie (\Gg )$ (where
$Z(\Gg
)$ denotes the center of $\Gg$).  This is Theorem \ref{abc} below.

\bigskip

We envision the category of presentable group sheaves as a generalisation
relative to an arbitrary base scheme $S$, of the category of algebraic Lie
groups
over $Spec (\cc )$.

We mention here a few questions related to the analogy with classical
algebraic
groups. Property 8 poses an obvious existence problem: given a Lie
algebra object
in the category of vector sheaves, does it come from a presentable
group sheaf with
vector sheaf center?  I don't know the answer to this question.
We do know,
however, that $Aut(L)$ is a presentable group sheaf (Lemma \ref{AutLie}).
Another
question is the existence of a ``universal covering'', i.e. a morphism
$\tilde{\Gg
}\rightarrow \Gg$ surjective with finite kernel such that for any other such
morphism $\Ff \rightarrow \Gg$ there is a factorization $\tilde{\Gg }
\rightarrow
\Ff \rightarrow \Gg$.  There are obvious questions about the generalisation of
the theory of representations to the case of presentable group sheaves.
The first
among these is whether there always exists a faithful representation into
$Aut(V)$
for $V$ a vector sheaf.  I suspect that the answer is no, but don't have a
counterexample.  For connected group sheaves this problem concerns only the
center, because we always have the adjoint representation of $\Gg$ on
$Lie (\Gg
)$. Beyond the question of the description of the representations, there
is also
the question of whether a suitable tannakian theory exists, namely given
a group
$\Gg \subset Aut (V)$, is $\Gg$ defined as the stabilizer of some
$\Gg$-invariant
sub-vector-sheaf $U$ in a tensor power of $V$?

The motivation for introducing presentable group sheaves comes from the
theory of
homotopy types over $Spec (\cc )$, or what Grothendieck
called ``schematization of
homotopy types'' in \cite{Grothendieck}.  We will discuss the application to
this theory at the end of the paper---note also that it is explained in
essentially the same way in \cite{kobe} where some applications to
nonabelian de
Rham cohomology are also announced.  Briefly, the considerations are as
follows.  A
homotopy type over $\Xx$ (which we call an  ``$n$-stack'') is a presheaf of
topological spaces on $\Xx$ satisfying a homotopic descent condition
(``fibrant''
in the terminology of Jardine \cite{Jardine1}, cf \cite{kobe}). This
condition is
the generalisation of the descent condition that goes into the definition of
$1$-stack. An $n$-stack or fibrant presheaf $T$ has homotopy sheaves
as follows.
First,  $\pi _0(T)$ is a sheaf of sets on $\Xx$. Then for $i\geq 1$ if
$S\in \Xx$
and $t\in T(S)$, $\pi _i (T|_{\Xx /S},t)$ is a sheaf of groups on $\Xx /S$
(abelian if $i\geq 2$). In the fibrant presheaf point of view, these
homotopy
sheaves are the sheafifications of the presheaves which one defines in
the obvious
way. These things satisfy the same sorts of properties as in the homotopy
theory of
spaces. In particular there are  notions of homotopy fiber products and (as
special cases) homotopy fibers and loop or path spaces. The homotopy
groups of the
homotopy fiber of a morphism fit into  the usual long exact sequence
(and there is
a similar exact sequence for homotopy fiber products in general).
There are also
notions of morphism spaces $Hom (T,T')$ which are spaces or $n$-groupoids
(depending on the point of view) and internal morphism objects
$\underline{Hom}(T,T')$ which are $n$-stacks whose global sections
are the morphism
spaces.

The main particularity of this situation is that $\pi _0(T)$ can be nontrivial
and not  just the union a set of points.  Because of this, one must consider
basepoints not only in $T(Spec (k ))$ but in $T(S)$ for any scheme $S$
(say, of
finite type) in order to get the full picture of $T$.  One is thus lead to
consider sheaves of groups on $\Xx /S$.

We would like to define a restricted class of $n$-stacks or
fibrant presheaves of
spaces  which we
will call {\em presentable}.  We would like this category to be closed under
homotopy fiber products and also under the truncation (or coskeleton)
operations of
eliminating the homotopy groups above a certain level.  From these
requirements it
follows that the condition for inclusion in the class of presentable presheaves
of spaces should be expressed solely in terms of the homotopy group sheaves.
From the exact sequences for homotopy fibers or more generally fiber products,
one can see that the category of group sheaves allowable as homotopy
group sheaves
of presentable spaces must be closed under kernel, cokernel and extension.
We would like our allowable group sheaves to be the algebraic Lie
groups when the
base space is $Spec (k )$, and of course for doing anything useful we need
notions of connectedness and an infinitesimal (Lie algebra) picture.  These are
the reasons which lead us to  look for a notion of sheaf of groups on
$\Xx /S$ with
the properties listed above.

I should add a note of caution about the terminology, for we propose the
terminology {\em presentable group sheaf} and also {\em presentable $n$-stack}.
If $\Gg$ is a sheaf of sets on $\Xx /S$ (i.e. a $0$-stack) which happens to have
a group structure, then the condition that $\Gg$ be presentable as a $0$-stack
is {\em not} the same as the condition that $\Gg$ be a presentable group sheaf
on $\Xx /S$.  The right way to think of a sheaf of groups is as
corresponding to a $1$-stack which we can denote $K(\Gg , 1)$ or $B\Gg$ (with a
morphism to $S$).  From this point of view the terminologies are compatible:
$\Gg$ is a presentable group sheaf over $S$ if and only if $K(\Gg , 1)$ is
a presentable $1$-stack.

Let's look more carefully at the reasoning that leads to our definition of
presentable $n$-stack.
 What are we going to do with a presentable $n$-stack
$T$?  If $W$ is (the  $n$-truncation
of) a finite CW complex considered as a constant $n$-stack on $\Xx$ then we can
look at the $n$-stack $Hom (W, T)$.  This is the {\em nonabelian cohomology of
$W$ with coefficients in $T$}.  If $T= K(\Oo , n)$ is the Eilenberg-MacLane
presheaf with homotopy group sheaf equal to the structure sheaf of rings $\Oo$
on $\Xx$ in degree $n$, then $\pi _0Hom (W, T)$ is just the cohomology $H^n(W,
\cc )$---or rather the sheaf on $\Xx$ represented by this vector space.
Similarly, if $G$ is a group scheme over $\cc$ then for $T=K(G, 1)=BG$ we get
that $Hom (W, G)$ is the moduli stack for flat principal
$G$-bundles or equivalently representations $\pi _1(W)\rightarrow G$.
We hope to
obtain an appropriate mixture of these cases by considering a more general
class
of $n$-stacks $T$.  In particular we would like to have a {\em Kunneth formula}
for two CW complexes $V$ and $W$,
$$
\underline{Hom} (U, \underline{Hom}(V,T))=\underline{Hom}(U\times V, T).
$$
One can imagine for example the problem of trying to compute the moduli stack of
flat principal $G$-bundles on $U\times V$ in terms of a Kunneth formula as
above.  One is forced to consider the cohomology of $U$ with coefficients in the
moduli stack $T'=\underline{Hom}(V,BG)$, and this stack is not necessarily
connected ($\pi _0(T')$ is roughly speaking the moduli space of principal
$G$-bundles).

The Kunneth formula is not an end in itself, as it is rare for a
space to decompose into a product.  It points the way to a
``Leray-Serre theory''
which could be more generally useful. If $W\rightarrow U$ is a morphism we would
be led to consider a relative morphism stack $T'=\underline{Hom}(W/U, T)
\rightarrow
U$ and then try to take the $n$-stack of sections $U\rightarrow T'$, a sort of
{\em nonabelian cohomology with twisted coefficients}.  I haven't fully
thought
about this yet (and in particular not about the de Rham
theory---see below---which
seems to be significantly more complicated than that which is needed in the
constant coefficient case, for example to make sense of the Kunneth formula).

The motivation for  all of this is to be able to do geometric versions of the
nonabelian cohomology in the case where $W$ is, say, a smooth projective
variety.  It is announced with some sketches of proofs in \cite{kobe}, how to
get a de Rham version of the morphism space $\underline{Hom}(W_{DR}, T)$ when
$T$ is a presentable $n$-stack.  We want of course to have the (analytic)
isomorphism between de Rham and Betti cohomology.  Needless to say, this will
not work for an arbitrary $n$-stack $T$ on $X$ (for example if one takes $T=W$
to a constant stack associated to a CW complex which is an algebraic variety
then
there will probably be nothing in $Hom (W_{DR}, W)$ corresponding  to the
identity in $Hom (W, W)$).  We need to impose conditions on $T$ which
guarantee
that it is reasonably close to the examples $K(\Oo , n)$ or $K(G, 1)$
given above
(in these cases, the de Rham-Betti isomorphism works as is already well
known).

As a first approach, the condition we want seems to be that the homotopy group
sheaves should be representable by group schemes over the base $S$.  In the case
where $T$ is the moduli stack of flat principal $G$-bundles on a space $V$,
encountered above when looking at the Kunneth formula, the $\pi _1$ sheaves are
indeed representable (the moduli stack is an algebraic stack).
Unfortunately the
condition of being representable is not stable under cokernels, but as explained
above this is important if we want our notion of good $n$-stack to be stable
under homotopy fiber products.

Before going directly to the conclusion that we need a category stable under
kernels, cokernels and extensions, we can analyze a bit more precisely just what
is needed.  Notice first of all that the algebraic de Rham theory is not
going to
work well in the case of higher cohomology with coefficients in the
multiplicative
group scheme, i.e. when $T= K({\bf G}_m, n)$ for $n\geq 2$.  I won't go into the
explanation of that here!  Thus, at least for the algebraic de Rham theory we
would like to have an appropriate notion of unipotent abelian group sheaf.  Not
yet having come up with a reasonable general theory of this, we can replace this
notion by the (possibly more restrictive) notion of {\em vector sheaf}.

The notion of vector sheaf  is explained in \S 4 below.
The reader may actually wish to start by reading this section, since the
theory of
vector sheaves
is in some sense a paradigm, applicable only for abelian group sheaves, of
what we are trying to do in general. The notion of vector sheaf was
introduced by
A. Hirschowitz \cite{Hirschowitz} who called it ``U-coherent sheaf''.
He defined
the category of U-coherent sheaves as the smallest abelian category
of sheaves of
abelian groups containing the coherent sheaves (note that the category of
coherent sheaves is not abelian on the big etale site or any big site).
We take a
more constructive approach, defining the notion of vector sheaf in terms of
presentations, although in the end the two notions are equivalent.
The notion of vector sheaf doesn't work too nicely in characteristic $p>0$,
basically because the Frobenius automorphism of the sheaf $\Oo$ is not
linear, so
the linear structure is no longer encoded in the sheaf structure.  As
we try in the
beginning of the paper to put off the hypothesis of characteristic zero as long
as possible, and as the notion of vector sheaf comes into the analysis
 at a later
stage (the infinitesima study related to properties 7 and 8 listed above), I
have decided not to put the section on vector sheaves at the beginning.
Still, it
is essentially self-contained for the reader who wishes to start there.

In considering the algebraic
de Rham theory we will only be looking at $n$-stacks $T$ with $\pi _i(T,t)$ a
vector sheaf on $S$ for $t\in T(S)$ and $i\geq 2$.  What does this mean for
our restriction on $\pi _1(T, t)$? Going back to the question of stability
under fiber products, we see from looking at the long exact homotopy sheaf
sequence that the minimum that is absolutely necessary is that our class of
group
sheaves $G$ be stable under central extension by a vector sheaf.  On the other
hand it also must be stable under taking kernels.

One could thus hope to make
good on a {\em minimalist approach} saying that we should look at the
category of
group sheaves generated by representable group sheaves (affine, say---this again
would be needed to make the de Rham theory work), under the operations
of kernel
and central extension by a vector sheaf.  A vector sheaf always has a
presentation
as the cokernel of a morphism of {\em vector schemes}, i.e. representable
vector-space objects over the base $S$ (these are sometimes called {\em linear
spaces} in the complex analytic category \cite{Grauert} \cite{Fischer}
\cite{Axelsson-Magnusson}). The most natural approach then is to say,
suppose a
group sheaf $G$ has a presentation as the cokernel of a morphism $F_2
\rightarrow
F_1$ of representable group sheaves, and suppose $E$ is a central
extension of $G$
by  a vector sheaf $U$ which is itself the cokernel of a morphism $V_2
\rightarrow
V_1$ of vector schemes.  Then try to combine these into a presentation of $E$
with, for example, a surjection $V_1\times F_1\rightarrow E$.  The
problem (which
I was not able to solve although I don't claim that it is impossible)
is then to
lift the multiplication of $E$ to an associative multiplication on $V_1\times
F_1$.

As I didn't see how to do this, a slightly more general approach was needed,
wherein we consider groups which have presentations by objects where the
multiplication lifts but not necessarily to a multiplication satisfying the
associativity property.  This is the reasoning that leads to the definition of
{\em $S$-vertical morphism:} \, a morphism where one can lift things such as
multiplications in a nice way cf \S 2. We finally come to the definition of
{\em presentable group sheaf} as a group sheaf $G$ which admits a vertical
surjection $X\rightarrow G$ from a scheme of finite type over $S$, and such
that there is a  vertical surjection $R\rightarrow X\times _{G}X$ again from
a scheme of finite type over $S$.

One could, on the other hand, take a {\em  maximalist approach} and try to
include anything that seems vaguely algebraic.  This would mean, for example,
looking at sheaves $G$ such that there are surjections (in the etale sense,
although not necessarily etale morphisms) $X\rightarrow G$ and $R\rightarrow
X\times _GX$ with $X$ and $R$ schemes of finite type over the base $S$.  We
call this condition P2.  This might also work (in fact it might even be the
case that a P2 group sheaf is automatically presentable).  However, I was not
able to obtain a reasonable infinitesimal analysis which could lead, for
example, to the notion of connected component---though again, I don't claim
that
this could never work.

In a similar vein, one might point out that there is a fairly limited range of
situations in which we use the lifting properties going into the definition of
verticality.
I have chosen to state the
condition of verticality in what seems to be the most natural setting,
but this leads to requiring that many more lifting properties be satisfied than
what we actually use. One could rewrite the definition of verticality
to include only those lifting properties that we use afterward.  It might be
interesting to see if this change makes any difference in which group sheaves
are allowed as presentable.

All in all, the definitions we give here of presentable group sheaf and of
presentable $n$-stack are first attempts at giving useful and reasonable
notions,
but is is quite possible that they would need to be altered
in the future in view
of applications.

A word about the characteristic of the ground field (or base scheme).  While our
aim is to work over a field of characteristic zero, there are certain parts of
our discussion valid over any base scheme, namely
those concerning the abstract method for defining conditions of presentability.
When it comes down to finding conditions which result in a nice theory (and in
particular which result in a theory having the required local structure) we must
restrict ourselves to characteristic zero.  It is possible that a variant could
work nicely in positive characteristic, so we will present the first part of the
argument concerning the definition of presentability (which is valid over
any base
scheme), in full generality (\S 1) before specifying in characteristic zero
which
morphisms we want to use in the presentations (\S 2). Actually the definition
given in \S 2 works in any characteristic but we can only prove anything about
local properties in characteristic zero (\S\S 4-9), so it is probably the
``right'' definition only in characteristic zero. With an appropriate
different
definition of verticality (certainly incorporating divided powers) what we do in
these later sections might be made to work in any characteristic.

\subnumero{Notations}
We fix a noetherian ground
ring $k$, for sections 1-3.  From section 4 on, we assume that $k$ is an
uncountable field of characteristic zero,  and we may when necessary assume that
the ground field is $k=\cc$.

Let $\Xx$ denote the site of noetherian schemes over $k$ with the etale topology
(this is known as the ``big etale site'').

If $S\in \Xx$ then we denote by $\Xx /S$ the site of schemes over $S$ (again
with the etale topology).

A {\em sheaf} on $\Xx$ means (unless otherwise specified) a sheaf of sets.  For
a sheaf of groups, we sometimes use the terminology {\em group sheaf}.

We will confuse notations between an object of $\Xx$ and the sheaf it
represents.

Denote by $\ast$ the sheaf on $\Xx$ with values equal to the one-point set; it
is represented by $Spec (k)$ (and we can interchange these notations at will).

If $S$ is a sheaf on $\Xx$ (most often represented by an object) then we have
the site $\Xx /S$ of objects of $\Xx$ together with morphisms to $S$.
There is an
equivalence between the notions of sheaf on $\Xx /S$ and sheaf on $\Xx$ with
morphism to $S$.  Since we will sometimes need to distinguish these, we
introduce
the following notations.

If $\Ff$ is a sheaf on $\Xx$ then its {\em restriction up} to $\Xx /S$ is
denoted
by $\Ff |_{\Xx /S}$, with the formula
$$
\Ff |_{\Xx /S}(Y\rightarrow S)= \Ff (Y).
$$

If $\Ff$ is a sheaf on $\Xx /S$ then we denote by $Res_{S/\ast}\Ff$ the
corresponding sheaf on $\Xx$ together with its morphism
$$
Res_{S/\ast}\Ff \rightarrow S.
$$
It is defined by the statement that $Res_{S/\ast}\Ff (Y)$ is equal to the set of
pairs $(a, f)$ where $a: Y\rightarrow S$ and $f\in \Ff (
Y\stackrel{a}{\rightarrow}
S)$. We call this the {\em restriction of $\Ff$ from $S$ down  to $\ast$}.

More
generally if $S'\rightarrow S$ is a morphism  and if $\Ff$ is a sheaf on $\Xx
/S'$ then we obtain a sheaf $Res _{S'/S}\Ff$ on $\Xx/S$ called the
{\em restriction of $\Ff$ from $S'$ down  to $S$}.

The operations of restriction up and restriction down are not inverses: we have
the formula, for a sheaf $\Ff $ on $\Xx /S$,
$$
Res _{S'/S}(\Ff |_{\Xx /S'}) = \Ff \times _SS' .
$$

On the other hand, suppose $p:\Ff \rightarrow S'$ is a morphism of sheaves on
$\Xx/S$.  Then we denote by $\Ff /S'$ the corresponding sheaf on $\Xx /S'$ (the
data of the morphism is implicit in the notation). It is defined by the
statement
that $\Ff /S' ( Y\rightarrow S')$ is equal to the set of $u \in \Ff
(Y\rightarrow
S)$ such that $p(u)\in S'(Y\rightarrow S)$ is equal to the given morphism
$Y\rightarrow S'$.   For another point of view note that there is a tautological
section of $(S'/S)|_{ \Xx /S'}$,  and $\Ff /S'$
is the preimage of this section in $\Ff |_{\Xx /S'}$.

As a special case we get that if $\Ff$ is a sheaf on $\Xx = \Xx /\ast$ with a
morphism $\Ff \rightarrow S$ then we obtain a sheaf $\Ff /S$ on $\Xx /S$.

The operations
$$
(\Ff \rightarrow S')\mapsto \Ff /S'
$$
from sheaves on $\Xx /S$ with morphisms to $S'$ to sheaves on $\Xx /S'$,
and
$$
\Ff ' \mapsto (Res _{S'/S}\Ff ' \rightarrow S'
$$
from sheaves on $\Xx /S'$ to sheaves on $\Xx /S$ with morphisms to $S'$, are
inverses. For this reason it is often tempting to ignore the strict notational
convention and simply use the same notations for the two objects.  This is not
too dangerous except in the last section of the paper where we will try to be
careful.

If a sheaf $\Ff$ on $\Xx$ is representable by an object $F\in \Xx$ and
if $F\rightarrow S$ is a morphism then $\Ff /S$ is representable by the same
object $F$ together with its morphism, considered as an object of $\Xx /S$.
For this reason we will sometimes drop the notation $\Ff /S$ and just denote
this as $\Ff$ when there is no risk of confusion (and in fact the attentive
reader will notice that even in the definition two paragraphs ago we have
written
$S'$ when we should have written $S'/S$ in the first sentence...but the second
version would have been impossible because not yet defined...!)

Finally there is another natural operation: suppose $\pi : S'\rightarrow S$ is a
morphism and $\Ff$ is a sheaf on $\Xx /S'$. Its {\em direct image} is the sheaf
$\pi _{\ast}(\Ff )$ defined by the statement that
$$
\pi _{\ast}(\Ff
)(Y\rightarrow S):= \Ff (Y\times _SS' \rightarrow S').
$$
This is {\em not} the same thing as the restriction down from $S'$ to $S$.
Think of the case where $S$ is one point and $S'$ is a collection of several
points.  The value of $Res_{S'/S}(\Ff )$ at $S$ is the {\em union} of the values
of $\Ff$ over the points in $S'$ whereas the value of $\pi _{\ast}(\Ff )$ at $S$
is the {\em product} of the values of $\Ff$ at the points in $S'$.

\numero{Presentability conditions for sheaves}

We will define several conditions, numbered $P1$, $P2$, $P4(\Mm )$,
$P5(\Mm )$ (whereas two other conditions $P3$ and $P3\frac{1}{2}$ will be
defined
later, in \S 2). The last two depend on a choice of a class $\Mm$ of morphisms
in $\Xx$ subject to certain properties set out below.   In the upcoming section
we then specify which class $\Mm$ we are interested in (at least in
characteristic zero), the class of {\em vertical morphisms}.  Since the
preliminary results depend only on the formal properties of $\Mm$ we thought it
might be useful to state them in general rather than just for the class of
vertical morphisms, this is why we have the seeming complication of introducing
$\Mm$ into the notations for our properties.

We also introduce {\em boundedness conditions} denoted $B1$ and $B2$.  These
conditions sum up what is necessary in order to be able to apply Artin
approximation.

Fix a base scheme $S\in \Xx$.
In what follows, we work in the category of sheaves on $\Xx /S$. Thus a sheaf
is supposed to be on $\Xx /S$ unless otherwise specified.

\noindent
{\bf P1.}\,\, We say that $\Ff$ is {\em P1} if there is a surjective morphism of
sheaves $X\rightarrow \Ff$ where $X$ is represented by a scheme of finite type
over $S$. We may assume that $X$ is affine.

\noindent
{\bf P2.}\,\, We say that $\Ff$ is {\em P2} if there are surjective morphisms of
sheaves $X\rightarrow \Ff$ and $Y\rightarrow X\times _{\Ff}X$ where $X$ and $Y$
are  represented
by  schemes of finite type over $S$. We may assume that $X$ and $Y$ are affine.

\begin{lemma}
\mylabel{I.t}
If $G$ is a sheaf of groups which is P1, and $G$ acts on a sheaf $F$ which is
P2, then the quotient sheaf $F/G$ is again P2.
\end{lemma}
{\em Proof:}
Choose surjections $\varphi :X\rightarrow F$ and $(p_1,p_2):Y\rightarrow X\times
_FX$. The action is a map $G\times F\rightarrow F$, and
we can choose a surjection $(q_1,q_2):W\rightarrow G\times X$ (with $W$ an
affine scheme, by condition P1 for $G$), such that the action lifts to a map
$m:W\rightarrow X$.  There is obviously a surjection $X\rightarrow F/G$.  We
have a map
$$
W\times _X Y\rightarrow X\times X
$$
(where the maps used in the fiber product are $m:W\rightarrow X$ and
$p_1:Y\rightarrow X$), defined by
$$
(w,y)\mapsto (q_2(w), p_2(y)).
$$
This map surjects onto the fiber  product $X\times _{F/G}X$. It clearly maps
into this fiber product.  The map is surjective because if $(x,x')\in X\times
X$ with $g\varphi (x)=\varphi (x')$ then for a point $w$ of $W$  lying above
$(g,x)$ we have $\varphi (m(w))= g\varphi (x)=\varphi (x')$; in particular there
is a point $y$ of $Y$ with $p_1(y)=m(w)$ and $p_2(y)=x'$, so the point
$(w,y)$ maps to $(x,x')$.
Our surjection
$$
W\times _X Y\rightarrow X\times _{F/G}X
$$
now shows that $F/G$ is P2.
\eop

{\em Remark:}
These conditions are independent of base scheme $S$ for finite-type morphisms.
More precisely if $S'\rightarrow S$ is a morphism of finite type and if
$\Ff '$ is
a sheaf on $\Xx /S'$ then denoting by $\Ff = Res _{S'/S}\Ff '$ its restriction
down to $S$ we have that $\Ff $ is $P1$ (resp. $P2$) if and only
if $\Ff '$ is $P1$ (resp. $P2$).

\subnumero{Boundedness conditions}

We consider the following boundedness conditions for a sheaves on $\Xx$. These
two conditions are designed to contain exactly the information needed to apply
the Artin approximation theorem \cite{Artin}.
\newline {\bf B1.} \,\,
We say that a sheaf  $\Ff$ is {\em B1} if, for
any $k$-algebra $B$, we have that
$$
\lim _{\rightarrow}\Ff (Spec (B')) \rightarrow \Ff (Spec (B))
$$
is an isomorphism, where the limit is taken over the subalgebras $B'\subset
B$ which are of finite type over $k$. This is equivalent to the local finite
type
condition of Artin \cite{Artin}.

\noindent
{\bf B2.} \,\,  We say that a sheaf $\Ff$ is {\em B2} if, for any complete local
ring $A$, we have that the morphism
$$
\Ff (Spec (A)) \rightarrow \lim _{\leftarrow} \Ff (Spec (A/{\bf m} ^i)
$$
is an isomorphism.

The {\bf Artin approximation theorem} (\cite{Artin}) can now be stated as
follows.

{\em
Suppose $\Ff$ is a sheaf of sets which is B1 and B2.  If $S=Spec (A)$ is an
affine scheme with point $P\in S$ corresponding to a maximal ideal ${\bf m}
\subset
A$ then for any
$$ \eta \in
\lim _{\leftarrow} \Ff (Spec (A/{\bf m} ^i))
$$
and for $i_0\geq 0$ there exists an etale neighborhood $P\in
S' \rightarrow S$ and
an element $\eta ' \in \Ff (S')$ agreeing with $\eta$ over  $Spec (A/{\bf m}
^{i_0})$.
}

\begin{lemma}
\mylabel{I.t.1}
1.\,\, If $\Ff$ and $\Gg$ are B1 (resp. B2) and $f,g$ are two morphisms from
$\Ff$ to $\Gg$ then the equalizer is again B1 (resp. B2).
\newline
2.\,\, Suppose $\Ff\rightarrow
\Gg$ is a surjective morphism of sheaves.  If $\Ff$ and $\Ff \times _{\Gg}\Ff$
are B1 then $\Gg$ is B1.
\end{lemma}
{\em Proof:}
Fix $S=Spec (A)$ and $\{ B_i\}$ our directed system of $A$-algebras.  Let $B=
\lim _{\rightarrow}B_i$.
Suppose $\eta \in \Gg (B)$.  There is a natural morphism
$$
\lim _{\rightarrow} \Gg (B_i)\rightarrow \Gg (B).
$$
First we prove injectivity. Suppose $\varphi , \psi \in \Gg (B_i)$ map to the
same element of $\Gg (B)$.  We may choose an etale surjection of finite type
$Spec (B'_i)\rightarrow Spec (B_i)$ such that the restrictions $\varphi '$ and
$\psi '$ lift to elements  $u,v\in \Ff (B_i)$.  Their images in $\Ff (B')$
give a
point $(u,v)_{B'}$ in $\Ff \times _{\Gg} \Ff (B')$ (here $B':= B\otimes
_{B_i}B'_i$).  By the condition B1 for the fiber product, there is a $j\geq i$
such that this point comes from a point $\eta \in \Ff \times _{\Gg}\Ff (B'_j)$.
On the other hand, note that the product $\Ff \times \Ff$ is B1.  The image of
$\eta$ in $\Ff \times \Ff (B')$ is the same as that of $(u,v)$; and
by the B1 condition for the product, there is $k\geq j$ such that the image of
$\eta$ in $\Ff \times \Ff (B'_k)$ is equal to the image of $(u,v)$.
In particular, $(u|_{Spec (B'_k)},v|_{Spec (B'_k)})$ lies in the fiber product
$\Ff \times _{\Gg}\Ff (B'_k)$. In other words, $u$ and $v$ have the same images
in $\Gg (B'_k)$.  These images are the restrictions of the original $\varphi ,
\psi$.  Since $Spec (B'_k)\rightarrow Spec (B_k)$ is an etale surjection, the
images of $\varphi$ and $\psi $ in $\Gg (B_k)$ are the same.  This proves the
injectivity.

Now we prove surjectivity.
Then there exists an etale surjection of
finite type
$$
Spec (B')\rightarrow Spec (B)
$$
such that $\eta |_{Spec (B')}$ comes
from an element $\rho \in \Ff (B')$.  The functor ``etale surjections of
finite type'' is itself B1, so there is an etale $Spec (B'_i)\rightarrow Spec
(B_i)$ inducing $B'$.  Then $B'=\lim _{\rightarrow} B'_j$ where $B'_j=
B_j\otimes _{B_i}B'_i$ for $j\geq i$.  By the property B1 for $\Ff$ there is
some $j$ such that $\rho$ comes from $\rho _j\in \Ff ( B'_j)$.  We obtain an
element $\eta '_j\in \Gg (B'_j)$ mapping to $\eta ':=\eta |_{Spec (B')}$. The
two
pullbacks of $ \eta '$ to $Spec (B'\otimes _BB')$ are equal.  Note that
$$
B'\otimes _BB' = \lim _{\rightarrow} B'_k\otimes _{B_k}B'_k,
$$
so by the above injectivity, there is some $k$ such that the two pullbacks of
$\eta _j|_{Spec (B'_k)}$ to $Spec (B'_k\otimes _{B_k}B'_k)$ are equal.  Now
the sheaf condition for $\Gg$ means that $\eta _j|_{Spec (B'_k)}$ descends to an
element $\eta _k\in \Gg (B_k)$.  The restriction of $\eta _k$ to $B'$ is equal
to the restriction of $\eta$, so the sheaf condition for $\Gg$ implies that
the restriction of $\eta _k$ to $Spec (B)$ is $\eta$.
\eop

{\em Remark:}  The direct product of a finite number of B1 (resp. B2) sheaves is
again B1 (resp. B2) so part 1 of the lemma implies that the properties B1 and B2
are maintained under fiber products.

\begin{theorem}
\mylabel{I.t.2}
Suppose $\Ff$ is a sheaf which is P2.  Then $\Ff$ is B1.  If the ground field
is uncountable, then $\Ff$ is B2.
\end{theorem}
{\em Proof:}
The condition B1 follows from the previous lemma.  Indeed, let $X\rightarrow
\Ff$ and $R\rightarrow X\times _{\Ff}X$ be the morphisms given by the property
P2, with $X$ and $R$ of finite type (in particular, B1).  Note that $R\times
_{X\times _{\Ff}X}R= R\times _{X\times X}R$ is a scheme of finite type, so the
lemma implies that $X\times _{\Ff}X$ is B1; another application of the lemma
then shows that $\Ff$ is B1.

For B2, let $S=Spec (A)$ with $A$ a complete local ring, and let $S_n:= Spec
(A/{\bf m}^{n+1})$.
Let $X\rightarrow
\Ff$ and $R\rightarrow X\times _{\Ff}X$ be the morphisms given by the property
P2, with $X$ and $R$ of finite type over $S$.  Schemes of finite type are B2.
We show surjectivity of the map
$$
\Ff (S)\rightarrow \lim _{\leftarrow} \Ff (S_n).
$$
Suppose $(\varphi _n )$ is a compatible system of elements of $\Ff (S_n)$. Let
$$
E_n := \{ x_n \in X(S_n):\;\;\; x_n \mapsto \varphi _n \} .
$$
Note that $E_n$ is a nonempty closed subset of the scheme $X(S_n )$ (that is,
the
scheme whose $Spec (k)$-valued points are $X(S_n)$).  Let
$$
E'_n:= \bigcap _{m\geq n} {\rm im}(E_m \rightarrow E_n);
$$
this is an intersection of a decreasing family of nonempty constructible subsets
of $E_n$.  Since $k$ is uncountable, this intersection is nonempty.  Indeed,
the closures of the images form a decreasing family of closed sets, which
stabilizes by the noetherian property of $E_n$; then within this closed subset,
the dense constructible subsets contain open sets which are complements of
proper closed subsets.  The union of countably many proper closed subsets is a
proper subset, so the intersection of the open complements is nonempty. (Note
however that $E'_n$ is not necessarily constructible).

The morphism $E'_{n+1} \rightarrow E' _n$ is surjective.  To see this, suppose
$u\in E' _n$. We can consider the subsets
$$
D_m := \{ v\in E_m, \;\; v\mapsto u\} .
$$
These are closed subsets of $E_m$, nonempty by the condition $u\in
E'_n$. Let $D' _{n+1}= \bigcap _{m\geq n+1} {\rm im} (D_m \rightarrow D_{n+1})$.
By the same proof as above, $D'_{n+1}$ is nonempty.  But it is contained in
$E'_{n+1}$ and maps to $u\in E'_n$.

The surjectivity of the maps implies that the inverse limit $\lim _{\leftarrow}
E'_n $ is  nonempty.  It is a subset of $\lim _{\leftarrow} X(S_n)=X(S)$,
consisting of elements mapping to $(\varphi _n)$ in $\lim _{\leftarrow}\Ff
(S_n)$. (In fact, this subset is equal to the inverse image of $(\varphi _n)$.)
We obtain an element of $X(S)$, hence an element of $\Ff (S)$, mapping to
$(\varphi _n)$.  This proves surjectivity.  Note that this part of the proof
only used property P1 for $\Ff$.

We now prove injectivity.  Note that $X\times _{\Ff}X$ is P1, so by the proof
above, we obtain surjectivity of the morphism
$$
X\times _{\Ff}X(S)\rightarrow \lim _{\leftarrow}X\times _{\Ff}X(S_n).
$$
Suppose two elements $u,v\in G(S)$ go to the same element of $G(S_n)$ for all
$n$ (we write this $u_n=v_n$).  We can lift them to elements $x,y\in X(S)$, and
we obtain a compatible sequence of elements $(x_n, y_n)\in X\times_{\Ff}X
(S_n)$.
By the surjectivity of the above morphism, there is an element $(x',y')\in
X\times _{\Ff}X(S)$ with $x'_n=x_n$ and $y'_n=y_n$. The images $u'$ and $v'$
of $x'$ and $y'$ in $\Ff (S)$ are equal. By
the B2 property for $X$, this implies that $x'=x$ and $y'=y$, which shows that
$u=v$.
\eop

We have the following Krull-type property.

\begin{lemma}
\mylabel{Krull}
Suppose $\Ff$ is a sheaf which is B1 and B2.  Then for any scheme $S$ of
finite type the natural morphism is an injection
$$
\Ff (S ) \hookrightarrow \prod _{{\rm Art.} \, S'\rightarrow S} \Ff (S' )
$$
where the product is taken over $S'\rightarrow S$ which are artinian and of
finite type.
\end{lemma}
{\em Proof:}
Suppose $f,f'\in \Ff (S)$ agree over all artinian subschemes.  Let $\Gg = S
\times _{\Ff } S$ be the fiber product where $f$ and $f'$ provide the two
morphisms from $S$ to $\Ff$.  Then $\Gg$ is B1 and B2 (by the remark following
Lemma \ref{I.t.1}).  But $\Gg$ has a (unique) section over any artinian
$S'\rightarrow S$ and applying B2, B1 and Artin approximation we obtain
sections of $\Gg$ over an etale covering of $S$.
\eop

\subnumero{Choice of a class of morphisms $\Mm$}

Fix a base scheme $S\in \Xx$.
We assume fixed for the rest of this section a subset $\Mm\subset
Mor (\Xx /S)$ of
morphisms in $\Xx /S$, containing the identities and closed under composition
(i.e.
$\Mm$ is the set of morphism of a subcategory of $\Xx /S$) subject to the
following
axioms:  \newline
{\bf M1}\,\, If $a$ and $b$ are composable morphisms such that $a$ and $ba$
are in
$\Mm$, and $a$ is surjective, then $b$ is in $\Mm$.
\newline
{\bf M2}\,\, Compatibility with base change: if $\Ff \rightarrow \Gg$ is an
$\Mm$-morphism and $\Hh \rightarrow \Gg$  any morphism, then $\Ff
\times _{\Gg}\Hh \rightarrow \Hh$ is an $\Mm$-morphism; and conversely if
$a:\Ff\rightarrow \Gg$ is a morphism such that $\Ff \times _\Gg
Y\rightarrow Y$ is
in $\Mm$ for every $S$-scheme and morphism $Y\rightarrow \Gg$, then $a$ is in
$\Mm$.
\newline
{\bf M3}\,\, Etale morphisms between schemes are in $\Mm$.

{\em Remark:}
It follows from these axioms that the direct product of morphisms in $\Mm$ is
again a morphism in $\Mm$.

In the next section we will specify a certain such subcategory $\Mm$, the class
of {\em vertical morphisms}, and show that it satisfies these axioms.  But there
may be other interesting examples of such a class of morphisms $\Mm$ to
which the
following definitions and lemmas could be applied.

We can now extend our list of presentability properties which refer to the
class $\Mm$. We use the notation $\Mm$-morphism for morphism lying in $\Mm$.
The gap in the numbering is to leave a place for the property $P3$ later.
This property (which is absolute rather than relative to a base scheme $S$)
will come up only at the end of the paper, but it turns out to be more
logical to
number it in between $P2$ and $P4$ (this is the numbering used in \cite{kobe}).

\noindent
${\bf P4(\Mm )}$\,\, We say that a sheaf $\Ff$ is $P4(\Mm )$ if there exist
surjective $\Mm$-morphisms
$$
X\rightarrow \Ff
$$
and
$$
R\rightarrow X\times _{\Ff}X
$$
with $X$ and $R$ represented by affine schemes of finite type over $S$.

\noindent
${\bf P5(\Mm )}$\,\,
We say that $\Ff$ is $P5(\Mm )$ if it is $P5(\Mm )$ and if, in addition, the
structural morphism $\Ff \rightarrow S$ is in $\Mm$.

\begin{lemma}
\mylabel{I.z.1}
If $\Ff$ and $\Gg$ are $P4(\Mm )$ (resp. $P5(\Mm )$) then so is
$\Ff \times _S\Gg$.
\end{lemma}
{\em Proof:}  The presentation is just the product of the presentations for
$\Ff$ and $\Gg$.
\eop

\subnumero{Kernels and extensions}

\begin{lemma}
\mylabel{I.1.a}
If $f,g:\Gg \rightarrow \Hh$ are two morphisms, and if $\Gg$ and $\Hh$ are
$P4(\Mm )$,
then the equalizer $\Ff$ is $P4(\Mm )$.
\end{lemma}
{\em Proof:}
Let $X\rightarrow \Hh$, $R\rightarrow X\times _{\Hh}X$, $Z\rightarrow \Gg$ and
$T\rightarrow Z\times _{\Gg}Z$ be  $\Mm$-morphisms with $X$, $R$, $Z$ and $T$
schemes of finite type over $S$.  Assume that we have liftings
$f',g': Z\rightarrow X$ of $f$ and $g$. Set
$$
W:= Z\times _{X\times _SX}R.
$$
It is a scheme of finite type over $S$. Note that the composed map
$Z\times _{\Gg}
\Ff \rightarrow Z\rightarrow X\times _SX$ factors through $X\times _{\Hh}X$, and
we have
$$
W= (Z\times _{\Gg} \Ff )\times _{X\times _{\Hh}X}R.
$$
From this and property $M2$, it is clear that the morphism $W\rightarrow \Ff$ is
surjective and in $\Mm$. Now set
$$
V:= (W\times _SW)\times _{Z\times _SZ}T.
$$
Again, this is of finite type over $S$.  We have
$$
W\times _{\Ff}W= W\times _{\Gg}W = (W\times _SW)\times _{Z\times _SZ}(Z\times
_{\Gg}Z).
$$
Therefore
$$
V= (W\times _{\Ff} W)\times _{Z\times _{\Gg}Z}T.
$$
From this and property $M2$ it is clear that the morphism
$V\rightarrow W\times _{\Ff}W$ is surjective and in $\Mm$.
We obtain the property $P4(\Mm )$ for $\Ff$.
\eop

\begin{corollary}
\mylabel{I.1.a.1}
If $\Ff\rightarrow \Hh$ and $\Gg \rightarrow \Hh$ are two morphisms between
$P4(\Mm )$
sheaves, then the fiber product $\Ff \times _{\Hh}\Gg$ is $P4(\Mm )$.
\end{corollary}
{\em Proof:}
The fiber product is the equalizer of the two morphisms $\Ff \times _S\Gg
\rightarrow \Hh$.
\eop

\begin{lemma}
\mylabel{I.1.b}
Suppose $\Hh$ is a group sheaf which is $P5(\Mm )$.  If $\Hh$ acts freely on a
sheaf $\Gg$ with quotient $\Ff = \Gg /\Hh$, then the morphism $\Gg \rightarrow
\Ff$ is in $\Mm$.
\end{lemma}
{\em Proof:}
Make a
base change by a scheme $Y\rightarrow \Ff$.  Let $\Gg '':= \Gg \times
_{\Ff}Y$.  Then $\Hh$ acts freely on $\Gg ''$ with quotient $Y$.  Since the
morphism $\Gg '' \rightarrow Y$ is surjective in the etale topology, we may
find an etale morphism (of finite type and surjective) $Y'\rightarrow Y$ such
that the base change $\Gg ^3$ of $\Gg ''$ to $Y'$ admits a section.  Then  $\Gg
^3= Y'\times _S\Hh$.  In particular, the morphism $\Gg ^3\rightarrow Y'$ is
in $\Mm$, hence also the morphism $\Gg ^3 \rightarrow Y$.  Finally, the morphism
$\Gg ^3 \rightarrow \Gg ''$ is surjective, since $Y'\rightarrow Y$ is
surjective,
and is an $\Mm$-morphism because it becomes an etale morphism after
base change to
any scheme.  By property $M1$, the morphism $\Gg '' \rightarrow Y$ is
in $\Mm$; then by $M2$ the morphism $\Gg \rightarrow \Ff$ is
in $\Mm$.
\eop

\begin{lemma}
\mylabel{I.1.c}
Suppose $\Gg$ is a $P4(\Mm )$ sheaf, and suppose $X\rightarrow \Gg$ is a
morphism with $X$ a scheme of finite type over $S$.
Then there exists a surjective $\Mm$-morphism $R\rightarrow X\times _{\Gg}X$
with $R$ a scheme of finite type over $S$.
\end{lemma}
{\em Proof:}
Let $Y\rightarrow \Gg$ and $Q\rightarrow Y\times _{\Gg}Y$ be the surjective
$\Mm$-morphisms.
There is an etale surjection $X'\rightarrow X$ such that the lifting
$X'\rightarrow Y$ exists.  Note that
$$
X'\times _{\Gg} X' = (X' \times _S X')\times _{Y\times _SY}(Y\times _{\Gg}Y).
$$
We get that
$$
R:= (X'\times _{\Gg}X')\times _{Y\times _{\Gg}Y}Q=
(X' \times _S X')\times _{Y\times _SY}Q
$$
is a scheme of finite type.  But also the morphism
$$
R= (X'\times _{\Gg}X')\times _{Y\times _{\Gg}Y}Q\rightarrow
X'\times _{\Gg}X'
$$
is  in $\Mm$, by property $M2$.  Finally,
$$
X'\times _{\Gg} X' = (X\times _{\Gg}X) \times _{X\times _SX} X'\times _SX'
$$
and $X'\times _SX'\rightarrow X\times _SX$ is an $\Mm$-morphism by $M3$ and the
remark following the properties $M$.  Thus $X'\times _{\Gg} X'\rightarrow
X\times
_{\Gg} X$ is in $\Mm$ (it is also surjective), so the surjection $R\rightarrow
X\times _{\Gg}X$ is in $\Mm$. \eop

\begin{theorem}
\mylabel{I.1.d}
Suppose $\Hh$ is a group sheaf which is $P5(\Mm )$, and suppose that $\Hh$
acts freely
on a sheaf $\Gg$ with quotient $\Ff = \Gg /\Hh $.
Then $\Ff$ is $P4(\Mm )$ (resp. $P5(\Mm )$) if and only if $\Gg$ is $P4(\Mm )$
(resp. $P5(\Mm )$).
\end{theorem}
{\em Proof:}
By the lemma, the morphism $\Gg \rightarrow \Ff$ is in $\Mm$.   If $\Gg$ is
$P4(\Mm )$
then there is a surjective $\Mm$-morphism $X\rightarrow \Gg$
with $X$ a scheme of finite type over $S$.  The morphism $X\rightarrow
\Ff$ is then surjective and in $\Mm$.
Let $Y\rightarrow \Hh$ be a surjective $\Mm$-morphism.
Now we have a surjective $\Mm$-morphism
$$
X\times _SY\rightarrow X\times _S\Hh = X\times _{\Ff}\Gg ,
$$
and another surjective $\Mm$-morphism
$$
X\times _{\Ff}\Gg \rightarrow \Ff \times _{\Ff}\Gg =\Gg .
$$
Apply the previous lemma to the composition of these two morphisms, using the
property $P4(\Mm )$ of $\Gg$.  We obtain the existence of a surjective
$\Mm$-morphism
$$
T\rightarrow (X\times _SY)\times _{\Gg} (X\times _SY)
$$
with  $T$ a scheme of finite type over $S$.
On the other hand, note that we have a surjective $\Mm$-morphism
$$
X\times _{\Ff}X\times _S\Hh=(X\times _{\Ff}\Gg )\times _{\Gg}
(X\times _{\Ff}\Gg )\rightarrow X\times _{\Ff} X,
$$
and a surjective $\Mm$-morphism
$$
(X\times _SY )\times _{\Gg}
(X\times _SY )\rightarrow
(X\times _{\Ff}\Gg )\times _{\Gg}
(X\times _{\Ff}\Gg ).
$$
Composing these three morphisms we obtain a surjective $\Mm$-morphism
$$
T\rightarrow X\times _{\Ff}X.
$$
This proves that $\Ff$ is $P4(\Mm )$.

Suppose now that $\Ff$ is $P4(\Mm )$.  Let
$$
X\rightarrow \Ff , \;\;\; R\rightarrow X\times _{\Ff} X
$$
be the presentation given by the property $P4(\Mm )$.  We may choose $X$ in
such a way
that there exists a lifting $X\rightarrow \Gg$ (the freedom to
replace $X$ by an etale cover comes from Property $M3$ and Lemma
\ref{I.1.c}).  This gives an isomorphism $X\times _{\Ff}\Gg\cong X\times _S\Hh$.
Let  $$
Y\rightarrow \Hh , \;\;\; W\rightarrow Y\times _{\Hh}Y
$$
be the presentation given by the property $P4(\Mm )$ of $\Hh$.  We obtain
surjective $\Mm$-morphisms
$$
X\times _SY\rightarrow X\times _S\Hh
$$
and (defining $U:= X\times _SW$)
$$
U:=X\times _SW \rightarrow (X\times _SY)\times _{X\times _S\Hh }(X\times _SY).
$$
Put $Z:= X\times _SY$.  Then we have  surjections in $\Mm$
$$
Z\rightarrow X\times _{\Ff}\Gg \rightarrow \Gg
$$
(giving the first part of property $P4(\Mm )$), and
$$
U\rightarrow Z\times _{X\times _{\Ff}\Gg}Z.
$$
Now,
$$
(X\times _{\Ff}\Gg )\times _{\Gg} (X\times _{\Ff}\Gg )=
$$
$$
X\times _{\Ff}(X\times _{\Ff}\Gg )=(X\times _{\Ff}X)\times _{\Ff}\Gg ,
$$
and we have an $\Mm$-surjection
$$
R\times _{\Ff}\Gg \rightarrow (X\times _{\Ff}X)\times _{\Ff} \Gg .
$$
Since $R\rightarrow \Ff$ lifts to $R\rightarrow \Gg$ we have $R\times
_{\Ff}\Gg=R\times _S\Hh$ and letting $V\rightarrow R\times _SY$ be an etale
surjection (needed for a certain step below), we obtain $\Mm$-surjections
$$
V\rightarrow R\times _SY \rightarrow R\times _S\Hh \rightarrow (X\times
_{\Ff}\Gg
)\times _{\Gg} (X\times _{\Ff}\Gg ).
$$
On the other hand,
$$
Z\times _{\Gg} Z= Z\times _{X\times _{\Ff}\Gg }((X\times _{\Ff}\Gg )
\times _{\Gg}
(X\times _{\Ff}\Gg ))\times _{X\times _{\Ff}\Gg }Z
$$
so we obtain a  surjection in $\Mm$
$$
Z\times _{X\times _{\Ff}\Gg }V\times _{X\times _{\Ff}\Gg }Z\rightarrow Z\times
_{\Gg}Z.
$$
We can assume (by choosing $V$ appropriately) that the morphism
$$
V\rightarrow (X\times _{\Ff}\Gg )\times
_{\Ff} (X\times _{\Ff}\Gg )
$$
lifts to a morphism
$$
V\rightarrow Z\times_{\Ff} Z.
$$
We then have an $\Mm$-surjection
$$
U\times _ZV\times _ZU
\rightarrow
(Z\times _{X\times _{\Ff}\Gg }Z)\times _ZV\times _Z
(Z\times _{X\times _{\Ff}\Gg }Z)
$$
(where the two maps from $V$ to $Z$ used in the fiber product are the two
projections composed with $V\rightarrow Z\times _{\Ff}Z$). Note that the right
hand side is equal to
$$
Z\times _{X\times _{\Ff}\Gg }V\times _{X\times _{\Ff}\Gg }Z,
$$
which admits, as we have seen above, an $\Mm$-surjection to $Z\times _{\Gg}Z$.
Since $U\times _ZV\times _ZU$ is a scheme of finite type over $S$, this
completes the verification of the property $P4(\Mm )$ for $\Gg$.

We have now shown the equivalence of the conditions $P4(\Mm )$ for $\Ff$ and
$\Gg$.
By the lemma, the morphism $\Gg \rightarrow \Ff$ is in $\Mm$.
By Property $M1$, the structural morphism $\Ff \rightarrow S$ is in $\Mm$
if and only if the structural morphism $\Gg \rightarrow S$ is. Given the
equivalence of the conditions $P4(\Mm )$, this gives equivalence of the
conditions $P5(\Mm )$.
\eop

Finally we give a lemma which allows us some flexibility in specifying
resolutions.

\begin{lemma}
\mylabel{I.1.j}
Suppose that $F$ is a sheaf on $S$ with surjective $\Mm$-morphisms
$X\rightarrow F$ and $R \rightarrow X\times _F X$ such that $X$ and $R$ are
$P4(\Mm )$.
Then $F$ is $P4(\Mm )$.
\end{lemma}
{\em Proof:}
Let $X'\rightarrow X$ and $Q\rightarrow X'\times _XX'$, and $R'\rightarrow R$
be the $\Mm$-surjections given by the hypotheses.
We obtain a surjection $X'\rightarrow F$ in $\Mm$.  On the other hand,
$R'\rightarrow X\times _FX$ is in $\Mm$ and surjective, so
$$
X'\times _XR'\times _XX' = R'\times _{X\times _FX}(X'\times _FX')\rightarrow
X'\times _FX'
$$
is an $\Mm$-surjection.  But the left side is equal to
$$
(X'\times _XX')\times _{X'}R' \times _{X'}(X'\times _XX')
$$
if we  choose (as we may assume is possible) a lifting $R'\rightarrow X'\times
_FX'$ over $X\times _FX$.  There is thus a surjection in $\Mm$
$$
Q\times _{X'}R'\times _{X'}Q\rightarrow X'\times _XR'\times _XX'.
$$
Composing we get the required
$$
Q\times _{X'}R'\times _{X'}Q\rightarrow X'\times _FX'.
$$
\eop

\subnumero{Stability of the condition $P5(\Mm )$}

In the following corollary and theorem we will make use of a supplementary
condition on the class $\Mm$:
\newline
{\bf M4}\,\, If  $f: \Ff \rightarrow \Gg$
is a surjective morphism of sheaves of groups, then $f$ is in $\Mm$.

\begin{corollary}
\mylabel{I.z}
Suppose $\Mm$ satisfies condition M4 in addition to the conditions M1-3.
If $\Gg$ is a $P4(\Mm )$ group sheaf then it is also $P5(\Mm )$.
\end{corollary}
Indeed, M4 applied with $\Gg = \{ 1\}=S$ gives that the
structural morphism $\Ff \rightarrow S$ for any sheaf of groups, is in $\Mm$.
\eop

\begin{theorem}
\mylabel{I.1.e}
Suppose $\Mm$ satisfies condition M4 in addition to the conditions M1-3.
Then if
$$
1\rightarrow \Ff \rightarrow \Ee \rightarrow \Gg \rightarrow 1
$$
is an extension of group sheaves and if any two of the elements are
$P5(\Mm )$, the
third one is too.
\end{theorem}
{\em Proof:} Suppose that $\Ff$ is $P5(\Mm )$.  Then $\Ee$ is $P5(\Mm )$ if and
only if $\Gg$
is $P5(\Mm )$ (by applying the previous theorem in view of the fact that $\Ff$
acts
freely on $\Ee$ with quotient $\Gg$).  The remaining case is if $\Ee$ and $\Gg$
are $P5(\Mm )$.  Then by Lemma \ref{I.1.a}, the kernel $\Ff$ (which is an
equalizer of
two maps $\Ee \rightarrow \Gg$) is $P4(\Mm )$.  By the above corollary,
$\Ff$ is
$P5(\Mm )$.
\eop

\numero{Lifting properties and verticality}

We now fill in what class of morphisms $\Mm$ we would like to use in the theory
sketched above.
We could, of course, take $\Mm = \Xx$ to be the full set of morphisms of $\Xx$.
This might well be a reasonable choice, but I don't see how to get a good
infinitesimal theory in characteristic zero out of this choice.  We could
also try,
for example, to take $\Mm$ as the class of flat (or maybe smooth) morphism
s.  But
then any non-flat group scheme over $S$ would be a counterexample to
property M4,
and as we have seen this property is essential to be able to specify a class of
presentable groups closed under kernels, cokernels and extensions. Thus we have
to work a little harder to find an appropriate class of morphisms.

We say that a morphism of sheaves $a:\Ff \rightarrow \Gg$, is {\em vertical} (or
{\em $S$-vertical}, if the base needs to be specified), if it satisfies the
following lifting properties for all $n\geq 1$:

Suppose $Y$
is a scheme with  $n$ closed subschemes $Y_i\subset Y$, with retractions
$r_i:Y\rightarrow Y_i$---commuting pairwise ($r_ir_j=r_jr_i$)---such that for
$j\leq i$, $r_i$ retracts $Y_i$ to $Y_j\cap Y_i$. Suppose given a morphism
$Y\rightarrow \Gg$, and liftings $\lambda _i:Y_i\rightarrow \Ff$ such
that $\lambda
_i|_{Y_i\cap Y_j}= \lambda _j|_{Y_i\cap Y_j}$. Then for any $P\in Y$
lying on at
least one of the $Y_i$ there exists an etale neighborhood $P\in Y' \rightarrow
Y$ and  a lifting $\lambda : Y' \rightarrow \Ff$  which agrees with the given
liftings $\lambda _i|_{Y_i\times _YY'}$ on $Y_i\times _YY'$.

For future reference we call this lifting property $Lift _n(Y; Y_i)$.

\begin{lemma}
\mylabel{I.u.1}
Suppose $f:\Ff \rightarrow \Gg$ is a morphism of sheaves which are P2.
Then $f$ is vertical if and only if $Lift _n(Y; Y_i)$ holds for all systems
$(Y; Y_i)$ with $Y$ (and hence $Y_i$) artinian.
\end{lemma}
{\em Proof:}  Suppose given a system $(Y, Y_i)$ which is not artinian. Choose a
point $y_0$ (in one of the $Y_i$) and try to find a lifting in an etale
neighborhood of $y_0$.  We can find liftings on $Y^{(n)}$ (the infinitesimal
neighborhoods of $y_0$) by hypothesis.  Using the P2 property of $\Ff$ and an
argument similar to that of Theorem \ref{I.t.2}, we can choose a compatible
sequence of liftings.  Since $\Ff$ is B2 we obtain a lifting over the spectrum
of the complete local ring, then by Artin approximation (using B1) we obtain a
lifting on an etale neighborhood of $y_0$.
\eop

\begin{theorem}
\mylabel{I.u}
We have the following statements: \newline
1. \,\, If $\Ff \rightarrow \Gg$ is vertical and if $\Hh \rightarrow
\Gg$ is any morphism of sheaves, then $\Ff \times _{\Gg}\Hh \rightarrow \Hh$ is
vertical.
\newline
2. \,\, If $a:\Ff \rightarrow \Gg$ is a morphism of sheaves such that for any
$S$-scheme $Y$ and morphism $Y\rightarrow \Gg$, we have that $\Ff \times
_{\Gg}Y\rightarrow Y$ is vertical, then $a$ is vertical;
\newline
3. \,\, If $a:\Ff \rightarrow \Gg$ and $b:\Gg \rightarrow \Hh$ are two morphisms
which are vertical, then $ba$ is vertical (also the identity is vertical); and
\newline
4. \,\,    If $a:\Ff \rightarrow \Gg$ and $b:\Gg \rightarrow \Hh$ are two
morphisms such that $a$ and $ba$ are vertical, and $a$ is surjective, then
$b$ is vertical.
\newline
5.\,\, The etale
surjections between schemes are vertical. \newline
6.\,\, Any injective morphism $\Ff \hookrightarrow S$ is vertical.
\newline
7.\,\, If  $f: \Ff \rightarrow \Gg$
is a surjective morphism of sheaves of groups, then $f$ is vertical.
\end{theorem}
{\em Proof:}
The lifting property concerns only maps from
schemes to $\Gg$, so it obviously satisfies parts 1 and 2.  For part 3, just
lift two times successively (for the identity the lifting property is
tautological).  For part 4, the proof is by induction on $n$.  Keep the
notations
$a$, $b$, $\Ff$, $\Gg$ and $\Hh$ of part 4. Suppose $n=1$.  Then  we
just have to
note that if we have a lifting $Y_1 \rightarrow \Gg$ for $b$, then since $a$ is
surjective, we can lift further to $Y_1\rightarrow \Ff$ (locally in the etale
topology).  The lifting for $ba$ gives $Y\rightarrow \Ff$ and we just
project back
to $\Gg$ to get the lifting for $\Gg$. This gives the case $n=1$. We may
assume that the present lemma is known when there are strictly fewer than  $n$
subschemes.  Suppose we have liftings $\lambda _i:
Y_i\rightarrow \Gg$;  in order to get a lifting $\lambda$, and using the lifting
property for the morphism $ba$, it suffices to choose liftings $\mu _i :
Y_i\rightarrow \Ff$ with $\mu _i|_{Y_i\cap Y_j}= \mu _j|_{Y_i\cap Y_j}$. We can
do this by induction.  Suppose we have chosen $\mu _1,\ldots , \mu _{k-1}$.
Since $k-1<n$, we know the lemma when there are $k-1$ subschemes; apply the
lifting property for the morphism $a$ with respect to the morphism
$Y_k\rightarrow \Gg$, with respect to the subschemes $Y_k\cap Y_i$,
$i=1,\ldots ,
k-1$, and with respect to the liftings $\mu _j|_{Y_k\cap Y_j}$. We obtain a
lifting $\mu _k:Y_k\rightarrow \Ff$ such that
$\mu _k|_{Y_k\cap Y_j}=\mu _j|_{Y_k\cap Y_j}$. By induction now we obtain all
of the liftings $\mu _1,\ldots , \mu _n$.  The lifting property for $ba$ gives
a lifting $\mu$ and we can set $\lambda := a\mu$.  This completes the
verification of part 4.

For the etale surjections (part 5), use the previous lemma. Suppose $i:\Ff
\rightarrow S$ is injective (part 6).  To verify the lifting property for
$Y\rightarrow S$ we just have to verify that this morphism factors through
$Y\rightarrow \Ff$.  For this, use the facts that $Y$ retracts onto $Y_1$ (over
$S$) and that the morphism $Y_1\rightarrow S$ factors through $Y_1\rightarrow
\Ff$.

Finally we verify $Lift _n(Y, Y_i)$ for the morphism $f:\Ff \rightarrow
\Gg$ in part 7.
Let $r_i: Y\rightarrow Y_i$ denote the retractions. Suppose given $\mu :
Y\rightarrow \Gg$ and $\lambda _i : Y_i \rightarrow \Ff$ satisfying the
necessary compatibility conditions. Since $f$ is surjective, we may suppose
that there is a lifting $\sigma : Y\rightarrow \Ff$ of $\mu$ (by
restricting to an etale neighborhood in $Y$).  We construct inductively $\phi
_i : Y\rightarrow \Ff$ lifting $\mu$, with $\phi _i |_{Y_j}=\lambda _j$ for
$j\leq i$.  Denote the multiplication operations in $\Ff$ or $\Gg$ by $\cdot$.
Let  $$
h_1:= (\lambda _1\cdot (\sigma |_{Y_1})^{-1})\circ r_1: Y\rightarrow \ker (f).
$$
Put $\phi _1 := h_1\cdot \sigma $.
Then $\phi _1$ restricts to $\lambda _1$ on $Y_1$, and lifts $\mu$.  Suppose we
have chosen $\phi _i$.  Let
$$
h_{i+1}:= (\lambda _{i+1}\cdot (\phi _i |_{Y_{i+1}})^{-1})\circ r_{i+1}
:Y\rightarrow \ker (f),
$$
and put $\phi _{i+1}:= \phi _i \cdot h_{i+1}$.  This lifts $\mu$ because
$h_{i+1}$ is a section of $\ker (f)$.  For $j\leq i$, $r_{i+1}$ maps
$Y_j$ to $Y_j\cap Y_{i+1}$, and there $\lambda _{i+1}=\lambda _j$ agrees with
$\phi _i$ so $h_{i+1}|_{Y_j}=1$.  We don't destroy the required property for
$j\leq i$.  On the other hand, we gain the required property for $j=i+1$, by
construction.  This completes the inductive step to construct $\phi _i$.
Finally, the $\phi _n$ is the lifting required for property $Lift _n(Y,Y_i)$.
This completes the proof of part 7.
\eop

From the above results, the class $\Mm$ of vertical morphisms satisfies the
axioms M1, M2, M3, {\em and} M4 of the previous section.
This is the principal class $\Mm$ to which we will refer, in view of which we
drop $\Mm$ from the notation when $\Mm$ is the class of vertical morphisms.
Thus
the conditions  $P4$ and $P5$ refer respectively to
$P4(\Mm )$ and $P5(\Mm )$ with $\Mm$ the class of vertical morphisms.

In particular we obtain the results \ref{I.z.1}, \ref{I.1.a}, \ref{I.1.a.1},
\ref{I.1.b}, \ref{I.1.c}, \ref{I.1.d}, \ref{I.z}, and \ref{I.1.e} for
the properties P4 and P5.

We have some further results about $P4$ and $P5$.

\begin{lemma}
\mylabel{I.x}
Suppose that $\Ff$ is $P4$.
In the situation of the lifting property $Lift_n(Y; Y_i)$ for the morphism
$X\rightarrow \Ff$ given by property $P4$, suppose that $Y$ is the
scheme-theoretic union of $Y_1,\ldots , Y_n$.  Then the lifting is unique.
\end{lemma}
{\em Proof:}
In effect, for morphisms $Y\rightarrow X$ with $X$ a scheme, if $Y$ is the
scheme theoretic union of the $Y_i$ then the morphism is determined by its
restrictions to the $Y_i$.
\eop

\begin{proposition}
\mylabel{aaa}
The property of being vertical is stable under base change
of $S$: suppose $p:S'\rightarrow S$ is a morphism of schemes.  If $f:\Ff
\rightarrow \Gg$ is vertical then  $p^{\ast}(f):p^{\ast}(\Ff)\rightarrow
p^{\ast}(\Gg )$ is vertical.  Furthermore if $\Hh \rightarrow \Kk$ is an
$S'$-vertical morphism of sheaves on $\Xx /S'$ then the restriction down to $S$,
$$
Res _{S'/S}(\Hh )\rightarrow Res _{S'/S}(\Kk )
$$
is $S$-vertical.
\end{proposition}
{\em Proof:}
This follows from the form of the lifting
properties.
\eop

{\em Remark:}  We often ignore the notation of ``restriction down'', then the
first part of the proposition states that if $\Ff \rightarrow \Gg \rightarrow
S$ with the first morphism being $S$-vertical, then
$\Ff \times _SS'\rightarrow \Gg \times _SS'$ is $S'$-vertical.
The last part of
the proposition states that if $\Hh \rightarrow \Kk \rightarrow S'$ with the
first morphism being $S'$-vertical, then it is also $S$-vertical.

\begin{corollary}
\mylabel{I.1.j.1}
Suppose $\Ff$ is a sheaf over $S$, and suppose $S'\rightarrow S$ is a
surjective etale morphism such that $\Ff |_{S'}$ is $P4$ over $S'$.
Then $\Ff$
is $P4$ over $S$.
\end{corollary}
{\em Proof:}
If $(Y,Y_i)$ is a system for the lifting property over $S$, then their
pullbacks $(Y',Y_i')$ form such a system over $S'$.  If a morphism of sheaves
over $S'$ satisfies the lifting property, then we can lift for the system
$(Y', Y'_i)$.  This gives a lifting over $Y'$ for the system $(Y,Y_i)$, that
is a lifting etale locally, thus satisfying the lifting property over $S$.
Thus a morphism which is $S'$-vertical is also $S$-vertical.  It follows that
$\Ff |_{S'}$ is $P4(\Mm )$ over $S$.  Now
$$
(\Ff |_{S'})\times _{\Ff} (\Ff |_{S'})= \Ff |_{S'}\times _{S'} (S'\times
_SS'),
$$
and $S'\times _SS'$ is $P4(\Mm )$ over $S$. Thus
$(\Ff |_{S'})\times _{\Ff} (\Ff
|_{S'})$ is $P4(\Mm )$ over $S$. We can now apply Lemma \ref{I.1.j} with
$X=\Ff |_{S'}$ and $R=(\Ff |_{S'})\times _{\Ff} (\Ff
|_{S'})$.
\eop

\subnumero{Presentable group sheaves}
In view of the nice properties of $P5$ group sheaves, we make the following
change
of notation.  A sheaf of groups $\Gg$ over $\Xx /S$ is a {\em presentable group
sheaf} if it is  $P5$.
Note that we use this terminology only for sheaves of groups.

\begin{corollary}
\mylabel{uvw}
A sheaf of groups which is representable by
a scheme of finite type $G$ over $S$, is presentable.
\end{corollary}
{\em Proof:}
This is because we can
take $X=G$ and $R$ equal to the diagonal $G$ in the definition of property
$P4$; and property $P5$ is then Corollary \ref{I.z}.
\eop

\begin{corollary}
\mylabel{vwx}
The category of presentable group sheaves contains the
category generated by representable group sheaves under the operations of
extensions, kernels, and division by normal subgroups.
\end{corollary}
{\em Proof:}
Theorem \ref{I.1.e}.
\eop

In particular, the category of presentable group sheaves is
much bigger than the category of representable group sheaves.  I  believe that
the category of presentable group sheaves is strictly larger than the category
generated generated by representable group sheaves under the operations of
kernel, cokernel and extension. For example the group sheaves $Aut (V)$ for a
vector sheaf $V$, which are presentable as shown below, are probably not
generated
from representable group sheaves by kernels, extensions and quotients
(although I
don't have a counterexample).  In an intuitive sense, however, the two
categories
are about the same.

The two previous corollaries would also hold for the category of $P5(\Mm )$
group
sheaves for any class $\Mm$ satisfying $M1$ through $M4$.

We now give the main argument where we use the lifting properties and the notion
of verticality, i.e. the special definition of $\Mm$.

\begin{lemma}
If $\Gg$ is a sheaf of groups and $X\rightarrow \Gg$ is a vertical
surjection, with identity section $e:S\rightarrow X$ then (choosing a point
$P$ on $e(S)$) there is a lifting of the multiplication to a map of etale germs
$$
\mu : (X,P)\times _S(X,P) \rightarrow (X,P)
$$
such that $\mu (x,e)=\mu (e,x)=x$.
\end{lemma}
{\em Proof:}
Let $Y=X\times _SX$ and $Y_1 = X\times _Se(S)\cong X$ and $Y_2 = e(S)\times _SX
\cong X$.
We  have retractions $Y\rightarrow Y_1$ and $Y\rightarrow Y_2$ as in the lifting
property.  The multiplication map $\Gg \times _S\Gg \rightarrow \Gg$ composes to
give a map $Y=X\times _SX\rightarrow \Gg$.  The identity gives liftings
$Y_1\rightarrow X$ and $Y_2\rightarrow X$ agreeing on $Y_1\cap Y_2 = e(S)\times
_S e(S)$.  By the definition of verticality of the morphism $X\rightarrow \Gg$,
there is an etale neighborhood $P\in Y'\rightarrow Y$ and a lifting to a map
$Y'\rightarrow X$ agreeing with our given lifts on $Y'_1$ and $Y'_2$.
This gives
the desired map (note that when we have written the product of two etale germs,
this means the germ of the product rather than the product of the two spectra of
henselian local rings).
\eop

We use this result in the following way.
A map $\mu : X\times X\rightarrow X$ (defined on germs at a point $P$) such that
$\mu (e,x)=\mu (x,e) = x$, gives rise to an exponential map
$T(X)_e^{\wedge}\rightarrow X$ where $T(X)_e$ is the tangent vector scheme (see
\S\S 5-8 below) to $X$ along the identity section $e$ and $T(X)_e^{\wedge}$
denotes the formal completion at the zero section.  To define this exponential
map note that the multiplication takes tangent vectors at $e$ to tangent vector
fields on $X$ which we can then exponentiate in the classical way.  The formal
exponential map is an isomorphism between $T(X)_e^{\wedge}$ and the completion
$X^{\wedge}$ along $e$. This is a fairly strong condition on $X$ which we will
exploit below, notably to get $Lie (\Gg )$ and to develop a theory of
connectedness.

In particular this technique allows us to prove directly (in \S 6 below) that
when $k$ is a field of characteristic zero, presentable group sheaves over $Spec
(k)$ are just algebraic Lie groups over $k$.

It is possible that in characteristic $p$ there would be an appropriate
notion of
verticality taking into account divided powers, which would have the same effect
of enabling a good infinitesimal theory. This is why we have left the class
$\Mm$
as an indeterminate in the first part of our discussion above.

\subnumero{The conditions $P3$ and $P3\frac{1}{2}$}

We now add the following two conditions, which will be used as conditions on
$\pi _0$ in the last section (in contrast to the condition $P5$ which is to be
used on $\pi _1$ and even $\pi _i$, $i\geq 2$).  These conditions depend on a
functorial choice of class $\Mm (Y)$ of morphisms of sheaves over $Y$ for each
$Y\in \Xx$.  We will leave to the reader the (easy) job of stating these
properties in this generality, and instead we will state them directly when $\Mm
(Y)$ is taken as the class of $Y$-vertical morphisms. Note that the
properties we
are about to state are {\em absolute} properties of sheaves on $\Xx$ rather than
relative properties of sheaves over some base $S$.

\noindent
{\bf P3.} \,\, A sheaf $\Ff$ on $\Xx$ is $P3$ if there is a surjection
$X\rightarrow \Ff$ from a scheme $X$ of finite type over $Spec (k)$, and
if there
is a surjection $\varphi : R\rightarrow X\times _{\Ff}X$ from a scheme $R$ of
finite type over $Spec (k)$ such that $\varphi$ is an $X\times X$-vertical
morphism.

\noindent
${\bf P3\frac{1}{2}.}$ \,\, A sheaf $\Ff$ on $\Xx$ is $P3\frac{1}{2}$ if
there is a surjection
$X\rightarrow \Ff$ from a scheme $X$ of finite type over $Spec (k)$, and
if there
is a surjection $\varphi : R\rightarrow X\times _{\Ff}X$ from a scheme $R$ of
finite type over $Spec (k)$ such that $\varphi$ is an $X$-vertical
morphism, where the map to $X$ is the first projection of $X\times _{\Ff}X$.

{\em Remark:} These properties seem almost identical.  The first was refered to
in \cite{kobe} (already as property $P3$).  However it will turn out that the
second version (which I hadn't yet thought of at the time of writing
\cite{kobe})
seems more useful---cf \S 10 below.  The author apologizes for this complication
of the notation!

{\em Remark:}
$$
P5 \Rightarrow P4\Rightarrow P3\frac{1}{2}\Rightarrow P3 \Rightarrow P2
\Rightarrow P1.
$$

These properties will not come into our study of group sheaves over a base $S$.
Rather, they come  in as conditions on $\pi _0$ of $n$-stacks on $\Xx$, in our
brief discussion at the end of the paper.  In fact we could have put off stating
these properties until \S 10, but the reader had probably been wondering for
some time already why we are skipping number $3$ in our list of properties.

We quickly give the analogues, for  $P3\frac{1}{2}$, of some of the basic
facts about our other properties. We leave to the reader the task of elicudating
the corresponding properties for $P3$.

\begin{lemma}
\mylabel{P3a}
Suppose $\Gg$ is $P3\frac{1}{2}$ and suppose $X$ is a scheme of finite type with
a morphism $X\rightarrow \Gg$.  Then there is a surjection from a scheme of
finite type  $R\rightarrow X\times _{\Gg} X$ which is vertical with respect to
the first factor $X$.
\end{lemma}
{\em Proof:}
Let $Y\rightarrow \Gg$ and $W\rightarrow Y\times _{\Gg}Y$ be the surjections
with the second one being vertical with respect to the first factor $Y$.
There is an etale covering $X' \rightarrow X$ and a lifting of our morphism to
$X'\rightarrow Y$.  Then
$$
X'\times _{\Gg}X'=(X'\times X') \times _{Y\times Y} (Y\times _{\Gg}Y)
$$
so
$$
R:=(X' \times X' )\times _{Y\times Y} W \rightarrow X'\times _{\Gg}X'
$$
is surjective.  It is vertical with respect to the first factor $Y$ and hence
vertical with respect to the first factor $X'$. Since $X'\rightarrow X$
is etale,
this morphism is also vertical with respect to $X$ (via the first factor).
The surjection
$$
X'\times _{\Gg}X' \rightarrow X\times _{\Gg}X
$$
is the pullback of the etale morphism $X'\times X'\rightarrow X\times X$ so it
is also vertical with respect to the first factor $X$.  Composing we obtain
$$
R\rightarrow X\times _{\Gg}X
$$
vertical with respect to the first factor.
\eop

\begin{corollary}
\mylabel{P3b}
Suppose $\Gg$ is $P3\frac{1}{2}$ and suppose $\Ff \subset \Gg$. If $\Ff$ is $P1$
then it is $P3\frac{1}{2}$.
\end{corollary}
{\em Proof:}
Let $X\rightarrow \Ff$ be a surjection from a scheme of finite type $Y$.
From the above lemma we get a surjection $R\rightarrow X\times _{\Gg}X$
which is vertical with respect to the first factor, but since $\Ff\rightarrow
\Gg$ is injective $X\times _{\Gg} X=X\times _{\Ff}X$ and we're done.
\eop

\begin{corollary}
\mylabel{P3c}
Suppose $\Gg$ is $P3\frac{1}{2}$ and $\Hh$ is $P2$, then the equalizer
$\Ff$ of any two morphisms $f,g: \Gg \rightarrow \Hh$ is again
$P3\frac{1}{2}$.
\end{corollary}
{\em Proof:}
By Lemma \ref{I.1.a} with ${\cal M}$ being the class of all morphisms,
we obtain that $\Ff$ is $P2$ and in particular $P1$. Since it is a subsheaf of
$\Gg$, the previous corollary applies to show that $\Ff$ is
$P3\frac{1}{2}$.
\eop

\begin{corollary}
\mylabel{P3d}
Suppose $\Ff \rightarrow \Hh$ and $\Gg \rightarrow \Hh$ are two morphisms such
that $\Hh$ is $P2$ and $\Ff$ and $\Gg$ are $P3\frac{1}{2}$. Then the fiber
product $\Ff \times _{\Hh} \Gg$ is $P3\frac{1}{2}$.
\end{corollary}
{\em Proof:}
The fiber product is the equalizer of two morphisms $\Ff \times \Gg
\rightarrow \Hh$.  Note that the product of two $P3\frac{1}{2}$ sheaves is
again $P3\frac{1}{2}$---this comes from the general statement that if $\Aa
\rightarrow \Bb$ is $S$-vertical and  if $\Aa '\rightarrow
\Bb'$ is $S'$-vertical then $\Aa \times \Aa'\rightarrow \Bb \times \Bb '$
is $S\times S'$-vertical (a direct consequence of the form of the lifting
properties).
\eop

Finally we have the analogue of one half of Theorem \ref{I.1.d}. I didn't quite
see how to do the other half.

\begin{proposition}
\mylabel{P3e}
Suppose $S$ is a scheme of finite type, and suppose $\Hh$ is a group sheaf over
$S$ which is $P5$. Suppose that $\Gg\rightarrow S$ is a sheaf and that $\Hh$
acts freely on $\Gg$ over $S$, with quotient $\Ff = \Gg /\Hh$. If $\Ff$ is
$P3\frac{1}{2}$ then $\Gg$ is $P3\frac{1}{2}$ (here $\Ff$ and $\Gg$ are being
considered as the restrictions down to $Spec (k)$ of the corresponding
sheaves over $S$).
\end{proposition}
{\em Proof:}
We follow the proof of the second half of Theorem \ref{I.1.d}.
Let
$$
X\rightarrow \Ff , \;\;\; R\rightarrow X\times _{\Ff} X
$$
be the presentation given by the property $P3\frac{1}{2}$.  We may choose
$X$ in
such a way
that there exists a lifting $X\rightarrow \Gg$,
giving an isomorphism
$X\times _{\Ff}\Gg\cong X\times _S\Hh$. Let
$$
Y\rightarrow \Hh , \;\;\; W\rightarrow Y\times _{\Hh}Y
$$
be the presentation given by the property $P4(\Mm )$ of $\Hh$.  We obtain
surjective $S$-vertical morphisms
$$
X\times _SY\rightarrow X\times _S\Hh
$$
and (defining $U:= X\times _SW$)
$$
U:=X\times _SW \rightarrow (X\times _SY)\times _{X\times _S\Hh }(X\times _SY).
$$
Put $Z:= X\times _SY$.  Then we have a surjection
$$
Z\rightarrow X\times _{\Ff}\Gg \rightarrow \Gg
$$
and an $S$-vertical surjection
$$
U\rightarrow Z\times _{X\times _{\Ff}\Gg}Z.
$$
Now,
$$
(X\times _{\Ff}\Gg )\times _{\Gg} (X\times _{\Ff}\Gg )=
$$
$$
X\times _{\Ff}(X\times _{\Ff}\Gg )=(X\times _{\Ff}X)\times _{\Ff}\Gg ,
$$
and we have a surjection vertical with respect to the first factor $X$,
$$
R\times _{\Ff}\Gg \rightarrow (X\times _{\Ff}X)\times _{\Ff} \Gg .
$$
Since $R\rightarrow \Ff$ lifts to $R\rightarrow \Gg$ we have $R\times
_{\Ff}\Gg=R\times _S\Hh$ and letting $V\rightarrow R\times _SY$ be an etale
surjection, we obtain surjections
$$
V\rightarrow R\times _SY \rightarrow R\times _S\Hh \rightarrow
(X\times _{\Ff}\Gg
)\times _{\Gg} (X\times _{\Ff}\Gg ).
$$
The first is etale, the second is $S$-vertical, and the third is $X$-vertical
for the first factor, so the composition is $X$-vertical.
As before
$$
Z\times _{\Gg} Z= Z\times _{X\times _{\Ff}\Gg }((X\times _{\Ff}\Gg )
\times _{\Gg}
(X\times _{\Ff}\Gg ))\times _{X\times _{\Ff}\Gg }Z
$$
so we obtain an $X$-vertical  surjection
$$
Z\times _{X\times _{\Ff}\Gg }V\times _{X\times _{\Ff}\Gg }Z\rightarrow Z\times
_{\Gg}Z.
$$
We can assume by choosing $V$ appropriately that the morphism
$$
V\rightarrow (X\times _{\Ff}\Gg )\times
_{\Ff} (X\times _{\Ff}\Gg )
$$
lifts to a morphism
$$
V\rightarrow Z\times_{\Ff} Z.
$$
We then have an $S$-vertical surjection
$$
U\times _ZV\times _ZU
\rightarrow
(Z\times _{X\times _{\Ff}\Gg }Z)\times _ZV\times _Z
(Z\times _{X\times _{\Ff}\Gg }Z).
$$
The right
hand side is equal to
$$
Z\times _{X\times _{\Ff}\Gg }V\times _{X\times _{\Ff}\Gg }Z,
$$
which admits, as we have seen above, an $X$-vertical surjection to $Z\times
_{\Gg}Z$. By composing we obtain an $X$-vertical, and hence $Z$-vertical
surjection
$$
U\times _ZV\times _ZU\rightarrow Z\times _{\Gg}Z.
$$
This completes the proof.
\eop

\numero{Functoriality}

Suppose $F$ is a sheaf over $S$, and suppose $\pi : S'\rightarrow S$ is a
morphism.  We denote by $\pi ^{\ast}(F)$ the restriction $F|_{\Xx /S'}$, which
is the sheaf associated to the presheaf $Y\rightarrow S' \mapsto F(Y\rightarrow
S)$. If $F$ is representable then $\pi ^{\ast}F$ is also representable by the
fiber product $F\times _SS'$.  In general, we allow ourselves to use the
notations $\pi ^{\ast}F$, $F\times _SS'$ and $F|_{S'}$ interchangeably.

We have  defined, for a sheaf $G$ on $S'$, the {\em restriction down $Res
_{S'/S}(G )$}.

Suppose $G$ is a
sheaf on $S'$.  We defined the direct image by
$$
\pi _{\ast}G(Y\rightarrow S):= G(Y\times _SS' \rightarrow S').
$$
The morphism $F(Y\rightarrow S) \rightarrow F(Y\times _SS' \rightarrow S)$
gives a natural morphism
$$
F\rightarrow \pi _{\ast} \pi ^{\ast}(F),
$$
and the morphism $G(Y\times _S S'\rightarrow S) \rightarrow G (Y\rightarrow
S')$ (coming from the graph morphism $Y\rightarrow Y\times _S S'$) gives a
natural morphism
$$
\pi ^{\ast}\pi _{\ast} (G)\rightarrow G.
$$
These functors are adjoints and the above are the adjunction morphisms. More
precisely, we have a natural isomorphism
$$
Hom (F, \pi _{\ast}G)\cong Hom (\pi ^{\ast}F,G).
$$
This may be verified directly.

{\em Remark:} If $f:A\rightarrow B$ is a vertical morphism over
$S'$ then  $\pi _{\ast}(f):\pi _{\ast}A\rightarrow \pi _{\ast}B$ is vertical
over $S$.  To see this, note that if $Y, Y^{(n)}$ is a collection of
$S$-schemes with retractions etc. as in the definition of verticality, then
$Y\times _S S', Y^{(n)} \times _SS'$ is a collection with retractions over
$S'$.  The verticality of $\pi _{\ast}(f)$ for the case of  $Y, Y^{(n)}$
follows from the verticality of $f$ for the case of
$Y\times _S S', Y^{(n)} \times _SS'$.

{\em Remark:}  Direct and inverse images are compatible with fiber products.
For inverse images this is easy.  For direct images, suppose we have morphisms
$A\rightarrow C$ and $B\rightarrow C$ on $S'$.  We obtain morphisms $A\times
_CB\rightarrow A$ and $A\times _CB$ satisfying a universal property.  These
give morphisms
$$
\pi _{\ast}(A\times _CB)\rightarrow \pi _{\ast} A \;\; (resp. \;\; \pi _{\ast}B
\, ).
$$
We show the universal property: suppose
$$
(u,v)\in  (\pi _{\ast}A\times _{\pi _{\ast}C}\pi _{\ast}B  )(Y),
$$
that is $u\in A(Y\times _SS') $ and $v\in B(Y\times _SS')$ with the same image
in $C(Y\times _SS')$.  We obtain a unique element of $(A\times _CB)(Y\times
_SS')$ mapping to $(u,v)$.  This gives the claim.

\begin{lemma}
\mylabel{I.1.g.2}
if $\Ff$ is a coherent sheaf on $S'$ and $\pi :S'\rightarrow S$ is a finite
morphism then $\pi _{\ast}(\Ff )$ is a coherent sheaf on $S$.
\end{lemma}
{\em Proof:}
We may assume $S$ and $S'$ affine, so that $S=Spec (A)$ and $S'=Spec (A')$ with
$A'$ a finite $A$-algebra. The coherent sheaf $\Ff$ corresponds to an
$A'$-module $M$.  This implies that
$$
\pi _{\ast}(\Ff ) (Spec (B)\rightarrow Spec (A)) = \Ff ( Spec (B\otimes _AA'))
$$
$$
= M\otimes _{A'}(B\otimes _AA') = M\otimes _AB.
$$
This formula means that $\pi _{\ast}(\Ff )$ corresponds to the same
module $M$ considered as an $A$-module; in particular it is coherent.
\eop

\begin{lemma}
\mylabel{I.1.h}
If $F$ is $P4$ (resp. $P5$) on $S$ then $\pi ^{\ast}F$ is $P4$ (resp. $P5$)
on $S'$.
\end{lemma}
{\em Proof:}
Note first of all that $S$-verticality of a morphism of sheaves over $S'$
implies $S'$-verticality.
Now if $F$ is $P4$, let $X\rightarrow F$ and $R\rightarrow X\times _FX$ be the
corresponding vertical surjections.  We get $\pi ^{\ast}(X) \rightarrow \pi
^{\ast}(F)$ and
$$
\pi ^{\ast}(R)\rightarrow \pi ^{\ast}(X\times _FX)=
\pi ^{\ast}(X)\times _{\pi ^{\ast}(F)}\pi ^{\ast}(X),
$$
surjective and $S$-vertical (hence $S'$-vertical) morphisms.  Note that $\pi
^{\ast}(X)$ and $\pi ^{\ast}(R)$ are schemes of finite type over $S'$, so we
obtain the proof for P3. For $P5$ note that $\pi ^{\ast}(F)=F\times _SS'$,
so by Theorem \ref{I.u}, $\pi ^{\ast}(F)\rightarrow S'$ is $S$-vertical; hence
it is $S'$-vertical as required.
\eop

\begin{lemma}
\mylabel{I.1.i}
Suppose $\pi :S'\rightarrow S$ is a finite morphism and suppose $G'$ is a $P4$
sheaf on $S'$.  Then $\pi _{\ast}(G')$ is a $P4$ sheaf on $S'$.
\end{lemma}
{\em Proof:}
Let $X'\rightarrow G'$ and $R'\rightarrow X'\times_{G'}X'$ be the surjective
vertical morphisms with $X'$ and $R'$ schemes of finite type over $S'$.
Let $G:= \pi _{\ast}(G')$ and similarly for $X$ and $R$.
By the above remark, we obtain vertical morphisms $X\rightarrow G$ and
$R\rightarrow X\times _GX$.  (Note that $X\times _GX= \pi _{\ast}(X'\times
_{G'}X'$ by above.)

In the case of a finite morphism $\pi : S'\rightarrow S$, note that if
$f:A\rightarrow B$ is surjective over $S'$ then  $\pi _{\ast}(f):\pi
_{\ast}A\rightarrow \pi _{\ast}B$ is surjective.  This is a general property of
sheaves on the etale topology, for which we sketch the proof (an application
of Artin approximation).  If $\eta \in \pi _{\ast}(B)(Y)$, this means $\eta :
Y\times _SS'\rightarrow B$.  For $y'\in Y\times _SS'$ there is an etale
neighborhood $U\rightarrow Y\times _SS'$ and a lifting $U\rightarrow A$.  We
need to find an etale neighborhood $V$ of the image $y\in Y$ and a  lifting
$V\times _SS' \rightarrow U$.  Define a functor $L(V/Y)$ to be the set of
liftings $V\times _SS' \rightarrow U$ over $Y\times _SS'$.  It is B1, and a
lifting exists on $\hat{V}= Spec (\Oo _{Y,y}^{\wedge})$, so by Artin
approximation there is an etale neighborhood $V$ with a lifting.

Applying this to our case, the morphisms $X\rightarrow G$ and
$R\rightarrow X\times _GX$ are surjective.  By Lemma \ref{I.1.j}, it suffices to
prove that $X$ and $R$ are $P4$.  Thus it suffices in general to show:
if $Z$ is a
scheme of finite type over $S'$ then $\pi _{\ast}(Z)$ is $P4$.  We make a
further
reduction: a scheme of finite type can be presented as the kernel of a morphism
${\bf A}^n \rightarrow {\bf A}^m$; the direct image is then the kernel of $\pi
_{\ast}{\bf A}^n \rightarrow  \pi _{\ast}{\bf A}^m$.  The kernel of a morphism
of $P4$ sheaves is again $P4$ (Lemma \ref{I.1.a}) so it suffices to
treat the case
$Z={\bf A}^n$.  But in this case, $Z$ is a coherent sheaf and its direct image
is also a coherent sheaf.  One can see directly that a coherent sheaf $\Ff$ is
$P4$ by using the fact that it has a resolution of the form
$$
\Oo ^a \rightarrow \Oo ^b \rightarrow \Ff \rightarrow 0
$$
(exact even on the big site $\Xx /S$), or by looking at Lemmas
\ref{I.1.g.1} (below) and \ref{I.1.g.2}.  This completes the proof.
\eop

\begin{lemma}
\mylabel{restrictionPreserves?}
Suppose $S'\rightarrow S$ is a morphism of schemes of finite type.  Suppose
$\Ff$ is a sheaf on $\Xx /S'$.  Then $Res _{S'/S}\Ff$ is $P3\frac{1}{2}$  if and
only if $\Ff$ is $P3\frac{1}{2}$.
\end{lemma}
{\em Proof:}
This is only a matter of terminology since, $P3\frac{1}{2}$ being a global
property, the statement that $\Ff$ is $P3\frac{1}{2}$ really means that the
restriction of $\Ff$ down to $Spec (k)$ is $P3\frac{1}{2}$.   It is obviously
equivalent to say this after first restricting down to $S$.
\eop

\numero{Vector sheaves}

With this section we begin the part of our study which requires working over a
ground field $k$ of characteristic zero.  From now on $\Xx$ denotes the big
etale
site of schemes over $Spec (k)$.  Before returning to the definition of
presentability and its infinitesimal study, we make a detour to discuss {\em
vector sheaves}.  These are objects which will be the linearizations of
presentable group sheaves---we are also interested in vector sheaves as
candidates for the $\pi _i (T,t)$ with $t\in T(S)$, for an $n$-stack $T$ on $\Xx
/S$ for $i\geq 2$.

To be slightly more precise,
suppose $S\in \Xx$ is a scheme over  $k$, and let $\Xx /S$ denote the
category of
schemes over $S$.  We will define a notion of {\em vector sheaf} on $\Xx /S$.

This notion is what was called ``U-coherent sheaf'' by Hirschowitz in
\cite{Hirschowitz}.
The particular case which we call ``vector scheme'' below has already been well
known for some time as the ``linear spaces'' of Grauert \cite{Grauert},
appearing
notably in Whitney's tangent cones \cite{Whitney}.

We feel that the terminology ``vector sheaf'' is more
suggestive.  Many of the results below seem to be due to Hirschowitz
\cite{Hirschowitz} (in particular, the observation that duality is involutive)
although some parts of the theory are certainly due to
\cite{Fischer}, \cite{Grauert}, \cite{Whitney}. We
have integrated these results into our treatment for the reader's convenience.
Essentially the only thing new in our treatment is the first lemma (and the
analogous statement about extensions).

Before starting in on the definition, I would like to make one note of caution.
The category of vector sheaves will not satisfy any nice  (ascending or
descending) chain condition.  This is one of the principal differences with
vector spaces or modules over a noetherian ring, and could in the long run pose
a major problem if one wants to consider an ``infinite dimensional'' version of
the theory such as by looking at $ind$- or $pro$- objects.

We have a sheaf of rings $\Oo$ on $\Xx$, defined by $\Oo (X):= \Gamma (X, \Oo
_X)$.  Note that it is represented by the affine line.

\begin{lemma}
\mylabel{I.a}
Suppose $F$ is a sheaf of abelian groups on $\Xx /S$, representable by a scheme
which is affine and of finite type over $S$.  If there exists a structure of
$\Oo$-module for $F$, then this structure is unique.  If $F$ and $G$ are two
such sheaves, and if $a:F\rightarrow G$ is a morphism of
sheaves of groups, then $a$ is a morphism of $\Oo$-modules.
\end{lemma}
{\em Proof:}
The first statement of the lemma follows from the second. For the second
statement, suppose $u\in F_X$.  Consider the element $tu\in F_{X\times {\bf
A}^1}$.  For any positive integer $n$ we have $tu|_{X\times \{ n\}}=u+\ldots +u$
($n$ times).  The same is true for the image $a(u)$.  Therefore
$$
a(tu)|_{X\times \{ n\}} =ta(u)|_{X\times \{ n\}} .
$$
We obtain two morphisms $X\times {\bf A}^1\rightarrow G$ which are equal on
the subschemes $X\times \{ n\}$; this implies that they are equal.
(Here is a proof of this:  we may suppose that $X$ and the base $S$ are
affine, so $X=Spec (A)$ and $G=Spec (B)$ and a morphism $X\times {\bf
A}^1\rightarrow G$  corresponds to a morphism $\phi  : B\rightarrow
A[t]$.
Pick any $b\in B$ and write
$$
\phi  (b)= \sum _{j=1}^m p_{j}t^j;
$$
but the matrix $a_{nj}= n^j$ for $n,j=1,\ldots , m$ is invertible as a matrix
with coefficients in $k$, so there is a matrix $c_{nj}$ with
$$
p_j= \sum _{n=1}^m c_{nj}\phi (b)(n) .
$$
Thus $\phi (b)$ is determined by the values at positive integers $\phi (b)(n)$.)
\eop

A {\em vector scheme over $S$} is a sheaf $V$ of abelian groups on $\Xx /S$
which
is a sheaf of $\Oo$-modules and such that there exists an etale covering $\{
S_{\alpha}\rightarrow S\}$ such that each $V|_{S_{\alpha}}$ is representable
by a scheme $F_{\alpha}$ which is affine of finite type over $S_{\alpha}$.
The
above lemma shows that the category of vector
schemes is a full subcategory of the category of sheaves of abelian groups on
$\Xx$.  In the complex analytic category these objects were called ``linear
spaces'' by Grauert and were studied in \cite{Grauert}, \cite{Fischer}.

The first remark is that, in fact, the locality in the  definition
of vector scheme was
extraneous.  In effect, since the representing schemes $F_{\alpha}$ are unique
up to unique isomorphism, they glue together to give a scheme $F$, affine and
locally of finite type over $S$.

\begin{lemma}
\mylabel{I.b}
Suppose $V$ is a vector scheme on $\Xx /S$, and suppose $S$ is affine. Then
there
is an exact sequence
$$
0\rightarrow V\rightarrow \Oo ^m \rightarrow \Oo ^n
$$
of abelian sheaves on $\Xx$.
\end{lemma}
{\em Proof:}
Write $S=Spec (A)$ and $V=Spec
(B)$.  The action of ${\bf G} _m$ gives a decomposition
$$
B = \bigoplus B^{\lambda}
$$
where $B^{\lambda}$ consists of functions $b$ such that $b(tv)=
t^{\lambda}b(v)$.
The sum is over $\lambda \geq 0$ (integers), since the action extends to an
action of the multiplicative monoid ${\bf A}^1$.  Furthermore, if $b\in B^0$
then $b(tv)=b(v)$ for all $t$ (including $t=0$), in particular $b(v)=b(0)$.
Thus $B^0= A$.  If $b\in B^{\lambda }$ for $\lambda >0$ then $b(0)=b(0\cdot
O)= 0$.  Thus the zero section corresponds to the projection onto $B^0=A$.
The decomposition is compatible with multiplication in $B$. It is also
compatible with the comultiplication $B\rightarrow B\otimes _A B$
corresponding to the addition law on $V$.  The comultiplication is
$$
B^{\lambda } \rightarrow \bigoplus _{\mu + \nu = \lambda} B^{\mu} \otimes _A
B^{\nu},
$$
and furthermore the coefficients $B^{\lambda} \rightarrow B^{\lambda}\otimes
_A B^0= B^{\lambda}$ and  $B^{\lambda} \rightarrow B^0\otimes
_A B^{\lambda}= B^{\lambda}$ are the identity (corresponding to the formula
$v+0=v=0+v$).  On the other hand, the composition $B\rightarrow B\otimes _A B
\rightarrow B$ corresponds to the map $v\mapsto v+v=2v$, which is also scalar
multiplication by $t=2$.  Thus the composition
$$
B^{\lambda } \rightarrow \bigoplus _{\mu + \nu = \lambda} B^{\mu} \otimes _A
B^{\nu}\rightarrow B^{\lambda}
$$
is equal to multiplication by $2^{\lambda}$.  The first and last terms in the
sum give a contribution of $b\mapsto 2b$ (by the observation $v+0=v=0+v$), so
for $\lambda \geq 2$, the composition
$$
B^{\lambda } \rightarrow \bigoplus _{\mu + \nu = \lambda , 0< \mu , \nu <
\lambda } B^{\mu} \otimes _A B^{\nu}\rightarrow B^{\lambda}
$$
is multiplication by $2^{\lambda}-2$, invertible.  Hence every element of
$B^{\lambda}$ is expressed as a sum of products of elements of $B^{\mu}$ and
$B^{\nu}$ for $\mu , \nu < \lambda$.  This proves that $B^1$ generates $B$ as
an $A$-algebra.  Since $B$ is of finite type over $A$ (a consequence of the
fact that we have supposed all of our schemes noetherian), we can choose a
finite
number of elements of $x_1,\ldots , x_m \in B^1$ which generate $B$ as an
$A$-algebra, and these elements give an embedding $V\subset \Oo ^m$.  This
embedding is linear, since the elements are elements of $B^1$ (from the above
discussion one sees that for $b\in B^1$ we have $b(u+v)=b(u)+ b(v)$). Write  $$
B= A[x_1,\ldots , x_m]/I
$$
for a homogeneous ideal $I=\bigoplus I^{\lambda}$.  We claim that $I$ is
generated as an ideal by $I^1$.  To see this, let $I'$ be the ideal generated
by $I^1$ and put $B'= A[x]/I'$.  Under the comultiplication of $A[x]$ we have
$$
I^1 \rightarrow I^1 \otimes _A A \oplus A \otimes _A I^1 ,
$$
so $I^1$ maps to zero in $B'\otimes _A B'$. Thus so does $I'$. We obtain a
comultiplication
$$
B'= A[x]/I'\rightarrow B'\otimes _A B',
$$
so $Spec (B')$ is a vector scheme too.  But the map $B'\rightarrow B$ is
surjective and an isomorphism on the pieces of degree $1$. It is compatible
with the comultiplication.  We claim that it is injective, showing this on the
part of degree $\lambda$ by induction on $\lambda$ (starting at $\lambda =2$).
If an element $b\in (B')^{\lambda}$ maps to zero in $B$, then by applying the
process given above (in the algebra $B'$) we can write $b= \sum b_{\mu} b_{\nu}$
for $\mu , \nu < \lambda$. But $b_{\mu}$ and $b_{\nu}$ map to the elements in
$B$
given by applying the same process to the image of $b$; as this image is $0$, so
are the images of $b_{\mu}$ and $b_{\nu}$.   By the induction hypothesis, the
map
is injective on the pieces of degrees $\mu , \nu$, so $b_{\mu}=b_{\nu}=0$,
giving $b=0$.  This induction shows that $B'\cong B$, so $I'=I$ is generated by
$I^1$.  Since $B$ is of finite type over $A$ (which is noetherian), $I$ is
generated by a finite number of elements.  This implies that it is generated by
a finite number of elements $y_1, \ldots , y_n$ of $I^1$.  These elements give
a linear map $\Oo ^n \rightarrow \Oo ^m$, and $V$ is the kernel.
\eop

We  come now to the main definition of this section.
A {\em vector sheaf on $S$} is a sheaf of abelian groups $F$ on $\Xx /S$ such
thatthere exists an etale covering $\{
S_{\alpha}\rightarrow S\}$ such that for each $\alpha$   there exists an exact
sequence
$$
U_{\alpha}\rightarrow V_{\alpha}\rightarrow F|_{S_{\alpha}}\rightarrow 0
$$
of sheaves of abelian groups, with $U_{\alpha}$ and $V_{\alpha}$ vector schemes
over $S_{\alpha}$.

Denote by $\Vv (S)$ the category of vector sheaves over $S$.

If $X\rightarrow S$ is an element of
$\Xx /S$, we denote by $F|_X$ the restriction of $F$ to the category $\Xx /X$.
It is a vector sheaf over $X$ (this is easy to see from the definitions).
If $F$ is a vector sheaf and $f\in F(Y)$ and $a:X\rightarrow Y$ is a morphism,
we denote the restriction of  $f$ to $X$ by $a^{\ast}(f)$ or just $f|_X$.

\begin{lemma}
\mylabel{I.c}
If $F$ is a vector sheaf, and $S$ is an affine variety, then the cohomology
groups $H^i (S, F)$ vanish for $i>0$.  If
$$
F_1\rightarrow F_2\rightarrow F_3
$$
is an exact sequence of vector sheaves (that is, an exact sequence in the
category of abelian sheaves on $\Xx$, where the elements are vector sheaves)
then for any $X$ over $S$ which is itself an affine scheme, the sequence
$$
F_1(X)\rightarrow F_2(X)\rightarrow F_3(X)
$$
is exact.
\end{lemma}
{\em Proof:}
Treat first the case where $F$ is a vector scheme.  We have an exact sequence
$$
0\rightarrow F\rightarrow \Oo ^a \rightarrow \Oo ^b
$$
by Lemma \ref{I.b}. Let $G$ be the kernel of the morphism $\Oo ^a \rightarrow
\Oo ^b$ on the small etale site over $S$.  It is a coherent sheaf.  Let $F'$ be
the sheaf on $\Xx$ whose value on $Y\rightarrow S$ is the space of sections of
the pullback (of coherent sheaves) of $G$ to $Y$.  There is a surjective
morphism $F'\rightarrow F$, which induces $F'(U)\stackrel{\cong}{\rightarrow}
F(U)$ for any $U$ etale over $S$ (or even any $U$ which is flat over $S$).
Let $K$ denote the kernel of $F'\rightarrow F$.  We claim that if $Y$ is any
scheme etale over $S$, then $H^i(Y, K)=0$.  Prove this by ascending induction
on $i$.  If the cohomology in degrees $<i$ of all fiber products of elements in
all etale covering families of $Y$ vanishes, then the degree $i$ sheaf
cohomology is equal to the degree $i$ \v{C}ech cohomology.  But the \v{C}ech
cohomology is calculated only in terms of the values of the sheaf on the fiber
products, and here the values of $K$ are zero.  Thus $H^i(Y,K)=\check{H}^i(Y,
K)=0$, completing the induction.  We obtain $H^i(S, F)= H^i(S, F')$.  But the
higher cohomology of a coherent sheaf on an affine scheme $S$ vanishes (even in
the big etale site).
We obtain the desired vanishing.
For the second part, suppose that $X=S$ is affine. The
restriction of the exact sequence to the small etale site (over $X$) remains
exact.  It can be completed to a $5$-term exact sequence where the first and
last terms are also coherent sheaves; then broken down into short exact
sequences.  The vanishing of $H^1$ of coherent sheaves on the small etale site
yields the desired exactness of all the short exact sequences of global
sections, and hence  the exactness of the sequence in question.
\eop

{\em Remark:} One can show that a vector sheaf $V$ over an affine $S$ has a
resolution by vector schemes, over $S$ rather than over an etale covering of
$S$ \cite{Hirschowitz}.

\begin{lemma}
\mylabel{I.d}
Suppose $F$ is a vector sheaf over $S$.  Then for any $X\in \Xx /S$ and $Y$ a
scheme of finite type over $k$, we have
$$
F(X\times _{Spec (k)}Y )=
F(X)\otimes _k \Oo (Y).
$$
The isomorphism is given by the pullback $F(X)\rightarrow F(X\times _kY)$
and the
scalar multiplication by the pullback of functions on $Y$.
\end{lemma}
{\em Proof:}
We first prove this when $F$
is a vector scheme. There is an exact sequence
$$
0\rightarrow F\rightarrow \Oo ^a \stackrel{M}{\rightarrow} \Oo ^b .
$$
We have $F(X)= \ker (M(X))$ and $F(X\times _kY)= \ker
(M(X\times _kY))$.  But $\Oo (X\times _kY)=\Oo (X)\otimes _k\Oo (Y)$,
and $M(X\times _kY)=M(X)\otimes 1$.  Since tensoring over $k$ is exact,
$$
\ker (M(X)\otimes 1) = \ker (M(X)) \otimes _k \Oo (Y)
$$
as desired.

Now suppose $F$ is a vector sheaf.  There is an
exact sequence
$$
U\rightarrow V\rightarrow F\rightarrow 0.
$$
If $Z$ is affine then the sequence
$$
U(Z)\rightarrow V(Z)\rightarrow F(Z)\rightarrow 0
$$
remains exact.  To see this,  replace $F$ by a coherent
sheaf $F'$ on the small etale site over $Z$.  The restriction of $F$ to the
small etale site over $Z$ is the quotient $F'$ of the restriction of
$U\rightarrow V$ to the small etale site over $Z$, that is to say the sections
of $F$ and $F'$ are the same on schemes etale over $Z$ (and in particular over
$Y$).  But if $Z$ is affine, then taking global sections preserves surjectivity
of a morphism of coherent sheaves. This gives the desired exact sequence
(proceed in a similar way for exactness at $V(Z)$).  Suppose now that $X$ and
$Y$ are affine.  Then applying the above to $Z=X$ and $Z=X\times _k Y$ we get
$$
U(X)\otimes _k \Oo (Y) \rightarrow V(X)\otimes _k\Oo (Y) \rightarrow
F(X\times _kY)\rightarrow 0 .
$$
The first morphism is the same as in the tensor product of
$$
U(X)\rightarrow
V(X)\rightarrow F(X)\rightarrow 0
$$
with $\Oo (Y)$, so the two quotients are isomorphic:  $F(X\times _kY)\cong
F(X)\otimes _k \Oo (Y)$. This completes the case where $X$ and $Y$ are
affine. But both sides of the equation have the property that they are
sheaves in each variable $X$ and $Y$ separately; thus we may first localize
on $X$ and then localize on $Y$, to reduce to the case where $X$ and $Y$ are
affine.

Finally, suppose $F$ is a
vector sheaf, and write $F= \bigcup _{i\in I} F_i$ as a directed union of vector
sheaves.  The tensor product of the union is equal to the union of the tensor
products:
$$
F(X)\otimes _k\Oo (Y) = \bigcup _{i\in I} F_i (X)\otimes _k\Oo (Y)  =
\bigcup _{i\in I} F(X\otimes _kY) = F(X\otimes _kY).
$$
Note that
the inclusion maps in the two directed unions are the same (since the
isomorphisms established above are uniquely determined by compatibility with the
morphisms $F_i (X)\rightarrow F_i (X\otimes _kY)$ and with scalar
multiplication by elements of $\Oo (Y)$). This completes the proof.
\eop

{\em Remark:}  We will mostly use this lemma in the following two cases.
Suppose $F$ is a vector sheaf over $S$.  Then for any $X\in \Xx /S$ we have
$F(X\times {\bf A}^1)= F(X)\otimes _k k[t]$. The isomorphism is given by the
pullback $F(X)\rightarrow F(X\times {\bf A}^1)$ and the scalar multiplication
by the pullback of the coordinate function $t$ on ${\bf A}^1$.
Similarly, $F(X\times {\bf G} _m) = F(X)\otimes _k k[t,t^{-1}]$, with
the  isomorphism uniquely determined by compatibility with the previous one
under
the inclusion ${\bf G} _m \subset {\bf A}^1$.

\begin{lemma}
\mylabel{I.e}
A vector sheaf has a unique structure of $\Oo$-module, and any morphism of
vector sheaves is automatically compatible with the $\Oo$-module structure.
\end{lemma}
{\em Proof:}
Suppose that $\phi :F\rightarrow G$ is a morphism of vector sheaves.  Suppose
$X\in \Xx /S$.  Suppose $f\in F(X\otimes {\bf A}^1)$.  The difference $g=\phi
(tf)-t\phi (f )$ is an element of $G(X\times {\bf A}^1)$ which restricts to zero
on  $X\otimes \{ n\}$ for any integer $n$.  We can write
$$
g=\sum _{i=1}^pg_i t^i
$$
with $g^i\in G(X)$ (by the previous lemma).  We know that
$$
g(n)=\sum _{i=1}^pg_i n^i = 0
$$
for any integer $n$.  But in $k$ the matrix $(n^i)_{1\leq n, i\leq p}$
has an inverse $(c_{ni})$, and we have
$$
g_i = \sum _{n=1}^p c_{ni}g(n) = 0.
$$
Therefore $\phi
(tf)-t\phi (f )=g=0$, for any $f$.  Thus $\phi $ is compatible with
multiplication by $t$.  Now suppose $\lambda \in \Oo (X)$.  This gives a
morphism $\gamma : X\rightarrow X\times {\bf A} ^1$ such that $\gamma ^{\ast}
(tp_1^{\ast}(f))= \lambda f$ for any $f\in F(X)$ or $G(X)$ (here $p_1:X\times
{\bf A}^1\rightarrow X$ is the projection). The fact that $\phi$ is a
morphism of sheaves means that it is compatible with $\gamma ^{\ast}$ and
$p_1^{\ast}$, so we have
$$
\phi (\lambda f)= \phi (\gamma ^{\ast} (tp_1^{\ast}(f)))=
\gamma ^{\ast} (\phi (tp_1^{\ast}(f)))
$$
$$
= \gamma
^{\ast}(t\phi (p_1^{\ast} (f)))= \gamma ^{\ast}(tp_1^{\ast}(\phi (f)))=
\lambda \phi (f).
$$
Thus $\phi$ is compatible with scalar multiplication.  This fact, applied to
the identity of $F$, implies that the scalar multiplication is unique if it
exists.

For existence, note that any morphism of vector schemes is automatically a
morphism of $\Oo$-modules, so the quotient has a structure of $\Oo$-module.
Thus any vector sheaf has a structure of $\Oo$-module.  If $F$ is a
vector sheaf expressed as a directed union $F= \bigcup _{i\in I} F_i$ of finite
vector sheaves, then the inclusions in the directed union are  compatible
with the $\Oo$-module structures; thus the union has an $\Oo$-module
structure.
\eop

The conclusion of this lemma is that the category of vector sheaves, with
morphisms equal to those morphisms of abelian sheaves compatible with the
$\Oo$-module structure, is a full subcategory of the category of sheaves of
abelian groups on $\Xx /S$.

Next we establish a Krull-type property.
\begin{lemma}
\mylabel{I.f}
Suppose that $F$ is a vector sheaf over $S$, with $f\in F(Y)$, and suppose that
for every $X\rightarrow Y$ where $X$ is an artinian scheme, $f|_X=0$.  Then
$f=0$.  Suppose $\phi : F\rightarrow G$ is a morphism of vector sheaves such
that for every $X\rightarrow S$ with $X$ artinian, $\phi |_X=0$. Then $\phi
=0$.
\end{lemma}
{\em Proof:}
We work with vector schemes over base schemes which are not necessarily of
finite type over $k$ (the definition is the same, but we require
additionally that the vector scheme be of finite type over the base). If
$U\rightarrow V$ is a morphism of vector schemes over a henselian local ring
$A$, and if $v$ is a section of $V$ over $A$ such that for each $n$ there exists
$u_n\in U(Spec (A/{\bf m}^n))$ with $u_n$ mapping to the restriction of $v$,
then
there exists a section $u$ of $U$ over $A$ which maps to $v$. This follows
from the strong Artin approximation theorem at maximal ideals, applied to
finding sections of the morphism $U\times _VSpec (A)\rightarrow Spec (A)$.

Now onto the proof of the lemma.
For the first statement, any section $f$ is contained in a vector
subsheaf of $F$, so we may suppose that $F$ is a vector sheaf.  Choose a
presentation
$$
U\rightarrow V \rightarrow F \rightarrow 0
$$
by vector schemes. We may replace $X$ by a covering, so we may suppose that our
section $f$ comes from a section $v$ of $V$.  From the previous paragraph, for
every henselized local ring $A$ of $X$, there exists a section $u$ of $U (Spec
(A)$ mapping to $v$.  But any such $A$---henselization at a point $P$---is the
direct limit of algebras $A_i$ etale of finite type over $X$ (which give etale
neighborhoods of $P$), and the space of sections is the direct limit:
$$
U(Spec (A))= \lim _{\rightarrow } U( Spec (A_i )).
$$
Thus there is a section $u_i$ over some $Spec (A_i)$ mapping to $v$.  Thus
every point $P$ of $X$ has an etale neighborhood on which there is a lifting
of $v$ to a section of $U$.  This implies that the image of $v$ in the cokernel
$F$ in the etale topology, is zero.  This gives the first statement, and
the second statement follows easily from this.
\eop

{\em Remark:} An alternative to the above proof is to use Lemma
\ref{Krull}.

The utility of this property comes from the following fact.
\begin{corollary}
\mylabel{I.g}
If $F$ is a vector scheme, and if $Y\rightarrow S$ is an element of
$\Xx /S$ with $Y$ artinian, then the functor $F_Y: Z\mapsto
F(Y\times _{Spec (k)}Z)$
from schemes over $Spec (k)$ to sets, is represented by an additive group
scheme (that is, a finite dimensional vector space) over $k$. This vector space
is the $k$-module $F(Y)$.
\end{corollary}
{\em Proof:}
By Lemma \ref{I.d}, we have $F_Y(Z)= F(Y)\otimes _k \Oo (Z)$ which is the space
of morphisms of schemes from $Z$ to the vector space $F(Y)$. Thus $F_Y$ is
represented by the vector space $F(Y)$.  Note that from the exact sequences
used in the proof of Lemma \ref{I.d}, $F(Y)$ is a finite-dimensional $k$-vector
space.
\eop

The group scheme ${\bf G} _m$ acts on every vector sheaf, by scalar
multiplication.  This action may be thought of as an action of the functor
${\bf G}
_m (X)$ on $F(X)$, or as an automorphism of $F(X\otimes {\bf G} _m )$
(multiplication by $t$) which is natural in $X$.  We have seen above that if
$F\rightarrow G$ is a morphism of sheaves of abelian groups between two vector
sheaves, then it is compatible with the ${\bf G} _m$ action.

Suppose $A$ is a vector scheme, and $F$ is a vector sheaf.
We look at $F(A)$, the space of sections over the scheme $A$. Let
$F(A)^{\lambda}$ denote the subgroup of elements $f\in F(A)$ such that
$f(ta)=t^{\lambda} f(a)$.  Here $a\mapsto ta$ is considered as a morphism
$A\times {\bf G} _m\rightarrow A$ over $S$, and $f(a)\mapsto
t^{\lambda}f(a)$ is the automorphism of $F(A\times {\bf G} _m)$ given by scalar
multiplication by $t^{\lambda} \in k[t,t^{-1}]$; the notation $f$ in the
second half of the formula actually denotes the pullback of $f$ to $A\times
{\bf G}
_m$.

\begin{lemma}
\mylabel{I.h}
With the above notations, $F(A)$ decomposes as a direct sum
$$
F(A) = \bigoplus _{\lambda \in \zz ,\lambda \geq 0}  F(A)^{\lambda} .
$$
This direct sum decomposition is natural with respect to morphisms
$F\rightarrow G$, and the linear piece $F(A)^1$ is exactly the space of
morphisms of vector sheaves $A\rightarrow F$.
\end{lemma}
{\em Proof:}
Recall that $F(A\times {\bf A}^1)= F(A)\otimes _kk[t]$, which we will just
write as $F(A)[t]$.  The morphism  of scalar multiplication $A\times{\bf
A}^1\rightarrow A$ gives $\Psi _t: F(A)\rightarrow F(A)[t]$ defined
by $(\Psi _tf)(a):= f(ta)$ (to be accurate, this should be defined in terms of
restriction maps for the morphisms involved, but we keep this notation for
simplicity). Then $F(A)^{\lambda}$ is the set of $f$ such that $\Psi _tf =
t^{\lambda}f$ in $F(A)[t]$.  Let $\Psi _s[t]: F(A)[t]\rightarrow F(A)[s,t]$
denote the extension of $\Psi _s: F(A)\rightarrow F(A)[s]$ to the polynomials in
$t$.  We have
$$
(\Psi _s[t]\Psi _tf)(a)= f(tsa)= (\Psi _{st}f)(a).
$$
Write
$$
\Psi _t(f)= \sum _{i=0}^{\infty} \psi _i(f)t^i,
$$
where $\psi _i(f)\in F(A)$ and for any $f$, there are only a finite number of
nonzero $\psi _i (f)$. Our previous formula
becomes
$$
\sum _{i,j} \psi _i (\psi _j (f))s^it^j = \sum _k \psi _k (f)(st)^k.
$$
Comparing terms we see that $\psi _i(\psi _j (f))=0$ for $i\neq j$ and $\psi
_i(\psi _i(f))=\psi _i (f)$.  But in general $f\in F(A)^{\lambda}$ if and only
if
$\psi _i(f)=0$ for $i\neq \lambda$ and $\psi _{\lambda}(f)=f$. Therefore $\psi
_i (f)\in F(A)^i$.  Restrict to $t=1$, and note that the composed morphism
$a\mapsto (a,1)\mapsto a$ is the identity  so $\Psi _1(f)=f$.  We get
$$
f= \sum _{i=0}^{\infty} \psi _i (f) ,
$$
and this sum is actually finite. Thus every element of $F(A)$ can be expressed
as a finite sum of elements of the $F(A)^{\lambda}$.  On the other hand, this
expression is unique:  if $f= \sum f_i$ with $f_i\in F(A)^i$ then
$$
\sum \psi _i(f)t^i=\Psi _t(f)=\sum \Psi _t(f_i)= \sum \psi _i (f_i)t^i = \sum
f_i t^i,
$$
and comparing coefficients of $t^i$ we get $f_i = \psi _i (f)$.
This completes the proof of the decomposition (note that in working with ${\bf
A}^1$ instead of ${\bf G} _m$ we obtain automatically that the exponents are
positive).

We have to show that $F(A)^1$ is equal to the space of linear morphisms from $A$
to $F$. A linear morphism gives an element of $F(A)^1$ (since it is compatible
with the action of $\Oo$ by Lemma \ref{I.e}), and the resulting map from the
space of morphisms to $F(A)^1$ is injective, since $F(A)$ is the space of
morphisms of functors $A\rightarrow F$.

Finally, we show surjectivity. For this, suppose
given an element $\phi \in F(A)^1$.
Suppose $Y$ is artinian, and $F$ is a vector sheaf over $S$. Then the
functor $Z\mapsto F(Y\times Z )$ is represented by a vector space $F_Y$ over
$k$ (Lemma \ref{I.g}).   Our element of $F(A)$ now gives a morphism of schemes
$\phi _Y: A_Y \rightarrow F_Y$ between these two vector spaces. It is compatible
with scalar multiplication, so it is linear. In particular, if $u,v\in A_Y(Spec
(k))= A(Y)$ then $\phi (u+v)= \phi (u)+\phi (v)$ in $F_Y(Spec (k))= F(Y)$.
Now suppose $X$ is any element of $\Xx /S$.  We show that $\phi : A(X)
\rightarrow
F(X)$ is a morphism of abelian groups.  Suppose $u,v\in A(X)$.  Let $f=\phi
(u+v)-\phi (u)-\phi (v)\in F(X)$.  By the previous paragraph, for any
$Y\rightarrow X$ with $Y$ artinian, we have $f|_Y=0$.  But the Krull property
of Lemma \ref{I.f} then implies that $f=0$.  This shows that $\phi$ is a
morphism
of sheaves of abelian groups.
\eop

If $F$ and $G$ are sheaves of abelian groups, we denote by $Hom (F,G)$ the
internal $Hom$, that is the sheaf of homomorphisms of sheaves of abelian
groups from $F$ to $G$.  The value $Hom (F,G)(X)$ is the space of morphisms of
sheaves of abelian groups from $F|_{\Xx /X}$ to  $G|_{\Xx /X}$ (this is already
a sheaf).

\begin{corollary}
\mylabel{I.i}
If $F\rightarrow G $ is a surjection of vector sheaves, and if $A$
is a vector scheme, then the morphism sheaves
$$
Hom (A, F) \rightarrow Hom (A, G)
$$
is surjective. If $X$ is affine then
$$
Hom (A,F)(X)\rightarrow Hom (A, G)(X)
$$
is surjective.
\end{corollary}
{\em Proof:}
It suffices to prove the second statement. We may assume that $X=S$.  We have
$$
Hom (A, G)(S)= G(A)^1
$$
by the last statement of the lemma.  Since $A$ is
affine, the morphism $F(A)\rightarrow G(A)$ is surjective, and by the previous
lemma this implies that $Hom (A,F)(S)=F(A)^1\rightarrow G(A)^1$ is surjective,
giving the corollary.
\eop

\begin{corollary}
\mylabel{I.j}
If $\phi : F\rightarrow G$ is a morphism of vector sheaves, then ${\rm
coker}(\phi )$ and ${\rm ker} (\phi )$ are vector sheaves.
\end{corollary}
{\em Proof:}
We may suppose that $S$ is affine and small enough. Choose presentations
by vector schemes (cf the remark before Lemma
\ref{I.d})
$$
U\rightarrow V\rightarrow F\rightarrow 0
$$
and
$$
0\rightarrow P \rightarrow R\rightarrow T \rightarrow G \rightarrow 0
$$
(note that the kernel $P$ is automatically a vector scheme).
The morphism $V\rightarrow G$ lifts to a morphism $V\rightarrow T$, by the
previous corollary, and we obtain a presentation
$$
R \oplus V \rightarrow T \rightarrow {\rm coker} (\phi ) \rightarrow 0.
$$
The fiber
products $V\times _T R$ and $V\times _T R$ are vector schemes,
and we have a presentation
$$
U\times _T R\rightarrow V\times _T R \rightarrow {\rm ker}(\phi )\rightarrow
0.
$$
\eop

Now we have shown that the category of vector sheaves is an abelian subcategory
of the category of sheaves of abelian groups on $\Xx /S$.

Suppose $A$ is a vector scheme and $F$ is a vector sheaf.  Let ${\bf 3}$ denote
the automorphism of $A$ obtained by multiplication by the scalar $3$ (any
integer $\neq 0, \pm 1$ will do).  We have
$$
F(A)=\bigoplus F(A)^{\lambda }
$$
(the decomposition given by Lemma \ref{I.h}) where $F(A)^{\lambda}$ may be
characterized as the subspace of elements $f$ such that ${\bf 3}^{\ast}(f)=
3^{\lambda}f$.  In particular, the linear subspace $Hom (A, F)= F(A)^1$ is
characterized as the subspace of elements $f$ such that
${\bf 3}^{\ast}(f)=
3f$.

\begin{theorem}
\mylabel{I.k}
Suppose $E$ and $G$ are vector sheaves, and
$$
0\rightarrow E \rightarrow F\rightarrow G \rightarrow 0
$$
is an extension in the category of sheaves of abelian groups on $\Xx /S$.  Then
$F$ is a vector sheaf.
\end{theorem}
{\em Proof:}
We proceed in several steps.   We may assume that $S$ is affine and small
enough.  Let
$$
\begin{array}{ccccccc}
              & V           &        &           &    & B           &     \\
              & \downarrow  &        &           &    & \downarrow  &     \\
              & U           &        &           &    & A           &     \\
              & \downarrow  &        &           &    & \downarrow  &     \\
0  \rightarrow &   E   &\rightarrow   & F & \rightarrow  & G & \rightarrow 0 \\
              & \downarrow  &        &           &    & \downarrow  &     \\
              & 0           &        &           &    & 0           &
\end{array}
$$
be presentations for $E$ and $G$.

{\em Step 1.}  {\em There exists a lifting of the morphism $A\rightarrow G$ to
an element $\phi \in F(A)$ with $({\bf 3}^{\ast} - 3)^2\phi =0$.}
The cohomology of $E$ over the affine $S$ is zero, so
$$
0\rightarrow E(A)\rightarrow F(A)\rightarrow G(A)\rightarrow 0
$$
is  exact.  Let $\alpha : A\rightarrow G$ denote the morphism in the
presentation above, and choose $f\in F(A)$ mapping to $\alpha$.  Then write
$$
({\bf 3}^{\ast} - 3)f = \sum e_i
$$
with $e_{\lambda} \in E(A)^{\lambda}$ (thus ${\bf 3}^{\ast} e_{\lambda} =
3^{\lambda}e_{\lambda}$). Let
$$
\phi = f-\sum c_{\lambda} e_{\lambda}
$$
for $c_{\lambda} = (3^{\lambda }-3)^{-1}$ when $\lambda \neq 1$ (and $c_1=0$).
We then have
\begin{eqnarray*}
({\bf 3}^{\ast} - 3)\phi &=& ({\bf 3}^{\ast} - 3)f -\sum c_{\lambda}({\bf
3}^{\ast} - 3) e_{\lambda}  \\
&=&\sum e_{\lambda}-\sum c_{\lambda}(3^{\lambda }-3)e_{\lambda} \\
&=& e_1.
\end{eqnarray*}
On the other hand, $({\bf 3}^{\ast} - 3)e_1=0$, so we get
$$
({\bf 3}^{\ast} - 3)^2\phi =0.
$$
In other words, $\phi$ is in the generalized eigenspace for the eigenvalue $3$
of the transformation ${\bf 3}^{\ast}$.

{\em Step 2.} {\em The extension $F$ satisfies the Krull property of Lemma
\ref{I.f}: if $f\in F(X)$ such that for any artinian $Y\rightarrow X$,
$f|_Y=0$, then $f=0$.}
Under these hypotheses, $f$ maps to an element $g\in G(X)$ satisfying the same
vanishing, so by Lemma \ref{I.f} we have $g=0$; thus $f$ comes from an element
$e\in E(X)$. This element again satisfies the same vanishing, so by Lemma
\ref{I.f}, $e=0$.

{\em Step 3.}  {\em If $A$ is a vector scheme and $F$ is an extension of two
vector sheaves, then any element $\phi \in F(A)$ with  $({\bf 3}^{\ast} -
3)^2\phi =0$ is a morphism of sheaves of abelian groups from $A$ to $F$.}
Suppose $Y$ is artinian, and $G$ is a vector sheaf over $S$. Then the
functor $Z\mapsto G(Y\times Z )$ is represented by a vector space $G_Y$ over
$k$.  If $F$ is an extension of two finite vetor sheaves $E$ and
$G$, then let $F_Y$ denote the functor $Z\mapsto F(Y\times Z)$. We obtain an
extension  $$
0\rightarrow E_Y \rightarrow F_Y \rightarrow G_Y \rightarrow 0
$$
in the category of sheaves of abelian groups over $Spec (k)$.  But since the
cohomology of the affine space $G_Y$ with coefficients in the additive group
$E_Y$ vanishes, there is a lifting of the identity to a section $u\in F_Y(G_Y)$.
Using $u$ we obtain an isomorphism of functors $F_Y \cong E_Y \times G_Y$, so
$F_Y$ is a scheme. Since $F_Y$ is a sheaf of abelian groups, $F_Y$ is an
abelian group-scheme over $k$.  Since it is an extension of two additive
groups, it is additive.  Our element of $F(A)$ now gives a morphism of schemes
$\phi _Y: A_Y \rightarrow F_Y$ between these two vector spaces.
We still have $({\bf 3}^{\ast}-3)^2)\phi _Y =0$.  But $F_Y(A_Y)$ decomposes
into eigenspaces
$$
F_Y(A_Y)=\bigoplus F_Y(A_Y)^{\lambda}
$$
where $f\in F_Y(A_Y)^{\lambda} \Leftrightarrow f(ta)= t^{\lambda}f(a)$.
In particular, $F_Y(A_Y)$ is the $3^{\lambda}$-eigenspace for ${\bf
3}^{\ast}$.  But since the space $F_Y(A_Y)$ is the direct sum of eigenspaces,
the generalized eigenspaces are equal to the eigenspaces, so $\phi _Y \in
F_Y(A_Y)^1= Hom (A_Y, F_Y)$. In particular, if $u,v\in A_Y(Spec (k))= A(Y)$
then $\phi (u+v)= \phi (u)+\phi (v)$ in $F_Y(Spec (k))= F(Y)$.

Now suppose $X$ is any element of $\Xx /S$.  We show that $\phi : A(X)
\rightarrow
F(X)$ is a morphism of abelian groups.  Suppose $u,v\in A(X)$.  Let $f=\phi
(u+v)-\phi (u)-\phi (v)\in F(X)$.  By the previous paragraph, for any
$Y\rightarrow X$ with $Y$ artinian, we have $f|_Y=0$.  But the Krull property
of Step 2 then implies that $f=0$.  This shows that $\phi$ is a morphism of
sheaves of abelian groups.

{\em Step 4.}  {\em There is a surjection from a vector scheme to $F$.}
The direct sum of the morphism $U\rightarrow F$ with our lifting $\phi :
A\rightarrow F$ gives a surjection $U\oplus A \rightarrow F \rightarrow 0$.
In fact, this fits into a diagram
$$
\begin{array}{ccccccc}
0  \rightarrow &U& \rightarrow & U\oplus A& \rightarrow  & A&\rightarrow 0\\
              & \downarrow  &     &  \downarrow &       & \downarrow &     \\
0 \rightarrow &  E   &\rightarrow   & F    & \rightarrow  & G & \rightarrow 0\\
              & \downarrow  &     &  \downarrow &       & \downarrow &     \\
              & 0           &        &   0        &    & 0    .       &
\end{array}
$$

{\em Step 5.}  {\em There is a surjection from a vector scheme to the
kernel of $U\oplus A \rightarrow F$ (proving the theorem).}
Taking the kernels along the top row of the above diagram  gives
$$
\begin{array}{ccccccc}
0  \rightarrow & K  & \rightarrow   & L & \rightarrow    & M   & \rightarrow 0\\
              & \downarrow  &     &  \downarrow &       & \downarrow &     \\
0  \rightarrow &U& \rightarrow & U\oplus A& \rightarrow  & A&\rightarrow 0\\
          & \downarrow  &     &  \downarrow &       & \downarrow &    \\
0 \rightarrow &  E   &\rightarrow   & F    & \rightarrow  & G & \rightarrow 0\\
              & \downarrow  &     &  \downarrow &       & \downarrow &     \\
              & 0           &        &   0        &    & 0    .       &
\end{array}
$$
But $K$ and $M$ are vector sheaves, and we have surjections
$V\rightarrow K \rightarrow 0$ and $B\rightarrow M \rightarrow 0$.  By
repeating the above argument in this case, we obtain a surjection
$$
V\oplus B \rightarrow L \rightarrow 0,
$$
finally giving our presentation
$$
V\oplus B \rightarrow U\oplus A \rightarrow F \rightarrow 0.
$$
Thus $F$ is a vector sheaf.
\eop

Our abelian category $\Vv$ of vector sheaves is therefore closed
under extensions of sheaves of abelian groups.

\subnumero{Duality}
Suppose $F, G$ are vector sheaves.  We have defined $Hom (F,G)$ which
is for now a sheaf of abelian groups.  Put
$$
F^{\ast} := Hom (F, \Oo ).
$$
If $\phi :F\rightarrow G$ is a morphism of vector schemes, then we obtain
a morphism $\phi ^t:G^{\ast} \rightarrow F^{\ast}$, and the construction $\phi
\mapsto \phi ^t$ preserves composition (reversing the order, of course).

\begin{lemma}
\mylabel{I.l}
{\rm (Hirschowitz \cite{Hirschowitz})}
Suppose
$$
0\rightarrow U\rightarrow V\rightarrow W \rightarrow F \rightarrow 0
$$
is an exact sequence with $U$, $V$ and $W$ vector schemes. Then taking the dual
gives an exact sequence
$$
0\rightarrow F^{\ast}\rightarrow W^{\ast}\rightarrow
V^{\ast} \rightarrow U^{\ast} \rightarrow
0.
$$
\end{lemma}
{\em Proof:}
Note first that the compositions are zero, since taking the dual is
compatible with compositions (and the dual of the zero map is zero!).

The map $F^{\ast}\rightarrow W^{\ast}$ is injective
because $W\rightarrow F$ is surjective (so any morphism $F\rightarrow \Oo$
restricting to $0$ on $W$, must be zero).

The morphism $V^{\ast }\rightarrow U^{\ast}$ is surjective:  if
$a:U\rightarrow \Oo$ is a morphism, it can be interpreted as a section of $\Oo
(U)^1$; but since $U\subset V$ is a closed subscheme, we can extend this to a
section $a'\in \Oo (V)$, then let $a''$ be the component of $a'$ in $\Oo
(V)^1$; restriction from $\Oo (V)$ to $\Oo (U)$  is compatible with the ${\bf G}
_m$ action, hence with the decomposition of Lemma \ref{I.h}, so $a''$ restricts
to $a$.

Suppose $b : W\rightarrow \Oo$ restricts to zero on $V$; then it factors
through the quotient sheaf  $F=W/V$, so it comes from $F^{\ast}$.  Thus the
sequence is exact at $W^{\ast}$.

We still have to prove exactness at $V^{\ast}$.  Choose embeddings
$U\hookrightarrow \Oo ^m$ and $W\hookrightarrow \Oo ^n$.  Then extend the first
to a function $V\rightarrow \Oo ^m$; combining with the second we obtain
$V\hookrightarrow \Oo ^{m+n}$, fitting into a diagram
$$
\begin{array}{ccccccc}
0\rightarrow & \Oo ^m& \rightarrow & \Oo ^{m+n} & \rightarrow  & \Oo ^n &
\rightarrow 0  \\
&\downarrow &&\downarrow && \downarrow & \\
0\rightarrow & U& \rightarrow & V & \rightarrow  & W & .
\end{array}
$$
Furthermore, $U= \Oo ^m \cap V$ as subschemes of $\Oo ^{m+n}$ (by the
injectivity of $W\rightarrow \Oo ^n$).    Given a linear map $\lambda :
V\rightarrow \Oo$ such that $\lambda |_{U}=0$, extend it to $\varphi : \Oo
^{m+n}\rightarrow \Oo$ such that $\varphi |_{\Oo ^{m}}=0$.  Replace $\varphi$
by its linear part under the decomposition of Lemma \ref{I.h} (this will
conserve the property $\varphi |_{\Oo ^{m}}=0$ as well as the property of
restricting to $\lambda$). Our $\varphi$ now descends to a map $\Oo
^n\rightarrow \Oo$, restricting to $\varphi |_W$ which extends $\lambda$.

Note in the previous paragraph, we have used the following general fact:  if
$X,Y\subset Z$ are closed subschemes of an affine scheme, and $\lambda \in \Oo
(X)$ such that $\lambda |_{X\cap Y}=0$, then there exists $\varphi \in \Oo (Z)$
such that $\varphi |_X=\lambda$ and $\varphi |_Y=0$.  To prove this, let $I_X$,
$I_Y$ and $I_{X\cap Y}$ denote the ideals of $X$, $Y$ and $X\cap Y$ in the
coordinate ring $\Oo (Z)$.  The definition of the scheme-theoretic intersection
$X\cap Y$ is that $I_{X\cap Y}= I_X+I_Y$, and our statement follows from the
translation that
$$
I_Y \rightarrow I_{X\cap Y}/I_X \subset \Oo (Z)/I_X
$$
is surjective.

We have completed the proof of the lemma.
\eop

\begin{corollary}
\mylabel{I.m}
The functor $F\mapsto F^{\ast}$ is an exact functor from the category of finite
vector sheaves, to the category of sheaves of abelian groups.
\end{corollary}
{\em Proof:}
Suppose
$$
0\rightarrow F' \rightarrow F \rightarrow F'' \rightarrow 0
$$
is an exact sequence of vector schemes.  Choose presentations
$$
0\rightarrow U'\rightarrow V'\rightarrow W'\rightarrow F' \rightarrow 0
$$
and
$$
0\rightarrow U''\rightarrow V''\rightarrow W''\rightarrow F \rightarrow 0,
$$
and combine these into a presentation
$$
0\rightarrow U\rightarrow V\rightarrow W\rightarrow F \rightarrow 0
$$
with $U=U'\oplus U''$, $V=V'\oplus V''$ and $W=W'\oplus W''$ (using the method
of Theorem \ref{I.k}, which is easier since we now have the required lifts
automatically).  These fit together into a diagram
$$
\begin{array}{ccccccc}
&0&&0&&0 & \\
 & \downarrow & & \downarrow && \downarrow \\
0\rightarrow &U'& \rightarrow &U& \rightarrow &U''& \rightarrow 0 \\
 & \downarrow & & \downarrow && \downarrow \\
0\rightarrow &V'& \rightarrow &V& \rightarrow &V''& \rightarrow 0 \\
 & \downarrow & & \downarrow && \downarrow \\
0\rightarrow &W'& \rightarrow &W& \rightarrow &W''& \rightarrow 0 \\
 & \downarrow & & \downarrow && \downarrow \\
0\rightarrow &F'& \rightarrow &F& \rightarrow &F''& \rightarrow 0 \\
 & \downarrow & & \downarrow && \downarrow \\
&0&&0&&0 &
\end{array}
$$
where all the rows and columns are exact.  Apply duality to this diagram; we
obtain a diagram with the arrows reversed, with the columns exact, by the lemma.
Furthermore, the same lemma shows that the upper three rows are exact (in fact,
this is easier because the rows in the original diagram are split, by
construction).  This implies that the bottom row is exact, as desired.
\eop

A {\em coherent sheaf} is a sheaf which (locally) has a presentation of the
form
$$
\Oo ^n\rightarrow \Oo ^m\rightarrow F\rightarrow 0.
$$
In particular, note that it is a vector sheaf.  This coincides with the
usual definition: if $S$ is affine and $X\rightarrow S$ is a morphism, then
$F(X)=F(S)\otimes _{\Oo (S)}\Oo (X)$ (this is because the same is true for
$\Oo$, and the presentation remains exact on the right after tensoring).

As usual, we can assume that a presentation as above exists globally over any
affine base.

\begin{corollary}
\mylabel{I.n}
The dual of a coherent sheaf is a vector scheme and vice-versa.
\end{corollary}
{\em Proof:}
Note that $\Oo ^{\ast}=\Oo$.  Taking the dual of a presentation of a coherent
sheaf gives
$$
0\rightarrow F^{\ast} \rightarrow \Oo ^m \rightarrow \Oo ^n,
$$
so $F^{\ast}$ is a vector scheme (the kernel here is a closed subscheme of $\Oo
^n$). Conversely, if $V$ is a vector scheme, take an exact sequence such as
given in Lemma \ref{I.b}, and apply the dual.  We obtain a presentation for
$V^{\ast}$ as a coherent sheaf.
\eop

\begin{corollary}
\mylabel{I.o}
The dual of a vector sheaf is again a vector sheaf.
\end{corollary}
{\em Proof:}
If $F$ is a  vector sheaf, choose a presentation
$$
U\rightarrow V\rightarrow F\rightarrow 0
$$
by vector schemes.  Taking the dual gives
$$
0\rightarrow F^{\ast}\rightarrow U^{\ast}\rightarrow V^{\ast}.
$$
By the previous corollary, $U^{\ast}$ and $V^{\ast}$ are coherent sheaves, in
particular  vector schemes.  Thus $F^{\ast}$ is the kernel of a morphism
of  vector sheaves, so $F^{\ast}$ is  a  vector sheaf.
\eop

\begin{lemma}
\mylabel{I.p}
If $F$ is a  vector sheaf, then $F^{\ast\ast}=F$ (via the natural
morphism).
\end{lemma}
{\em Proof:}
If $F$ is a vector scheme, this follows from the construction given in
Corollary \ref{I.n}: write $F=\ker (M)$ as the kernel of a matrix
$M:\Oo ^m\rightarrow \Oo ^n$; then $F^{\ast} = {\rm coker} (M^t)$ is
the cokernel
of the transpose matrix (and this $M^t$ is really just the transpose,
keeping the
same coefficients as in $M$).  Finally, $F^{\ast\ast}=\ker (M^{tt})$, but the
transpose of the transpose is the same matrix $M=M^{tt}$, so $F=F^{\ast\ast}$.
(The same argument works for coherent sheaves, of course).  If $F$ is any
 vector scheme, choose a presentation
$$
U\rightarrow V\rightarrow F\rightarrow 0
$$
and take the double dual.  Since $U^{\ast\ast}=U$ and $V^{\ast\ast}=V$ we get
$$
U\rightarrow V\rightarrow F^{\ast\ast}\rightarrow 0,
$$
so $F^{\ast\ast}=F$.
\eop

We have now shown that duality is an exact contravariant involution on the
category $\Vv $ of  vector sheaves, interchanging vector schemes and
coherent sheaves.

\begin{lemma}
\mylabel{I.q}
The vector schemes are projective objects in $\Vv $, and the coherent
sheaves are injective objects.  There exist enough projectives and injectives
(assuming that $S$ is affine).
\end{lemma}
{\em Proof:}
The argument given above shows that a vector scheme $A$ is a projective object:
if $F\rightarrow G$ is a surjection of  vector sheaves then, since $A$
is affine, $F(A)^1\rightarrow G(A)^1$ is surjective.  By definition, every
 vector sheaf admits a surjection from a vector scheme, so there are
enough projectives.  By duality, the coherent sheaves are injective and there
are enough injectives.
\eop

Taking the dual of the three step resolution by vector schemes shows that every
 vector sheaf $F$ admits a resolution
$$
0\rightarrow F \rightarrow U\rightarrow V\rightarrow W\rightarrow 0,
$$
with $U$, $V$ and $W$ coherent sheaves (in particular, injective).

\subnumero{Internal $Hom$ and tensor products}

We begin with a corollary to the last lemma.
\begin{corollary}
\mylabel{I.r}
If $A$ is a vector scheme, then the functor $V\mapsto Hom (A,V)$ from $\Vv $
to the category of abelian sheaves, is exact.  If $F$ is a coherent sheaf, then
the functor $V\mapsto Hom (V,F)$ is exact.
\end{corollary}
{\em Proof:}
If $S$ is affine, the functors $V\mapsto Hom (A,V)(S)$ and $V\mapsto Hom
(V,F)(S)$ are exact, by the lemma.  But the restriction of a vector scheme or a
coherent sheaf, to any object $X\in \Xx /S$ is again a vector scheme or coherent
sheaf over $X$, so we obtain exactness over every affine object; and since
exactness is a local condition, we get exactness.
\eop

\begin{lemma}
\mylabel{I.s}
If $F$ and $G$ are  vector sheaves, then $Hom (F,G)$ is a  vector
sheaf.
\end{lemma}
{\em Proof:}
Suppose $F$ and $G$ are vector schemes. Then the exact sequence
$$
0\rightarrow G\rightarrow \Oo ^a \rightarrow \Oo ^b
$$
yields an exact sequence
$$
0\rightarrow Hom (F,G)\rightarrow Hom (F, \Oo ^a) \rightarrow Hom (F,\Oo ^b);
$$
but the middle and right terms are direct sums of the dual $F^{\ast}$ which is
a  vector sheaf, so the kernel $Hom (F,G)$ is a  vector sheaf.
Now suppose $F$ is a vector scheme and $G$ is a  vector sheaf; resolving
$G$ by vector schemes we obtain a resolution of $Hom (F,G)$ by  vector
sheaves, from the previous sentence.  Thus $Hom (F,G)$ is a  vector sheaf
in this case too. Now suppose $F$ is a  vector scheme, and choose a
resolution
$$
U\rightarrow V\rightarrow F\rightarrow 0
$$
by vector schemes.  The functor $W\mapsto Hom (W,G)$ is contravariant and left
exact for any $G$, so we obtain an exact sequence
$$
0\rightarrow Hom (F,G)\rightarrow Hom (V,G)\rightarrow Hom (U,G).
$$
The middle and right terms are  vector sheaves by the previous arguments,
so the kernel is also.  This completes the proof in general.
\eop

We now define the {\em tensor product} $F\otimes ^{\Vv} G$ of two vector
sheaves to be
$$
F\otimes _{\Oo}G:= (Hom (F, G^{\ast}))^{\ast}.
$$
Beware that this is not just the tensor product of sheaves of $\Oo$-modules
(although this will be the case if $F$ and $G$ are coherent sheaves).
We can also define the {\em cotensor product}
$$
F\otimes ^{\Oo} G:= Hom (F^{\ast} , G).
$$
Again, beware here that this is not equal to the tensor product.  The
difference is seen in noting that the tensor product is right exact as usual,
whereas the cotensor product is left exact.  (These exactness statements hold
in both variables since the tensor and cotensor products are commutative, as we
see below).  Duality permutes the tensor and cotensor products:
$$
(F\otimes _{\Oo}G)^{\ast}= F^{\ast}\otimes ^{\Oo}G^{\ast}
$$
and
$$
(F\otimes ^{\Oo}G)^{\ast}= F^{\ast}\otimes _{\Oo}G^{\ast}.
$$

Define recursively
$$
V_1\otimes \ldots \otimes V_n := V_1\otimes (V_2\otimes \ldots \otimes V_{n})
$$
for either one of the tensor products.

By {\em multilinear form} $V_1\times \ldots V_n\rightarrow W$ we mean simply a
multilinear morphism of sheaves of groups. In the same way as above for the
linear morphisms, we obtain a vector sheaf $Mult(V_1\times \ldots \times V_n ,
W)$ of multilinear forms (denoted $Bil (\;\;\; )$ when $n=2$).

\begin{proposition}
\mylabel{I.s.1}
1. \,\, There is a natural isomorphism
$\alpha _{U,V}:Hom (U^{\ast}, V)\cong Hom (V^{\ast}, U)$ and $\alpha
_{U,V}\alpha _{V,U}$ is the identity.
\newline
2. \,\, There is a natural isomorphism
$$
Multi (V_1\times \ldots \times V_n, W)\cong
Hom (W^{\ast}, Multi (V_1\times \ldots \times V_n, \Oo ).
$$
3. \,\, There is a natural isomorphism
$$
Multi (V_1\times \ldots \times V_n , W)\cong Hom (V_1, Multi (V_2\times\ldots
\times V_n, W)).
$$
\end{proposition}
{\em Proof:}
In each case one defines natural maps in both directions and checks that the
two compositions are the identity.
\eop

\begin{theorem}
\mylabel{I.s.2}
Suppose $V_i$ are vector sheaves, $i=1,\ldots , n$.
There is a multilinear form
$$
\mu : V_1\times \ldots V_n \rightarrow V_1 \otimes _{\Oo} \ldots \otimes _{\Oo}
V_n
$$
which is universal in the sense that if
$$
\phi: V_1\times \ldots \times V_n\rightarrow W
$$
is a
multilinear form then there is a unique morphism
$$
\psi : V_1 \otimes _{\Oo} \ldots \otimes _{\Oo}
V_n \rightarrow W
$$
such that $\phi = \psi \circ \mu $.
\end{theorem}
{\em Proof:}
Note first that for $n=2$ there is a natural bilinear map $U\times V
\rightarrow Hom (U, V^{\ast})^{\ast}= U\otimes _{\Oo}V$. Inductively this
gives the multilinear map for any $n$.   The universal property says that the
induced map
$$
Hom (V_1\otimes _{\Oo}\ldots \otimes _{\Oo}V_n,W)\rightarrow Multi (V_1,\ldots
, V_n , W)
$$
should be an isomorphism. We prove this by induction on $n$, so we may suppose
it is true for $n-1$.  By the definition of the multiple tensor product, the
quantity on the left is
$$
Hom (Hom (V_1, (V_2\otimes _{\Oo}\ldots \otimes _{\Oo}V_n)^{\ast})^{\ast},W).
$$
By part 1 of the proposition, this is equal to
$$
Hom (W^{\ast}, Hom (V_1, (V_2\otimes _{\Oo}\ldots \otimes
_{\Oo}V_n)^{\ast})).
$$
By induction, $(V_2\otimes _{\Oo}\ldots \otimes
_{\Oo}V_n)^{\ast}=Multi (V_2\times \ldots \times V_n,\Oo )$. Coupled with part
3 of the proposition we get
$$
Hom (W^{\ast}, Multi (V_1\times \ldots \times V_n,\Oo )
$$
which then is equal to the right hand side above, by part 2 of the proposition.
\eop

\begin{corollary}
\mylabel{I.s.3}
The tensor and cotensor products have natural commutativity and associativity
isomorphisms satisfying the usual constraints.
\end{corollary}
{\em Proof:}
For the tensor product this follows from the universal property and the fact
that the notion of multilinear form is independent of the order of the
variables. For the cotensor product this follows because it is the dual of the
tensor product.
\eop

We can define symmetric and exterior powers, either with respect to the tensor
product or with respect to the cotensor product.  Let $S_n$ denote the
symmetric group on $n$ objects.  Let $V^{\otimes _{\Oo}n}$ (resp.
$V^{\otimes ^{\Oo}n}$) denote the tensor product (resp. cotensor product) of $n$
copies of a vector sheaf $V$. Then $S_n$ acts on $V^{\otimes _{\Oo}n}$
(resp. $V^{\otimes ^{\Oo}n}$) because of the commutativity and associativity.
This representation is completely reducible (this can be see object-by-object).
The components are vector sheaves; this can be seen by noting that the fixed
part is the kernel of a morphism of vector sheaves (the direct sum of $1-\gamma
$ for $\gamma \in S_n$); to get the other components, apply that to tensor
products with the irreducible representations of $S_n$.  The trivial component
of
$V^{\otimes _{\Oo}n}$ (resp. $V^{\otimes ^{\Oo}n}$) is denoted by $Sym
_{\Oo}^n(V)$, the symmetric power (resp. $Sym ^{\Oo}_n(V)$, the symmetric
copower).  The component corresponding to the sign representation is denoted by
$\bigwedge _{\Oo}^n(V)$, the exterior power (resp. $\bigwedge ^{\Oo}_n(V)$, the
exterior copower).

{\em Remark:}
There is a natural morphism $U\otimes _{\Oo}V\rightarrow U\otimes ^{\Oo}V$,
that is to say
$$
Hom (U,V^{\ast})^{\ast}\rightarrow Hom (U^{\ast}, V).
$$
However, this is not an isomorphism. A counterexample can be constructed by
looking for a case where the cotensor product is left exact but not right exact
(and noting that the tensor product is right exact), or vice-versa.

We have the expression $U\otimes ^{\Oo}V= Bil (U^{\ast}\times V^{\ast},\Oo)$.
The inequality mentionned above implies in particular that there are
bilinear functions on $U^{\ast}\times V^{\ast}$ which are not sums of tensors
$u\otimes v$.  This is a big difference from the case of schemes (for example
if $U$ and $V$ are coherent so that $U^{\ast}$ and $V^{\ast}$ are vector
schemes, then the bilinear functions {\em are} sums of tensor  products.

\subnumero{Automorphisms of vector sheaves}

We end our discussion of vector sheaves by showing how they give examples of
presentable group sheaves.

\begin{lemma}
\mylabel{I.1.g.1}
If $V$ is a vector sheaf over $S$, then $V$ is presentable.
\end{lemma}
{\em Proof:}
Suppose $V$ is a vector scheme.  Then taking $X=V$ and $R= V\times _VV=V$
we obtain the required presentation (note that the identity morphisms are
vertical)---so $V$ is $P4$, and then $P5$ by Corollary
\ref{I.z}.  It follows from Theorem \ref{I.1.d} that the quotient of one vector
scheme by another is again $P5$; and finally that the quotient of a
vector scheme by such a quotient is $P5$.  In view of the 3-stage resolution of
any vector sheaf by vector schemes, we obtain the lemma.
\eop

One of the main examples of presentable group sheaves is given by the following
theorem.

\begin{theorem}
\mylabel{I.1.g}
Suppose $V$ is a vector sheaf over $S$.  Then the group sheaf $Aut (V)$ is a
presentable.
\end{theorem}
{\em Proof:}
By the previous lemma and Lemma \ref{I.s}, $Hom (V,V)$ is $P5$. We can express
$$
Aut (V)\subset Hom (V,V)\times Hom (V,V)
$$
as the equalizer of the two morphisms
$$
\begin{array}{ccc}
Hom (V,V)\times Hom (V,V)&\rightarrow &Hom (V,V)\times Hom (V,V)\\
(a,b) & \mapsto & (ab,ba) \\
(a,b)&\mapsto & (1,1).
\end{array}
$$
Apply Lemma \ref{I.1.a} to obtain that $Aut (V)$ is $P4$, and then Corollary
\ref{I.z} to obtain that it is $P5$.
\eop

A particular case of this construction is when $V$ is a coherent sheaf which we
denote by $\Ff$.  There is a presentation
$$
U_2 \stackrel{\phi}{\rightarrow} U_1 \rightarrow \Ff \rightarrow 0
$$
where $U_i = \Oo ^{a_i}$. Let $Aut (U_2, U_1, \phi )$ denote the group sheaf of
automorphisms of the morphism $U_2 \rightarrow U_1$.  Any such automorphism
gives
an automorphism of $\Ff$ so we have a morphism
$$
Aut (U_2, U_1, \phi )\rightarrow Aut (\Ff ).
$$

\begin{lemma}
\mylabel{surjection}
This morphism is a surjection onto $Aut (\Ff )$, and $Aut (U_2, U_1, \phi )$
is represesentable by a group scheme over $S$.
\end{lemma}
{\em Proof:}
The representability by a group scheme is clear, since $Aut (U_i)$ are group
schemes (isomorphic to $GL(a_i)$) and the condition of compatibility with $\phi$
is a closed condition so $Aut (U_2, U_1, \phi )$ is a closed subscheme of $Aut
(U_1)\times Aut (U_2)$.

Suppose $S' \rightarrow S$ is a scheme and $P\in S'$ is a point.  Suppose
$\eta :
\Ff |_{S'}\rightarrow \Ff |_{S'}$ is an automorphism. Let
$$
U'_2 \stackrel{\phi '}{\rightarrow }U'_1 \rightarrow \Ff |_{S'}\rightarrow 0
$$
be a minimal resolution of $\Ff |_{S'}$ at the point $P$ (that is to say that
the
value $\phi '(P)$ is identically zero and the rank of $U'_2$ is minimal).  Then
there are locally free $W_i\cong \Oo ^{b_i}$ on $S'$ such that $U_i|_{S'} \cong
U'_i \oplus W$ and such that the map $\phi |_{S'}$ can be written in block
form with respect to this decomposition, with a morphism
$\psi '$ in the block of the $W_i$ and
the map $\phi '$ in the block of the $U'_i$, such that $\psi '$ is surjective.
Our
morphism $\eta$ extends to a morphism of resolutions
$U'_{\cdot} \rightarrow U'_{\cdot}$ which is an isomorphism near $P$ by the
minimality of the resolution (in fact the values $U'_i(P)$ are the $Tor
^i_{\Oo _{S'}}(\Ff |_{S'}, k_P)$ and an isomorphism of $\Ff |_{S'}$ induces an
isomorphism on the $Tor ^i$).  We can complete this with the identity in the
block of the $W_i$ to get an isomorphism of resolutions $U_i |_{S'}$ inducing
$\eta$.  This gives the desired surjectivity.
\eop

{\em Question:}  Does a similar result hold for the automorphisms of any vector
sheaf?

\numero{Tangent sheaves of presentable sheaves}

Suppose $S'\rightarrow S$ is an $S$-scheme.  Put
$$
Y:= S' \times Spec (k[\epsilon _1
,\epsilon _2, \epsilon _3]/(\epsilon _i^2, \epsilon _i\epsilon _j ))
$$
with the subschemes
$$
Y_i:= S' \times Spec (k[\epsilon _i ]/(\epsilon
_i^2))
$$
and
$$
Y_{ij}:= S' \times Spec (k[\epsilon _i ,\epsilon _j]/(\epsilon
_i^2, \epsilon _j^2, \epsilon _i\epsilon _j )).
$$
Note that $Y=Y_1\cup Y_2\cup Y_3$, and $Y_{ij}=Y_i\cup Y_j$, as well as
$Y_i\cap Y_k =
S'\subset Y$ and
$Y_{ij}\cap
Y_{jk}=Y_j$ (for $i\neq k$). It
should be stated explicitly that $Y_{ij}$ is the closed subscheme defined by the
ideal $(\epsilon _k)$, $k\neq i,j$; and $Y_i$ is the closed subscheme defined
by
the ideal $(\epsilon _j,\epsilon _k)$, $j,k\neq i$.

We need a weaker version of the notion of verticality.
We say that a morphism $\Ff \rightarrow \Gg$ of sheaves is {\em $T$-vertical}
if it satisfies the lifting property $Lift_2(Y_{ij};Y_i ,Y_j)$
and $Lift _3(Y; Y_{12},Y_{23}, Y_{13})$ (for any $S'$). Note that these systems
satisfy the retraction hypotheses in the lifting property, so the property of
$T$-verticality is weaker than the property of verticality.

The result of Theorem \ref{I.u} holds also for $T$-verticality, so the class
$\Tt$ of $T$-vertical morphisms satisfies the axioms M1-M4. In particular
the properties $P4$ and $P5$ imply $P4(\Tt )$ and $P5(\Tt)$ respectively.

The advantage of th weaker property of $T$-verticality is that if
$X\rightarrow Z$
is a morphism of schemes over $S$, then it is $T$-vertical.  To prove this, note
that  the properties $Y=Y_1\cup Y_2\cup Y_3$, $Y_{ij}=Y_i\cup Y_j$,
$Y_i\cap Y_k =
S'\subset Y$ and $Y_{ij}\cap
Y_{jk}=Y_j$ mean that for defining morphisms from $Y$ to a scheme (or from
$Y_{ij}$ to a scheme) it suffices to have compatible morphisms on the $Y_{ij}$
or on the $Y_i$.

({\em Caution:} We did not include the lifting condition $Lift _1(Y_1; S')$
in the notion of $T$-verticality; morphisms of schemes do not necessarily
satisfy this lifting property!)

The conclusion of the previous paragraph and property $M1$ for
$T$-verticality is that if $\Ff$ is a $P4(\Tt )$ sheaf then the structural
morphism $p:\Ff \rightarrow S$ is $T$-vertical; thus $P4(\Tt )\Leftrightarrow
P5(\Tt)$.

\begin{lemma}
\mylabel{I.1.e.1}
Suppose $f:\Ff\rightarrow \Gg $ is a morphism of $P4$ sheaves.
Then $f$ is $T$-vertical.
Furthermore, the
liftings in the lifting properties for $f$, for   the systems
$(Y_{ij};Y_i ,Y_j)$
and $(Y; Y_{12},Y_{23}, Y_{13})$, are unique.
\end{lemma}
{\em Proof:}
For $T$-verticality, we can choose vertical surjections $X\rightarrow \Ff$ and
$Y\rightarrow \Gg$ so that there is a lifting $X\rightarrow Y$.  This lifting
is $T$-vertical since it is a morphism between schemes (cf the above remark).
Hence the composition $X\rightarrow \Gg$ is $T$-vertical.  By Theorem \ref{I.u},
part 4 for $T$-verticality, applied to the composition $X\rightarrow \Ff
\rightarrow \Gg$, we obtain $T$-verticality of the morphism $f$.

To prove the uniqueness, note that liftings to schemes are unique
since  $Y=Y_1\cup Y_2\cup Y_3$ and $Y_{ij}=Y_i\cup Y_j$.  Then descend the
uniqueness down from $X$ to $\Ff$ where $X\rightarrow \Ff$ is the vertical
(hence
$T$-vertical) morphism provided by the property $P4$.  This descent of the
uniqueness property is immediate from the lifting property for $X\rightarrow
\Ff$.
\eop

In the statement of the following theorem, the condition is $P4$ and not
$P4(\Tt )$ (i.e. that isn't a misprint).

\begin{theorem}
\mylabel{I.1.f}
Suppose $\Ff\rightarrow \Gg $ is a morphism of $P4$ sheaves on $S$.  Suppose $u
:S\rightarrow \Ff$ is a section.  Then the relative tangent sheaf $T(f )_{u}$
over $S$, defined by
$$
T(\Ff )_{u} (b:S'\rightarrow S):= \{ \eta :S' \times Spec (k[\epsilon
]/(\epsilon ^2))\rightarrow \Ff\;\; :\;\;\;\; f \eta = fubp_1  \;\; \mbox{and}
\;\; \eta |_{S'}= ub \} ,
$$
has a natural structure of sheaf of abelian groups making it a vector sheaf.
\end{theorem}
{\em Proof:}
We first define the natural abelian group structure on this sheaf.  Suppose
$$
\eta _i: S'\times Spec(k[\epsilon _i]/(\epsilon _i^2))\rightarrow \Ff
$$
are sections of $T(f )_u$ over $S'$ ($i=1,,\ldots , 3$). ({\em Nota:} for the
definition of the group law we only need $i=1,2$; we need $i=1,2,3$ only to
check that it is associative.) Here (and below) we attach various subscripts to
the variables $\epsilon$. Use the notations established above:
$$
Y:= S' \times Spec (k[\epsilon _1
,\epsilon _2, \epsilon _3]/(\epsilon _i^2, \epsilon _i\epsilon _j ))
$$
with the subschemes
$$
Y_i:= S' \times Spec (k[\epsilon _i ]/(\epsilon
_i^2))
$$
and
$$
Y_{ij}:= S' \times Spec (k[\epsilon _i ,\epsilon _j]/(\epsilon
_i^2, \epsilon _j^2, \epsilon _i\epsilon _j )).
$$
Note that $Y=Y_1\cup Y_2\cup Y_3$, and $Y_{ij}=Y_i\cup Y_j$.  Again $Y_{ij}$ is
the closed subscheme defined by the ideal $(\epsilon _k)$, $k\neq i,j$; and
$Y_i$ is the closed subscheme defined by the ideal $(\epsilon _j,\epsilon _k)$,
$j,k\neq i$.   The systems $(Y_{ij};Y_i ,Y_j)$ and $(Y; Y_{12},Y_{23}, Y_{13})$
satisfy a unique lifting property for the morphism $f$ (Lemma \ref{I.1.e.1}).
Note that $Y_{ij}\cap Y_{jk}=Y_j$ (for $i\neq k$). We apply this first to the
system $(Y_{ij}; Y_i,Y_j)$. There is a unique morphism
$$
\eta _{ij}: Y_{ij}\rightarrow \Ff
$$
over the base morphism $Y_{ij}\rightarrow S\rightarrow \Gg$ and
agreeing with $\eta _i$ (resp. $\eta _j$) on $Y_i$ (resp. $Y_j$). Let
$$
\delta _{ij}: S' \times Spec (k[\epsilon ]/(\epsilon ^2))\rightarrow
Y_{ij}
$$
be the diagonal and---for future use---let
$$
\delta _{123}: S' \times Spec (k[\epsilon ]/(\epsilon ^2))\rightarrow
Y
$$
be the triple
diagonal.
Then we put
$$
\eta _i+\eta _j:= \eta _{ij} \circ \delta _{ij} .
$$
This gives a composition which is obviously commutative (the definition is
symmetric in the two variables).  To check that it is associative, apply unique
lifting for $(Y,Y_{ij})$ to get a unique $\eta _{123}: Y\rightarrow \Ff$
restricting to the $\eta _{ij}$ on $Y_{ij}$. Next,  note that the triple
diagonal
is equal to the composition of $1\times \delta _{23}$ with the diagonal
$$
S' \times Spec (k[\epsilon _0]/(\epsilon _0^2))\rightarrow
Spec (k[\epsilon _1,\epsilon ]/(\epsilon _1^2, \epsilon ^2, \epsilon
_1\epsilon )).
$$
Using this, we get
$$
\epsilon _1+(\eta _2 +\eta _3)= \eta _{123}\circ \delta _{123}.
$$
Similarly, we have
$$
(\epsilon _1+\eta _2) +\eta _3= \eta _{123}\circ \delta _{123},
$$
giving associativity.

The identity element (which we denote by $0$) is the composition
$$
S'\times Spec (k[\epsilon ]/(\epsilon ^2))\rightarrow S \rightarrow \Ff .
$$

This construction is natural: if
$$
\begin{array}{ccc}
\Ff & \rightarrow & \Ff ' \\
\downarrow &&\downarrow \\
\Gg & \rightarrow & \Gg '
\end{array}
$$
is a diagram with vertical arrows vertical, and if $u:S\rightarrow \Ff$
is a section  projecting to $u':S\rightarrow \Ff '$, then composition with the
morphism $\Ff \rightarrow \Ff '$ respects the conditions in the definition of
the tangent sheaves, and so it gives a morphism $T(f)_u\rightarrow
T(f')_{u'}$.  The addition we have defined is natural, so this morphism
of tangent sheaves respects the addition (it also respects the identity).

The inverse is obtained by applying the
automorphism $\epsilon \mapsto -\epsilon$.
This completes the construction of the natural structure of sheaf of abelian
groups.

Next, we show that if
$$
\Ff \stackrel{a}{\rightarrow }\Gg \stackrel{b}{\rightarrow }\Hh
$$
is a sequence of morphisms of $P4$ sheaves, and if $u:S\rightarrow \Ff$ is
a section, then we have an exact sequence
$$
0\rightarrow T(a)_u\rightarrow T(ba)_u \rightarrow T(b)_{au}
$$
We certainly get such a sequence with the
composition being zero. Furthermore, $T(a)_u$ is the subsheaf of $T(ba)_u$
consisting of those elements projecting to zero in $T(b)_{au}$ (this follows
immediately from the definition).

Furthermore, if $a$ is vertical, then the sequence is exact on the right. This
follows from the lifting property in the definition of vertical, in view of the
fact that $S'$ is a retraction of $S'\times Spec (k[\epsilon ]/(\epsilon ^2))$.
(Note that we have not required this lifting property in the definition of
$T$-verticality.)

Let $p: \Ff \rightarrow S$ denote the structural morphism for a $P4$ sheaf
$\Ff$,
and define the tangent sheaf $T(\Ff )_u:= T(p)_u$.  If $f:\Ff \rightarrow
\Gg$ is
a morphism of $P4$ sheaves, the exact sequence of the previous paragraph becomes
$$
0\rightarrow T(f)_u\rightarrow T(\Ff)_u \rightarrow T(\Gg )_{fu}.
$$
Again, if $f$ is vertical then this sequence is exact on the right also.

Finally, we show that if $\Ff$ is $P4$ then $T(\Ff )_u$ is a vector sheaf.  The
above exact sequence implies that if $f$ is a morphism of $P4$ sheaves then
$T(f)_u$ is a vector sheaf. Let $f: X\rightarrow \Ff$ be the vertical morphism
given by the property $P4$. Since the question is etale local on $S$, we may
assume that our section $u: S\rightarrow \Ff$ lifts to a section $v:
S\rightarrow X$.  We have an exact sequence
$$
0\rightarrow T(f)_v \rightarrow T(X)_v\rightarrow T(\Ff )_u\rightarrow 0.
$$
Note that $T(X)_v$ is a vector scheme (an easy thing to see---it is given by the
linear parts of the equations of $X$ at the section $v$).

Let $g:R\rightarrow X\times _{\Ff} X$ be the other vertical morphism given by
the property $P4$.

We claim that we have an exact sequence
$$
0\rightarrow T(X\times _{\Ff}X)_{(v,v)} \rightarrow T(X)_v \oplus T(X)_v
\rightarrow T(\Ff )_u \rightarrow 0.
$$
To see this, note that an element of $T(X\times _{\Ff}X)_{(v,v)}$ consists of
an element of $T(X\times _SX)_{(v,v)}$ mapping to $T(\Ff )_u\subset T(\Ff
\times _S\Ff )_{(u,u)}$.  Note that
$$
T(X\times _SX)_{(v,v)}=T(X)_v \oplus T(X)_v,
$$
and
$$
T(\Ff \times _S\Ff )_{(u,u)}=T(\Ff )_u\oplus T(\Ff )_u
$$
with the map from $T(\Ff )_u$ being the diagonal.  The quotient of
$T(\Ff \times _S\Ff )_{(u,u)}$ by the diagonal $T(\Ff )_u$ is thus isomorphic
to $T(\Ff )_u$ and we obtain the exact sequence in question.  The surjectivity
on
the right is from surjectivity of $T(X)_v\rightarrow T(\Ff )_u$.

Lift $(v,v)$ to a section $w:S\rightarrow R$.
The exact sequence for $g$ gives a surjection
$$
T(R)_w \rightarrow T(X\times _{\Ff}X)_{(v,v)}\rightarrow 0.
$$
Combining this with the above exact sequence, we obtain the right exact sequence
$$
T(R)_w \rightarrow T(X)_v \oplus T(X)_v
\rightarrow T(\Ff )_u \rightarrow 0.
$$
Since $T(R)_w$ and $T(X)_v$ are vector schemes, this shows that $T(\Ff )_u$
is a vector sheaf.
\eop

\numero{The case $S=Spec (k)$}
We now analyse the definitions of the previous sections in the case where the
base scheme is $S=Spec (k)$ (a hypothesis we suppose for the rest of this
section).

{\em Caution:}  We will use throughout this section certain properties of
vertical morphisms etc.  which hold only in the context $S=Spec (k)$.  The
reader should not extrapolate these properties to other cases.

Our first lemma is a preliminary version of the next lemma which we include
because the argument may be easier to understand in a simpler context.

\begin{lemma}
\mylabel{I.1.k}
Suppose $f: X\rightarrow Spec (k)$ is morphism of finite type.
Then $f$ is vertical if and only if $f$ is a smooth morphism.
\end{lemma}
{\em Proof:}
Suppose $X$ is smooth.  Then the required lifting properties hold.  Indeed, $X$
is etale locally a vector space, and Theorem \ref{I.u} (part 7) implies that
vector
spaces are vertical over $Spec(k)$.

Conversely, suppose $f$ is vertical, and suppose $x\in X$.  The first claim is
that for any $v\in T(X)_x$ there is a smooth germ of curve $(C,0)$ mapping to
$(X,x)$ with tangent vector $v$ at the origin.  Since $X$ is of finite type,
and by Artin approximation, it suffices to construct a compatible family of
morphisms
$$
\gamma _n:Spec (k[t]/t^n)\rightarrow X
$$
sending $Spec(k)$ to $x$ and with tangent vector $v$ (that is, the map
$\gamma _2$ represents $v$).  Before starting the construction, choose a
morphism
$$
\mu : X\times X\rightarrow X
$$
with $\mu (x,y)=\mu (y,x)=y$ for any $y$ (the possibility of finding $\mu$
follows from the definition of verticality). We now construct $\gamma _n$ by
induction, starting with $\gamma _2$ given by $v$.  Suppose we
have constructed $\gamma _{n}$ by the inductive procedure.  Let $Y(n):= Spec
(k[r]/r^{n})\times Spec (k[s]/s^2)$. The composition gives a morphism
$$
\phi _n:= \mu \circ (\gamma _n , \gamma _2): Y(n)\rightarrow X.
$$
We will show that $\phi _n$ factors through the morphism
$$
d: Y(n)\rightarrow Spec (k[t]/t^{n+1})
$$
which is dual to the morphism
$$
k[t]/t^{n+1} \rightarrow k[r,s]/(r^n,s^2)
$$
$$
t\mapsto r+s.
$$
We will then choose $\gamma _{n+1}$ equal to the resulting morphism
$Spec (k[t]/t^{n+1})\rightarrow X$, that is with $\phi _n =\gamma _{n+1}d$.
Since $\gamma _n$ restricts to $\gamma _{n-1}$, and since we have chosen
$\gamma _n$ by the inductive procedure, we have that
$$
\phi _n|_{Y(n-1)}=\phi _{n-1} = \gamma _n d.
$$
Writing $X=Spec (A)$ (in a neighborhood of $x$) the morphism $\phi _n$
corresponds to
$$
\phi _n^{\ast}: A\rightarrow k[r,s]/(r^n,s^2).
$$
We have that $\phi _n^{\ast} (a)$ reduces modulo $r^{n-1}$ to
$d^{\ast}\gamma _n^{\ast} (a)$.  Writing
$$
\gamma _n^{\ast}(a)= \sum _{j=0}^{n-1}b_jt^j
$$
we
have
$$
\phi _n^{\ast} (a)= \sum _{j=0}^{n-1}b_j (r+s)^j + \alpha r^{n-1} +
\beta r^{n-1}s.
$$
Write, on the other hand, the equation $\phi _n |_{Spec (k[r]/r^n)} = \gamma
_n$.  We get that
$$
\phi _n^{\ast} (a) \sim \sum _{j=0}^{n-1}b_j r^j \;\; \mbox{mod} (s).
$$
This gives $\alpha = 0$ in the above equation.  Finally, note that $(r+s)^n=
nr^{n-1}s$ modulo $(r^n, s^2)$.  Thus we may set $b_n:= \beta / n$ and
get
$$
\phi _n^{\ast} (a)= \sum _{j=0}^{n}b_j (r+s)^j .
$$
Put
$$
\gamma _{n+1}^{\ast} (a):= \sum _{j=0}^{n}b_j t^j ,
$$
and we get the desired factorization $\phi _n = \gamma _{n+1}d$.  This
completes the inductive step for the construction of the $\gamma _n$.  We
obtain the desired formal curve and hence a curve $(C,0)$ as claimed.

{\em Remark:}  Intuitively what we have done above is to integrate the vector
field on $X$ given by the tangent vector $v$ and the multiplication $\mu$.
Of course, the curve $C$ is an approximation to the integral curve, which
might only exist formally.

The next step in the proof of the lemma is to choose a collection of vectors
$v_1,\ldots , v_m$ generating $T(X)_x$, and to choose resulting curves $C_1,
\ldots ,  C_m$.  Using the map $\mu$ in succession (or applying directly the
definition of verticality) we obtain a map
$$
\Phi : (U,0):=(C_1\times \ldots \times C_m , 0)\rightarrow (X,x),
$$
inducing the given morphisms on the factors $C_i$ (considered as subspaces of
the product by putting the origin in the other places).  By construction the
differential $d\Phi _0$ is given by the vectors $v_1,\ldots , v_m$, in
particular it gives a surjection
$$
d\Phi _0: T(U)_0 \rightarrow T(X)_x \rightarrow 0.
$$
Note that $U$ is smooth of dimension $m$.
We claim that this implies $dim _x (X) \geq dim T(X)_x$.  To see this, let $d:=
dim _x(X)$ and $n:= dim T(X)_x$.  By semicontinuity, the dimension of the fiber
$\Phi ^{-1}(x)$ at the origin is at least equal to $m-d$.  In particular, the
tangent space to the fiber has dimension at least $m-d$; but this gives a
subspace of dimension $m-d$ of $T(U)_0$ which maps to zero in $T(X)_x$; by the
surjectivity of $d\Phi _0$ we get $n \leq m-(m-d)=d$, the desired inequality.

Finally, it follows from this inequality that $X$ is regular at $x$ and hence
smooth at $x$ (and, of course, the inequality is an equality!). This proves
the lemma.
\eop

\begin{lemma}
\mylabel{smooth}
Suppose $X$ and $Y$ are schemes of finite type over $k$ and $f: X\rightarrow Y$
is a morphism. Then $f$ is $Spec (k)$-vertical if and only if $f$ is smooth.
\end{lemma}
{\em Proof:}
Note first that if $f$ is smooth then it is etale-locally a product with affine
space so we get all of the lifting properties.
Suppose now that $f$ is vertical. If $Q\in Y$ and $P\in f^{-1}(Q)$
then $Lift _1(Y, Q)$ implies that, after replacing $Y$ by an etale neighborhood
of $Q$ we may suppose that there is a section $\sigma : Y\rightarrow X$
with $\sigma (Q)=P$.  Let $T(X/Y)_{\sigma}$ denote the relative tangent vector
scheme along the section $\sigma$.
It is easy to see that the morphism $T(X/Y)_{\sigma}\rightarrow Y$ is
$Spec(k)$-vertical.  We then obtain that the morphism
$$
\Gamma (Y, T(X/Y)_{\sigma})\rightarrow (T(X/Y)_{\sigma})_Q=T(f^{-1}(Q))_P
$$
is surjective, and this then implies that $T(X/Y)_{\sigma}$ is a vector bundle
over $Y$.  The same argument as in the previous lemma allows us to
``exponentiate'' in a formal neighborhood of $P$, to get a map $\varphi$ from
$T(X/Y)_{\sigma}^{\wedge}$ (the formal completion in a neighborhood of
$0(Q)$)  to $X$, which sends the zero section $0$ to $\sigma$ and whose tangent
map is the identity along $\sigma$.

We claim that if $S'$ is artinian local with a morphism
$S''\rightarrow X$ sending the origin to $P$, then the morphism factors via
$\varphi$ through a map $S'\rightarrow T(X/Y)_{\sigma}^{\wedge}$ sending the
origin to $0(Q)$. Prove this claim using the standard deformation theory
argument by induction on the length of $S'$: suppose $S''\subset S'$ is defined
by an ideal $I$ annihilated by the maximal ideal, and suppose we know the claim
for $S''$.  Then there exists a map $S'\rightarrow T(X/Y)_{\sigma}^{\wedge}$
extending the known map on $S''$ since $T(X/Y)_{\sigma}^{\wedge}$ is a vector
bundle over $Y$. The space of such extensions is a principal homogeneous space
over $I\otimes _k (T(X/Y)_{\sigma})_Q$ whereas the space of extensions of
$S''\rightarrow X$ to morphisms $S'\rightarrow X$ is a principal homogeneous
space over $I\otimes _kT(f^{-1}(Q))_P$. The map $\varphi$ induces an isomorphism
$$
(T(X/Y)_{\sigma})_Q\cong T(f^{-1}(Q))_P
$$
so there is an extension to a map $S' \rightarrow T(X/Y)_{\sigma}^{\wedge}$
which projects to our given map $S'\rightarrow X$. This proves the claim.

Now we can prove that $X\rightarrow Y$ is formally smooth at $P$. If $S''\subset
S'$ are artinian local and if $a:S'\rightarrow Y$ is a map lifting over $S''$ to
a map $b:S'' \rightarrow X$ sending the origin to $P$, then we get
(from the previous claim) that
the map $b$ factors through a map $S'' \rightarrow
T(X/Y)_{\sigma}^{\wedge}$. Since $T(X/Y)_{\sigma}^{\wedge}$ is a vector bundle
and in particular smooth over $Y$, this extends to a map
$S' \rightarrow
T(X/Y)_{\sigma}^{\wedge}$. This extension projects into $X$ to an extension
$S' \rightarrow X$ of the map $b$. This shows formal smoothness.   Since $X$
and $Y$ are of finite type, $f$ is smooth.
\eop

\begin{corollary}
\mylabel{I.1.l}
Suppose $G$ is a presentable  sheaf of groups on $\Xx /Spec (k)$ (which is
equal to
$\Xx$ in this case), and suppose $f:X\rightarrow G$ is a vertical morphism.
Then
$X$ is smooth over $Spec (k)$.
\end{corollary}
{\em Proof:}
The morphism $G\rightarrow Spec (k)$ is vertical by Theorem \ref{I.u} (7).
The composed morphism $X\rightarrow Spec (k)$ is vertical  hence smooth
by Lemma \ref{I.1.k}.
\eop

\begin{theorem}
\mylabel{I.1.m}
If $G$ is a presentable group sheaf on $\Xx /Spec (k)$ then it is represented
by a
smooth separated scheme of finite type over $k$ (in other words it is an
algebraic Lie group over $k$).
\end{theorem}
{\em Proof:}
We assume $k=\cc$ for this proof.
Choose vertical surjections
$f:X\rightarrow G$ and $R\rightarrow X\times _GX$.  Note that $R\rightarrow G$
is
vertical, so $X$ and $R$ are smooth schemes of finite type.  By adding some
factors of affine spaces we can assume that the components of $X$ and $R$ all
have the same dimension.
By the previous
section, the morphism $df:T(X)\rightarrow f^{\ast}T(G)$ is a morphism of vector
sheaves on $X$, hence it is a morphism of vector bundles.  It is surjective, so
the kernel is a strict sub-vector bundle $\Ff \subset T(X)$.
For each $x\in X$ we have
$$
\Ff _x:= \ker (T(X)_x \rightarrow T(G)_{f(x)}).
$$
The morphism $p_1: R\rightarrow X$ is vertical (since $X\times _GX\rightarrow
X$ is the pullback of the vertical $X\rightarrow G$ by the morphism
$X\rightarrow G$, and $p_1$ is the composition of the vertical $R\rightarrow
X\times _GX$ with this projection).  Therefore, by Lemma \ref{smooth} $p_1$
is smooth.  Suppose $r\in R$ maps to $(x,y)\in X\times X$. Let $g\in G$ denote
the common image of $x$ and $y$.  We have an exact sequence
$$
T(R)_r \rightarrow T(X)_x \oplus T(X)_y \rightarrow T(G)_g \rightarrow 0.
$$
From this we get that the image of the map on the left always has the same
dimension; in particular this shows that the map $T(R)\rightarrow (p_1,
p_2)^{\ast}T(X\times X)$ is strict. For any point $g$ in $G$ we can identify
$T(G)_g\cong T(G)_1$ by left multiplication. The morphism on the right in the
exact sequence then comes from a morphism of the form $p_1^{\ast}(\alpha
)-p_2^{\ast}(\alpha )$ where $\alpha : T(X) \rightarrow T(G)_1$ is obtained
from the differential of $f$ by the left-multiplication trivialization.  This
morphism is a morphism of vector bundles from the tangent bundle of $X\times X$
to the constant bundle $T(G)_1$,
so its kernel is a distribution in the tangent
bundle of $X\times X$.  The image of $R$ is an integral leaf of this
distribution.  In particular, the image of $R$ is a smooth complex submanifold
of $X\times X$ (note that the map from $R$ to the leaf is smooth since, by the
above exact sequence, the differential is surjective at any point---this implies
that the image is open in the leaf).

Choose a subvariety $X'\subset X$ which is everywhere transverse to the
distribution $\Ff$, and which meets every subvariety of $X$ of the form
$p_2(p_1^{-1}(x))$ for $p_i$ denoting the projections $R\rightarrow X$.  We may
assume that $X'$ is of finite type.  Let $R'$ be the intersection of $X'\times
X'$ with the image of $R$ in $X\times X$.  We claim that the morphism
$X'\rightarrow G$ is surjective and vertical, and that $R'= X'\times _GX'$.  To
see this, note that by hypothesis $X'\times _X R\rightarrow X$ is surjective on
closed points.  By our transversality assumptions this morphism is also smooth.
Thus any point in $X$ is equivalent via $R$ (etale-locally) to a point in
$X'$.
For verticality, it suffices to prove that $X'\times _GX \rightarrow X$ is
vertical (Theorem \ref{I.u}, parts 3 and 4).  And for this it suffices to note
that  $X'\times _X R \rightarrow X' \times _GX$ is surjective and vertical
(being the pullback of $X\times _XR\rightarrow X\times _GX$ by $X'\times
_GX\rightarrow X\times _GX$), that $X'\times _XR\rightarrow X$ is smooth and
hence vertical, and to apply Theorem \ref{I.u}, part 4.  We get $X'\rightarrow
G$ surjective and vertical.  If we put $R'' $ equal to the pullback of $R$ to
$X'\times X'$ then $R'' \rightarrow X'\times _G X'$ is surjective and vertical
(it being also the pullback of $R$ via $X'\times _GX' \rightarrow X\times
_GX$). The previous proof applied to this case shows that $R''$ is smooth over
its image $R'$, and that $R'$ is a smooth subvariety of $X'\times X'$.  But
now, by our previous transversality assumptions, the projections $R'\rightarrow
X'$ are etale.

We can now conclude that $G$, which is the quotient of $X'$ by the equivalence
relation $R'$, is a smooth algebraic space.  We will find an open subset
$U\subset G$ which is a smooth variety over $k$.  In order to do this,
let $d$ be the maximum number of  points in the fibers of $X'\rightarrow G$.
The fiber through a point $x$ is equal to $p_2(p_1^{-1}(x))$ where $p_i: R'
\rightarrow X'$ here denote the projections.  Let $W\subset X$ be the set of
points $x$ where the maximum number $d$ of points in the fiber $p_1^{-1}(x)$ is
achieved.  Since the morphism $p_1: R'\rightarrow X$ is etale, it is easy to see
that $W$ is an open subset, and that if we let $R'_W $ denote $p_1^{-1}(W)$
then $R'_W\rightarrow W$ is a finite etale morphism of degree $d$.  On the
other hand, if $x\in W$ and $y$ is in the fiber through $x$ then $y$ is also in
$W$.  This means that $p_2(R'_W)\subset W$.  The correspondence
$$
x\mapsto p_2(p_1^{-1}(x))
$$
gives a morphism $\chi$ from $W$ to the symmetric product $W^{(d)}$ having image
in the complement of the singular locus.  Then $W\times _{W^{(d)}}W= R'_W$.  In
particular, the quotient of $W$ by the equivalence relation $R'_W$ is the
image of $\chi$.  Note that $\chi $ is etale over its image, which is thus a
locally closed subscheme of $W^{(d)}$.  This shows that the quotient of $W$ by
the equivalence relation is a scheme $U$ of finite type. It is also smooth.
The morphism $W\rightarrow G$ factors through $U\rightarrow G$.

We claim
that the morphism $U\rightarrow G$ is an open subfunctor, that is for any
$Y\rightarrow G$ the fiber product $U\times _GY$ is an open subset of $Y$.
The fiber product is the quotient of $W\times _GY$ by the
induced equivalence relation; and the quotient of $X'\times _GY$ by the
equivalence relation is equal to $Y$.  Choosing local liftings
$Y\rightarrow X'$ we find that $X'\times _GY$ is the image of $R'\times
_{X'\times X'}(X'\times Y)\rightarrow X'\times Y$, that is it is the
pullback of $R'$. In particular it is a subscheme of $X'\times Y$.  This
subscheme surjects to $Y$ by a vertical morphism, a morphism which is hence
smooth.  The image of the open subset $W\times _GY$ (which is the
intersection of $X'\times _GY$ with $W\times Y$) is therefore an open set in
$Y$.  This shows that $U\subset G$ is an open subfunctor.

We can choose a finite number of elements $g_i \in G(S)$ such that $g_i\cdot U$
cover $G$.  For the finiteness use the surjection $X\rightarrow G$ with $X$ of
finite type (in particular, quasi-compact).

We now apply Grothendieck's theorem about representability which says
that if a
functor $G$ is a sheaf (in the Zariski topology, which is the case here since
Zariski is coarser than etale), and if it is covered by a finite number of open
subfunctors $G_i$ which are representable by schemes, then the functor $G$ is
representable by a scheme (the union of the schemes $G_i$).  In our case the
$G_i$ are the $g_i\cdot U$, representable by $U$.  Since $U$ is of finite type,
the union of a finite number of copies is again of finite type.

We obtain that $G$ is a scheme of finite type.  Note that $U$ is smooth so $G$
is smooth (alternatively, use that any group scheme is smooth).  To complete
the proof we just have to show that $G$ is separated. Note first that all
connected components of $G$ must have the same dimension, so we can speak of
the dimension of $G$ without problem.
Let $\Delta \subset G\times G$ denote the diagonal.  It is preserved by the
diagonal left action of $G(k)$ on $G\times G$ (that is, the action $g(a,b)=(ga,
gb)$).  The
complement $K:=\overline{\Delta}-\Delta$ is a closed subset of $G\times G$, of
dimension strictly smaller than the dimension of $G$. But $K$ is invariant
under the diagonal left action of $G(k)$, so its image $pr_1(K)\subset G$ is
invariant by the left action of $G(k)$.  Since $dim (K)< dim (G)$ the image $pr
_1(K)$ (which is a constructible subset of dimension $\leq dim (K)$) is not
dense in $G$.  On the other hand, if $K$ were nonempty then this image, being
left invariant, would contain a right translate of $G(k)$ which is Zariski
dense.  This contradiction implies that $K$ is empty, in other words $G$ is
separated.  This completes the proof of the theorem.
\eop

{\em Application:}
Suppose $S$ is any base scheme of finite type over $Spec (k)$ now, and
suppose $S'\rightarrow S$ is an artinian scheme of finite type.  Let $\pi :
S' \rightarrow Spec (k)$ denote the structural morphism.  If $G$ is a
presentable
group sheaf over $S$ the pullback $G|_{S'}$ is presentable (Lemma
\ref{I.1.h}) and the direct image $\pi _{\ast} (G|_{S'})$ is presentable over
$Spec (k)$ (Lemma \ref{I.1.i}).  By Theorem \ref{I.1.m}, $\pi _{\ast}(G|_{S'})$
is represented by a group scheme of finite type which we denote  $G_{S'}$ over
$k$.  We have
$$
G(S')= G_{S'}(Spec (k)).
$$
Furthermore, if $X\rightarrow G$ is a vertical surjection then we obtain a
scheme of finite type $X_{S'}= \pi _{\ast}(X|_{S'})$ with a morphism
$X_{S'}\rightarrow G_{S'}$.  This morphism is smooth.

\numero{Local study of presentable group sheaves}
In this section we return to the case of general base scheme $S$ (in
particular, the hypothesis $S=Spec (k)$ is no longer in effect).

First we establish some notations for formal completions.
Suppose $G$ is a presentable group sheaf.  Let $\widehat{G}$ denote the sheaf
which associates to $Y\in \Xx$ the set of values in $G(Y)$ which restrict to
the identity on $Y^{\rm red}$.  More generally, use the same notation
$\widehat{\Ff}$ whenever $\Ff$ is a sheaf with a given section playing the role
of the identity section (usually the section in question is understood from the
context).

\subnumero{Local structure}
\begin{lemma}
\mylabel{I.1.n}
Suppose $G$ is a presentable group over a base $S$. Suppose $Z\rightarrow G$
is a vertical  surjection with $Z$ an affine scheme of finite type over $S$. Let
$T(Z)_e\rightarrow S$ be the tangent vector scheme at a lift $e$ of the identity
section. For any $s\in S$ there is an etale neighborhood
$$
e(s)\in W \stackrel{p}{\rightarrow} Z
$$
and an etale $S$-morphism $q:U\rightarrow TZ$, such that $q=p$ over the
section $e$ (which maps to the zero section of $TZ$).
\end{lemma}
{\em Proof:}
Verticality of $Z\rightarrow G$ means that
we can choose a lifting of the multiplication of $G$ to $m: Z\times Z
\rightarrow Z$ such that $m(x,e)=x$ and $m(e,y)=y$.

Let $Q: Z\rightarrow Z$ be the automorphism $Q(x):= m(x,x)$. It has the
effect of multiplication by $2$ on the tangent scheme $TZ$ at the identity
section, because
$$
\frac{\partial }{\partial x}m(x,x)(e)= \frac{\partial }{\partial x}m(x,e)+
\frac{\partial }{\partial x}m(e,x)(e)= 2 \frac{\partial x}{\partial x} =2.
$$
If we embedd $Z\subset {\bf A}^N_S$ as
a closed subscheme with the identity section going to the origin-section,
then we
may extend $Q$ to a morphism $Q': {\bf A}^N_S\rightarrow {\bf A}^N_S$ such
that $Q'$ acts by multiplication by two on the tangent space at the
origin.  Let $\widehat{{\bf A}^N_S}$ denote the formal completion of the
affine space along the origin-section.  Then $Q'$ induces an
automorphism of   $\widehat{{\bf A}^N_S}$, and it is well known---and easy
to see using power series---that such an automorphism is conjugate to its linear
part (since the eigenvalues are  different
from $1$).
We obtain an automorphism $F:
\widehat{{\bf A}^N_S}\rightarrow \widehat{{\bf A}^N_S}$ such that
$F^{-1}\circ Q'\circ F = 2$.  Let $\widehat{Z}\subset \widehat{{\bf A}^N_S}$
be the closed formal subscheme obtained by completing $Z$ at the identity
section. Note that $\widehat{Z}$ is preserved by $Q'$. Thus the image
$F(\widehat{Z})$ is a formal subscheme which is preserved by
multiplication by $2$.  It follows that it is a cone, and in particular
that the linear parts of the equations defining $F(\widehat{Z})$ vanish
on $F(\widehat{Z})$.  This means that $F(\widehat{Z})$ is included in its
tangent scheme $T(F(\widehat{Z}))$ along the identity section.
Translating back by $F$ we obtain an immersion
$$
\widehat{Z}\hookrightarrow TZ
$$
which is the identity on the tangent space at the identity section.
The image is a closed formal subscheme preserved by scalar multiplication.

For any artinian scheme $S'$ over $S$, $Z(S')$ is a
smooth scheme over $Spec (k)$ and $\widehat{Z(S')}\subset TZ(S')$ is
a closed formal subscheme at the origin, with the same Zariski tangent space,
and which is  formally preserved by scalar multiplication. Therefore
$\widehat{Z(S')}\cong \widehat{TZ(S')}$. Now  $\widehat{Z}(S')$ is the
inverse image of $e\in \widehat{Z(Spec (k))}$ via the map
$$
\widehat{Z(S')}\rightarrow \widehat{Z(Spec (k))}.
$$
The same is true for the tangent scheme $TZ$. From these properties we get
that $\widehat{Z}(S')\rightarrow \widehat{TZ}(S')$ is an isomorphism for any
$S'$.

As that holds true for all artinian schemes
$S'$ over $S$ we get that the morphism $\widehat{Z} \rightarrow \widehat{TZ}$
is an isomorphism.
Artin approximation now
gives the existence of such an isomorphism (inducing the same map on tangent
schemes along the identity section) over an etale neighborhood in $Z$, as
required for the lemma.
\eop

\subnumero{Theory of the connected component}
  We need to develop a suitable theory of the connected component of a
presentable group sheaf $G$.
\begin{theorem}
\mylabel{I.1.o}
If $G$ is a presentable group sheaf over $S$, then there is a unique subsheaf of
groups $G^0\subset G$ such that $G^0$ is presentable and such that for any
artinian $S$-scheme $S'$, we have $G^0(S')$ equal to the connected component of
$G(S')$ (when these are considered as algebraic groups over the ground field of
$S'$---cf the application at the end of the section on the situation over
$Spec (k)$).
\end{theorem}
{\em Proof:}
We first show existence.  Let $Z\rightarrow G$ be a vertical surjection with $Z$
a scheme of finite type.  Let $\sigma : S\hookrightarrow Z$ be the identity
section.  We claim that there is an open neighborhood $U\subset Z$ of $\sigma
(S)$ such that for any artinian $S$-scheme $S'$, $U(S')$ is connected.    By
Lemma \ref{I.1.n}, there is an etale neighborhood of the zero section $W
\rightarrow TZ$ and another etale morphism $W\rightarrow Z$ giving an etale
neighborhood of the section $\sigma$.  We claim that (possibly  throwing out a
closed subset of $W$ not meeting the section) we can assume that the $W(S')$
are connected. In what follows we refer to the lifting of the zero section of
$TZ$ as the section $\sigma$ of $W$.

For any given $S'$, artinian located at $s\in S$, there is a
surjection of vector spaces
$$
(TZ)(S')\rightarrow V_i \subset (TZ)(s),
$$
for some subspace $V_i$ which depends on $S'$.
If $W\rightarrow TZ$ is our etale morphism, then we have
$$
W(S')=W(s)\times _{TZ(s)}(TZ)(S')
= W(s)\times _{TZ(s)}V_i \times_{V_i}(TZ)(S'),
$$
since a point $S'\rightarrow TZ$ has a unique lifting to $W$ once the lifting
is specified on the closed point.  Thus $W(S')$ is connected if and only if,
for all subspaces $V_i \subset (TZ)(s)$ we have that $W(s)\times _{TZ(s)}V_i$
is connected.

Let $Gr (TZ)\rightarrow S$ be the disjoint union of the grassmanian schemes of
subspaces of different dimensions. It is proper over $S$. We have a universal
subscheme
$$
\Vv \subset Gr (TZ)\times _S TZ.
$$
Note that the map $\Vv \rightarrow TZ$ is proper.
Let $\tilde{W}:= W\times _{TZ} \Vv $; this is an etale covering of $\Vv$, and
is proper over $W$. Let $\tilde{W}^N\subset \tilde{W}$ be the union of the
connected components in fibers which do not pass through the section $\sigma$
(relative to $Gr (TZ)$).   Note that $\tilde{W}^N$ is a constructible subset of
$\tilde{W}$ (one can see this by noetherian induction). Let $W^N\subset W$
be the image of $\tilde{W}^N$.  It is again a constructible subset.  A point
$w\in W$ is in $W^N$ if and only if there exists a vector subspace $V_i \subset
(TZ)(s)$ such that $w$ is in a different connected component of $V_i \times
_{TZ}W$ from $\sigma (s)$.  In particular, if we choose an analytic neighborhood
of the section $\sigma$ which is isomorphic to a tubular neighborhood of the
zero-section of $TZ$, then this analytic neighborhood doesn't meet $W^N$.
Thus there is a Zariski open neighborhood of $\sigma$ not meeting $W^N$.
Since taking a Zariski open subset doesn't affect connectivity (the schemes
$W_{S'}$ in question being smooth), we may replace $W$ by this open subset
and hence assume that $W^N$ is empty.  From the discussion of the previous
paragraph, this implies that the $W_{S'}$ are connected, proving the first
claim.

Let $U$ be the image of $W$ in $Z$. Note that the set-theoretic image is an
open set and is equal to the image of the functor, since $W\rightarrow Z$ is
etale.

Let $\eta : Z\times _SZ \rightarrow Z$ be a lifting of the multiplication map
$(g,h)\mapsto gh$ such that $\eta (z, 1)= z$ and $\eta (1,z)=z$.
We claim that the  composition law $Z\times _SZ \rightarrow G$ is a vertical
morphism. Note that $Z\times _SZ\rightarrow G\times _SG$ is vertical, so it
suffices to prove that the composition $G\times _SG\rightarrow G$ is
vertical.  For this, notice that there is an isomorphism $G\times _SG\cong
G\times _S G$ sending $(a,b)$ to $(ab,b)$, and which interchanges the
multiplication and the first projection.  Since the first projection is
vertical (this comes from the fact that $G\rightarrow S$ is vertical), we
obtain that the composition law is vertical, yielding the claim.

By Lemma \ref{I.1.c}, there exists a vertical surjection
$$
R\rightarrow (Z\times _SZ )\times _G (Z\times _S Z)
$$
with $R$ a scheme of finite type.
Let $G^0\subset G$ be the image of the morphism $U\times _SU\rightarrow G$.
Then the morphism $U\times _S U\rightarrow G^0$ is a vertical surjection, and
we have a vertical  surjection
$$
R'\rightarrow (U\times _SU )\times _{G^0} (U\times _S U)
$$
obtained by letting $R'$ be the inverse image of
$(U\times _SU )\times _{G^0} (U\times _S U)$ in $R$.
Note that $R'$ is also equal to the fiber product
$$
U\times _SU\times _SU\times _SU\times _{Z\times _SZ\times _S Z\times _S Z}R,
$$
so $R'$ is a scheme of finite type over $S$.

We
claim that for any artinian $S'$, the $G^0(S')$ is equal to the connected
component of $G(S')$.  To see this,  note first of all that $G^0(S')$ is
connected (since it is the image of $U(S')\times U(S')$ which is connected).
And secondly, note that the morphism
$$
Z(S')\rightarrow G(S')
$$
is an open map (this is a map of smooth varieties---cf the section on what
happens over a field and in particular the application at the end).  Therefore
the image of $U(S')$ is an open subset $V\subset G(S')$.  It is connected
since $U(S')$ is connected.
The image of $(U\times _SU)(S')$ is equal to the
image of the multiplication map $V\times V\rightarrow G(S')$.  It is easy to
see that if $V$ is a connected Zariski open subset of an algebraic group over a
field (containing the identity), then the  image of the multiplication map is a
subgroup.  Thus $G^0(S')$ is a subgroup of $G(S')$.  It contains an open
neighborhood of the  identity and it is connected, so it is equal to the
connected component.  We claim now that $G^0$ is a sheaf of subgroups of $G$.
If $g,h\in G^0(S')$ then the product $gh$ restricts into $G^0(S'')$ for any
artinian ring $S''$ over $S'$.  The sheaf $G^0$ is P2, hence it is B1 and B2
(Theorem \ref{I.t.2}). The inverse image of the section $gh$ by the morphism
$G^0\rightarrow G$ is again B1 and B2. This inverse image is nonempty
artinian $S''$.  By Artin approximation, the inverse image has a section
locally over $S'$, and since this section is unique if it exists, it gives
a section $gh\in G^0(S')$.

We have now shown existence of $G^0$ as required by the theorem.  For
uniqueness, suppose that $G^1$ were another candidate.  Then $G^0$ and $G^1$
are both B1 and B2 subsheaves of $G$ having the same points over artinian
$S'$.  Artin approximation implies that they are equal.
\eop

We say that a presentable group sheaf $G$ is {\em connected} if $G= G^0$.
The above theorem immediately gives the characterization that $G$ is connected
if and only if $G(S')$ is connected for all artinian $S'$.

\begin{corollary}
\mylabel{connex} We have the following properties.
\newline
1. \, If $G$ is connected then any quotient group of $G$ is
connected;
\newline
2.\, Of $G$ and $H$ are connected then any extension of $G$ by $H$ is connected;
\newline
3. \, If $G$ is a connected group sheaf over a base $S$ and if $Y\rightarrow S$
is any morphism of schemes then $G|_{\Xx /Y}$ is a connected group sheaf over
$Y$; and
\newline
4. \, If $f:Y\rightarrow S$ is a finite morphism and if $G$ is a connected group
sheaf over $Y$ then $f_{\ast}(G)$ is a connected group sheaf over $S$.
\newline
5.\, If $G$ is any presentable group sheaf then the connected component $G^0$ is
the largest connected presentable subgroup.
\end{corollary}
{\em Proof:}
Items 1-3 are immediate from the characterization.  To prove 4 note that if
$S'\rightarrow S$ then $f_{\ast}(G)(S')= G(Y\times _SS')$ and $Y\times _SS'$ is
artinian, so this latter group is connected, thus by the above characterization
$f_{\ast}(G)$ is connected. To prove 5 note that if $H$ is any connected
subgroup of $G$ then $H(S') \subset G^0(S')$ for all artinian $ S'$, hence
$H\subset G^0$.
\eop

\subnumero{Finite presentable group sheaves}

We say that a presentable group sheaf $G$ is {\em finite} if $G^0=\{ 1\}$.
If $G$ is any presentable group sheaf, then the connected component $G^0$ is a
normal subgroup sheaf, and the quotient $C:=G/G^0$ is again presentable. Over
artinian $S'$, this quotient is just the group of connected components, in
particular the connected component is trivial.  Thus $C$ is finite.

\begin{lemma}
\mylabel{I.1.p}
If $G$ is a finite presentable group sheaf, then there is an integer $N$
such that for any henselian local $S$-scheme $S'$ (with algebraically closed
residue field), the number of elements in $G(S')$ is less than or equal to $N$.
\end{lemma}
{\em Proof:}
We first treat the case where $S'$ is artinian local with algebraically closed
residue field. Let $Z\rightarrow G$ and $R\rightarrow Z\times _GZ$ be the
vertical surjections given by the fact that $G$ is $P4$. There is an etale
neighborhood $U\rightarrow Z\times _SZ$ of the diagonal such that $U$ is
isomorphic to an etale neighborhood of the zero section in the total scheme $TZ$
(and this isomorphism is compatible with the first projection to $Z$). This is
seen as in the argument above.  Furthermore, as above we may assume that the
fibers of the first projection $U\rightarrow Z$ are connected (over any artinian
scheme).  Then for any artinian scheme $S'\rightarrow U$, the two elements of
$G(S')$ obtained from the two projections $U\rightarrow Z$ are the same, by the
hypothesis that $G$ is finite.  (To see this, compare $(a,b): S'\rightarrow U$
with $(a,a): S\rightarrow U$; they are in the same fiber over $a$, and this
fiber is connected, so they have to have the same image in $G(S')$.)  Thus, any
artinian subscheme of $U$ lifts into $R$.  This implies that there is (locally
in the etale topology) a lifting $U\rightarrow R$.  Let $V\subset Z\times _SZ$
be the image of $U$.  It is a Zariski neighborhood of the diagonal, and locally
there is a lifting from $V$ into $R$.  Let $F\subset Z\times _SZ$ be the
reduced closed subscheme corresponding to the closed subset which is the
complement of $V$.  Suppose $Y\rightarrow S$ is an artinian local scheme (with
acrf).  If $(\alpha _1,\ldots , \alpha _n)$ is an $n$-tuple of distinct
points of $G(Y)$, then there is a  lifting $(a_1,\ldots , a_n) \in Z\times _S
\ldots \times _S Z(Y)$ such that for any $i,j$ we have that $(a_i, a_j):
Y\rightarrow Z\times _SZ$ is not contained in $V$. In particular, the
reduced point $(a_i,a_j)^{\rm red}$ is contained in $F$. Thus the reduced
point
$(a_1,\ldots , a_n)^{\rm red}$ is contained in the closed subscheme
$$
F^{(n)}:= \bigcap _{i,j} pr_{ij}^{-1}(F)\subset X\times _S\ldots \times _SX.
$$
We claim that there is an $n$ such that $F^{(n)}$ is empty.  For any
$(x_1,\ldots , x_k)\in F^{(k)}$, let
$$
\Phi (x_1,\ldots , x_k):= \{ y\in X, \;\; \pi (y)= \pi (x_i)\in S,\;\;
(y,x_1,\ldots , x_k)\in F^{(k+1)}\} .
$$
Note that these are closed subschemes of $X$ with strict inclusions
$$
\Phi (x_1,\ldots , x_k) \subset \Phi (x_1,\ldots , x_{k-1}).
$$
Furthermore, $\Phi (x_1,\ldots , x_k)$ varies algebraically with $(x_1,\ldots ,
x_k)$.

Let $d=dim (X)$ and let $\Lambda = \nn ^d$ with the lexicographic ordering
giving the most importance to the $d$th coordinate.  For any algebraic set
$Y$ of dimension $\leq d$, let $\lambda (Y)= (\lambda _1, \ldots , \lambda _d)$
be defined by setting $\lambda _d$ equal to the number of irreducible components
of dimension $d$.  Note that if $Y'\subset Y$ is a strict inclusion of a
closed subset then $\lambda (Y')< \lambda (Y)$.  Let $\Lambda ^{(k)}$ be the
finite set of all  $\lambda (\Phi (x_1,\ldots , x_k))$ for $(x_1,\ldots ,
x_k)\in F^{(k)}$ (it is finite because
$\Phi (x_1,\ldots , x_k)$ varies algebraically with $(x_1,\ldots , x_k)$).
Introduce an order relation on subsets $\Sigma \subset \Lambda$ by saying
$$
\Sigma < \Sigma \; \Leftrightarrow \forall \sigma \in \Sigma ,\, \exists \sigma
' \in \Sigma ',\;\; \sigma < \sigma ' .
$$
Then the sequence $\Lambda ^{(k)}$ is a sequence of finite subsets which is
strictly decreasing for this order relation.  We claim that this implies (by
combinatorics) that one of the $\Lambda ^{(k)}$ is empty.  To see this, assume
that the combinatorial claim is true for $d-1$.  We will show that the set of
upper bounds for $\lambda _d$ on $\Lambda ^{(k)}$ doesn't stabilize.  If it
were to stabilize after $k_0$ at a certain $y$, then for $k\geq k_0$ we could
let
$A^k\subset \nn ^{(d-1)}$  be the subset of elements $(a_1,\ldots, a_{d-1})$
such that   $(a_1,\ldots , a_{d-1},y)\in \Lambda ^{(k)}$.  We obtain a strictly
decreasing sequence of subsets for the case of $d-1$, so it is eventually
empty, meaning that in fact the upper bound for $\lambda _d$ didn't stabilize.
A decreasing sequence which doesn't stabilize can't exist, so eventually there
is
no upper bound, in other words $\Lambda ^{(k)}$ becomes empty. This gives the
claim.

Since one of the $\Lambda ^{(k)}$ is empty, one of the $F^{(k)}$ is empty.  Let
$N$ be chosen so that $F^{(N)}$ is empty (and consequently $F^{(k)}$ is empty
for
$k\geq N$). Then
by the above argument, if $(\alpha _1,\ldots , \alpha _n)$ is an $n$-tuple of
distinct points of $G(Y)$, we must have $n<N$. This gives the theorem in the
case of an artinian local $Y$.

Now suppose $A$ is a henselian local ring and $S'= Spec (A)$.  Let $S'_n:=Spec
(A/{\bf m}_A^n)$. In the inverse system $\lim _{\leftarrow} G(S'_n)$ we have
that all of the $G(S'_n)$ have cardinality bounded by $N$.  In particular, the
cardinality of the inverse limit is bounded by $N$.  Now suppose that there are
$N+1$ distinct points $y_i$ in $G(S')$.  Two of the points go to the same point
in $\lim _{\leftarrow} G(S'_n)$, which means that for two of the points, the
liftings $z_i,z_j\in Z(S')$ give a point $(z_i,z_j)$ in $Z\times _SZ$ which
lifts, over any $S'_n$, into $R$. By strong artin approximation (check here
!!!), the point $(z_i,z_j)$ must lift into $R$ so the two points in
$G(S')$ are equal, a contradiction.
This completes the proof of the lemma.
\eop

\begin{corollary}
\mylabel{I.1.q}
If $G$ is a presentable group sheaf, then $G$ is finite if and only if $G(S')$
is finite for any henselian (resp. artinian) $S'$.
\end{corollary}
{\em Proof:}
The lemma provides one direction.  For the other, note that if $G(S)$ is finite
for artinian $S'$ then $G^0(S')=\{ 1\}$.  By unicity in the characterization of
$G^0$ we get $G^0= \{ 1\}$.
\eop

\numero{Local study of presentable subgroups}

In this section we show that if $H\subset G$ is a presentable subgroup of a
presentable group $G$ then locally at the identity, in an appropriate sense,
$H$ is defined by the vanishing of a section of a vector sheaf.
This is a generalisation of the basic result that a subgroup of an algebraic
group is smooth, and hence a local complete intersection---cut out by a
section of its normal bundle.  We obtain this result only in a ``neighborhood of
the identity'', or more precisely upon pullback by a vertical morphism
$X\rightarrow G$ such that $X$ admits a lift of the identity.  If $Y$ is a
scheme with morphism $Y\rightarrow G$ such that $P\in Y$ maps to the identity
section in $G$, then there will be an etale neighborhood of $P\in Y$ lifting
to $X$ (which is why we can think of $X$ as a neighborhood of the identity).

This result will be used in a future study of de Rham cohomology (results
announced in \cite{kobe}).  There, it will be important to have a structure
theory for presentable subgroups because of the general principle that if $G$
is a presentable group sheaf then $G/Z(G)\subset Aut (\Ll )$ where $\Ll = Lie
(G)$ is the Lie algebra vector sheaf of $G$ (see \S 9 below). A good
understanding of the structure of  presentable subgroups will allow us to
reduce to looking at de Rham cohomology with coefficients in $Aut (\Ll )$ for
$\Ll$ a vector sheaf, and here we have a more concrete hold on what happens.

\begin{theorem}
\mylabel{D.1}
Suppose $G$ is a connected presentable group sheaf over $S$, and suppose
$H\subset G$ is a  presentable subgroup sheaf.  Suppose that $X_1\rightarrow
G$ is a vertical morphism with lift of the identity
section $e:S\rightarrow X_1$. Suppose $P\in S$. Then there is an etale
neighborhood $X\rightarrow X_1$ of $e(P)$ with a lift of the identity $e:
S\rightarrow X$ (possibly after localizing in the etale topology of $S$ here)
and an etale morphism $\rho : X\rightarrow TX_e$  sending $e$ to the zero
section, such that
$$
X\times _GH =  \rho ^{-1}(TX_e \times _{TG_e} TH_e).
$$
In particular, there is a vector sheaf $V$ over $S$ and a section $\sigma :
X\rightarrow V$ such that $X\times _GH= \sigma ^{-1}(0)$.
\end{theorem}
{\em Proof:}
Let $X_1\rightarrow G$ be a surjective vertical morphism with $X_1(S')$
connected for all artinian $S'$ (with $X_1$ a scheme of finite type).

Put $Y_1:= X_1\times _GH$. It is a subsheaf of $X_1$.

We can choose a vertical surjection $Z_1\rightarrow Y_1$ (with $Z_1$ a
scheme of finite type over $S$) together with a lift of $(e,e)$ also denoted by
$e$.  Note that the morphism  $Z_1\rightarrow H$
is also vertical (using the composition property of vertical morphisms and
the fact that the morphism $Y_1\rightarrow H$ is vertical by the pullback
property).

There is an etale neighborhood of $(e,e)\in X_1\times _SX_1$ denoted by
$U_1\rightarrow X_1\times _SX_1$ together with a lifting $\psi : U_1\rightarrow
X_1$ of the multiplication in $G$, such that $\psi $ restricted to the inverse
images of $\{ e\} \times _SX_1$ or  $X_1\times _S\{ e\} $ are the identity.
We obtain a morphism
$$
U_1\times _{X_1\times _SX_1}(Y_1\times _SY_1)\rightarrow Y_1
$$
compatible with the multiplication in $H$ and again having the property that
the restrictions to the inverse images of the two ``coordinate axes'' are the
identity.  Now pull back our multiplication to
$$
U_1\times _{X_1\times _SX_1}(Z_1\times _SZ_1)
$$
and note that $Z_1\rightarrow Y_1$ being vertical,  there is an etale
neighborhood of the identity section (all of this is local on $S$!)
$$
V_1\rightarrow U_1\times _{X_1\times _SX_1}(Z_1\times _SZ_1)
$$
(which we can consider just as an etale neighborhood $V_1\rightarrow Z_1\times
_SZ_1$) and a good lift of our multiplication
$$
V_1\rightarrow Z_1
$$
restricting to the identity on the inverse images of the ``coordinate
axes''.
We obtain in this way morphisms on the etale germs
$$
{\bf 2}_{Z_1}: (Z_1,e)\rightarrow (Z_1,e)
$$
and
$$
{\bf 2}_{X_1}: (X_1,e)\rightarrow (X_1,e)
$$
compatible with the morphism $Z_1\rightarrow X_1$.
These morphisms induce multiplication by $2$ on the tangent vector schemes.
There are unique analytic isomorphisms of complex analytic germs
$$
(X_1,e)^{\rm an}\cong (T(X_1)_e,0)^{\rm an}
$$
and
$$
(Z_1,e)^{\rm an}\cong (T(Z_1)_e,0)^{\rm an}
$$
transforming the automorphisms ${\bf 2}$ into multiplication by $2$ and
inducing the identity on tangent spaces at the identity section. (To see
uniqueness, note that over artinian bases these are germs of vector spaces, and
any germ of automorphism $f$ of a vector space, such that $f(2x)=2f(x)$, is
linear; hence fixing it at the identity fixes it.)

By uniqueness, these isomorphisms are compatible with the morphism
$Z_1\rightarrow X_1$.

On the formal level, we have an etale morphism of formal germs
$$
\hat{\varphi}:\widehat{T(X_1)_e}\rightarrow X_1
$$
such that $\widehat{T(Z_1)_e}$ maps into $Y_1$.  The {\em first claim} is that,
in fact, this gives a map
$$
Spec (\widehat{\Oo} _{T(X_1)_e,e(S)}) \rightarrow X_1
$$
such that
$$
T(Z_1)_e \times _{T(X_1)_e}Spec (\widehat{\Oo} _{T(X_1)_e,e(S)})
$$
maps into $Y_1$.
We can now apply Artin approximation to find an etale neighborhood
$W_1\rightarrow T(X_1)_e$ of the identity section (of course locally on $S$)
together with a morphism $W_1\rightarrow X_1$ inducing the identity on tangent
vector schemes at the identity section, and sending
$$
T(Z_1)_e\times _{T(X_1)_e}W_1\rightarrow Y_1.
$$
We can suppose that the morphism $W_1\rightarrow X_1$ is etale. In particular
the morphism $W_1\rightarrow G$ is vertical. We obtain two subsheaves
$$
im (T(Z_1)_e\times _{T(X_1)_e}W_1\stackrel{pr_2}{\rightarrow} W_1)
\subset
W_1\times _{X_1} Y_1 \subset W_1.
$$
They have the same tangent subsheaves at the identity.

Our {\em main claim} is that by taking an open subset of $W_1$ (still a
neighborhood of $e(P)$ for a given basepoint $P\in S$) we can assume that these
two subsheaves are equal.

The first subsheaf is given by the vanishing of the morphism
$$
W_1\rightarrow T(X_1)_e /T(Z_1)_e = T(G)/T(H),
$$
while the second subsheaf is equal to $W_1\times _GH$.  Setting $X=W_1$ we
obtain the result of the theorem.
We just have to prove the {\em first claim} and the {\em main claim}.

{\em Proof of the first claim:}
By the sheaf condition and the finite type condition B1 and B2 for $Y_1$, it
suffices to prove that for any artinian $S'$, we have
$$
T(Z_1)_e \times _{T(X_1)_e}Spec (\widehat{\Oo} _{T(X_1)_e,e(S)})(S')
$$
mapping into $Y_1(S')$.  That is to say, we have to prove that for any point
$S'\rightarrow T(Z_1)_e$ mapping to a point of $T(X_1)_e$ located near the
origin (that is to say factoring through
$Spec (\widehat{\Oo} _{T(X_1)_e,e(S)})$), this point maps into $Y_1(S')$.

We change to an algebraic notation.  We can suppose that $S=Spec(A)$,
$T(X_1)_e=Spec(B)$ and $T(Z_1)_e=Spec(C)$.  Further we can suppose that
$S'=Spec
(K)$ with $K$ artinian (although not necessarily of finite type). We have
$C\rightarrow K$. Since $T(Z_1)_e$ is a vector scheme we have a map
$C\rightarrow C[t]$ corresponding to multiplication by $t$ (and compatible with
the same map on $B$). Let $\hat{B}$ denote the completion of $B$ around the
zero section (which corresponds to an ideal ${\bf b}\subset B$).  We are
provided with a factorisation $B\rightarrow \hat{B}\rightarrow K$.

We can assume that $K$ is of finite type over $\hat{B}$, and in particular
that $K$ is the total fraction ring of a subring $R\subset K$ such that $R$
is finite over $\hat{B}$.  Let  ${\bf
r}\subset R$ denote the ideal corresponding to ${\bf b}\subset B$ (note that
$R$ is complete with respect to ${\bf r}$). Let $K\{ t\} \subset K[[ t]]$ denote
the set of formal series of the form $\sum a_it^i$ such that there exists $\eta
\in R$ such that $\eta a_i \in {\bf r}^i$.  With the same notations for $B$,
multiplication by $t$ provides a map  $\hat{B}\rightarrow B\{ t\}$ compatible
with the map $B\rightarrow B[t]$, hence we get a map $\hat{B}  \rightarrow K\{
t\}$.  On the other hand we get a map $C\rightarrow C[t]\rightarrow
K[t]\rightarrow K\{ t\}$.  Putting these together we get a map
$$
\hat{B}\otimes _BC \rightarrow K\{ t\}
$$
corresponding to multiplication by $t$.  There is an evaluation at $t=1$ which
is a map $K\{ t\} \rightarrow K$ (this summability of the formal series comes
from the definition of $K\{ t\}$ and the completeness of $R$),
and the above map
is compatible with this and with the map $\hat{B}\otimes _BC\rightarrow K$ given
at the start.  All in all we obtain a map
$$
Spec (K\{ t\}) \rightarrow
T(Z_1)_e \times _{T(X_1)_e}Spec (\widehat{\Oo} _{T(X_1)_e,e(S)})
$$
which induces on the subscheme $Spec (K)\rightarrow Spec (K\{ t\})$
(evaluation at $t=1$) the original inclusion.  Now compose with the projection
into $G$.  We obtain a morphism
$$
Spec (K\{ t\} )\rightarrow G
$$
which sends $Spec (K[[t ]])$ into $H$ (this comes from the condition that
$\widehat{T(Z_1)}$ maps into $Y_1$ together with B1 and B2 for $Y_1$ or $H$)
and we would like to show that it sends $Spec (K)$ (at $t=1$) into $H$.
It suffices to show that $Spec (K\{ t\} )\rightarrow H$.

By Noether normalization there is a morphism $R'\rightarrow R$ such that
$R'$ is integral and $R$ is finite over $R'$.  Let $K'$ be the total fraction
ring of $R'$: it is a field, and $K$ is finite over $K'$.  There is an ideal
${\bf r'}\subset R'$ which induces ${\bf r}$, and $K\{ t\}$ is finite over the
ring $K'\{ t\}$ defined in the same way as above with respect to this ideal.
Let $G'$ and $H'$ denote the direct images to $Spec (K')$ of the groups $G$
and $H$ pulled back to $K$.  We have that $H'$ is a presentable subgroup of
the presentable group $G'$ (Lemma \ref{I.1.i}),
but since $K'$ is a field, $H'\subset G'$ is a closed subgroup of the
algebraic group $G'$ over $K'$.  Since $K$ is finite over $K'$ we have
$$
K\{ t\} = K'\{ t\} \otimes _{K'}K,
$$
whence our point  $Spec (K\{ t\} )\rightarrow G$ gives a point
$Spec (K'\{ t\} )\rightarrow G'$ sending $Spec (K'[[ t]])$ into $H'$.  Now
since $H'$ is a closed subgroup of $G'$ both of which are algebraic groups (of
finite type) over $K'$, we get that $Spec (K' \{ t\} )\rightarrow H'$, meaning
that $Spec (K\{ t\} )\rightarrow H$.
This completes the proof of the first claim.
\eop

{\em Proof of the main claim:}
Suppose that the main claim is not true.  Note that there is a scheme of finite
type surjecting to $W_1\times _{X_1}Y_1$.  The falsity of the main claim means
that the morphism from this scheme to $T(X_1)_e/T(Z_1)_e$ is nonzero on any
subset of the form pullback of an open subset of $W_1$ containing $P$.  In
particular we can find a (possibly nonreduced) curve inside this scheme, such
that the section pulls back to something nonzero on the generic (artinian)
point, but such that the image of the curve in $W_1$ contains $P$ in its
closure.  We get an $S$-scheme $S'$ with reduced scheme equal to a curve, and a
morphism $\psi :S'\rightarrow W_1\times _GH$ such that the projection into
$T(X_1)_e/T(Z_1)_e$ is nontrivial at the generic point of $S'$, such that $P$
is in the closure of the image of $S'$.

Let $\overline{S}'$ be a closure of $S'$ relative to $W_1$ obtained by adding
one point over $P$.  Call this point $P'$.  Then for any $n$ there is an etale
neighborhood of $P\in W_1$ on which the squaring map $n$-times is defined.
We obtain an etale $\overline{S}'_n\rightarrow \overline{S}'$ on which the
squaring map $n$-times is defined.  We may assume that $\overline{S}'_n$
consists of an etale morphism $S'_n\rightarrow S'$, union one point $P'_n$ over
$P'$. Denote by $\psi _n:
\overline{S}'_n\rightarrow W_1$ the result of the squaring operation
iterated $n$
times.  There is an analytic isomorphism of a neighborhood of
$P'_n$ in $\overline{S}'_n$ with a neighborhood of $P'$ in $S'$, and an
analytic trivialization of a neighborhood of $P$ in $W_1$ (isomorphism with the
tangent vector scheme) such that $\psi _n= 2^n\psi$ as analytic germs
around the
point $P'_n$.

{\em Step 1.} There is an $n_0$ such that for any $n\geq n_0$, the projection of
$S'_n$ into $T(X_1)_e/T(Z_1)_e$ is nontrivial at the generic point of $S'_n$.
In particular for any $m$ the projection of
$S'_{mn_0}$ into $T(X_1)_e/T(Z_1)_e$ is nontrivial at the generic point of
$S'_{mn_0}$.

Let $v: W_1\rightarrow T(X_1)_e/T(Z_1)_e$ denote our section.  With respect
to our analytic trivialization of $W_1$ where the squaring map becomes
multiplication by $2$, can  take a Taylor expansion for $v$ around the identity
section of $W_1$,
$$
v= v_1 + v_2+ v_3 + \ldots  + v_{i-1} + w_i,
$$
with $v_j(2x)= 2^jv(x)$ and $w_i$ vanishes to order $i$ along $e$;
this notion can be defined by considering $w_i$ as a section  of a coherent
sheaf $\Ff$ which contains $T(X_1)_e/T(Z_1)_e$.  By hypothesis the restriction
of $v$ to $S'$ is nonzero at the generic point of $S'$. Let $\Gg _{S'}$ be the
quotient of $\Ff |_{S'}$ by the ``torsion'' subsheaf (i.e. the subsheaf
of sections supported in dimension zero).  That a
section is nonzero at the generic point means that its projection into $\Gg
_{S'}$ is nonzero.  We may choose $i$ big enough so that $v$ is nonzero in $\Gg
_{S'}$ modulo the image of sections which vanish to order $i$ along $e$.
Let
$\overline{v}_j$ denote the projection of $v_j$ into the space of sections of
$\Gg_{S'}$ modulo the image of the sections vanishing to order $i$.
At least one of the $\overline{v}_j$ is nonzero.
Now notice that the projection of $v(2^nx)$ is equal to
$$
\overline{v(2^nx)} = 2^n\overline{v}_1(x) + 2^{2n}\overline{v}_2(x)
+ \ldots + 2^{(i-1)n}\overline{v}_{i-1}(x).
$$
A little $2$-adic argument shows that there is $n_0$ such that for $n\geq n_0$
this quantity must be nonzero.  We obtain that
$\overline{v(2^nx)}\neq 0$ and hence that $v(2^nx)=v(\psi_nx)$ is nonzero at
the generic point of $S'_n$, as claimed for Step 1.

{\em Step 2.}   The Zariski closure of the union of the images of the $\psi
_{mn_0}$ contains the zero-section.  To prove this, note that in the formal
completion at $P$, the union of the closures of the $S'_{mn_0}$ is a subset
stable
under multiplication by $2^{n_0}$, hence its Zariski closure is stable under
(fiberwise) multiplication by $2$, hence it is fiberwise homogeneous and thus
contains the zero-section.  The completion of the Zariski closure contains the
Zariski closure of the intersection with the completion, so the zero-section is
in the closure.

{\em Step 3.}  Over the generic  point of $S$, the zero section is  in the
Zariski closure of the $S'_{mn_0}$.  Otherwise we would obtain a function
nonvanishing on the zero section and vanishing on the $S'_{mn_0}$; clearing
denominators this function can be assumed defined over $S$ rather than the
generic point of $S$, and since (we may assume) the $S'_{mn_0}$ are all schemes
of
pure dimension $1$ dominating $S$, this function defined over $S$ which
vanishes generically on the $S'_{mn_0}$, must vanish identically on the
$S'_{mn_0}$.
This would contradict the fact that the zero section is in the Zariski closure
globally over $S$.

{\em End of proof of claim:}  Now we work over the generic geometric artinian
point of $S$.  Change notations now to suppose that $S$ is artinian and $S'=S$;
we note the schemes $S'_{mn_0}$ by $S_{mn_0}$ (they are all isomorphic to $S$)
with $S'=:S_1$.   We have points $S_{mn_0}\rightarrow W_1\times _GH$ all mapping
to something nonzero in $T(X_1)_e/T(Z_1)_e$.

Note, as a bit of a detour, that the connected component of the identity in
$W_1\times _GH (S)$, must map to zero in $T(X_1)_e/T(Z_1)_e(S)$.  This is
because $T(X_1)_e/T(Z_1)_e(S)= T(X_1)_e(S)/T(Y_1)_e(S) = T(G)_e(S)/T(H)_e(S)$,
whereas verticality of $X_1\rightarrow G$ implies that
$X_1(S)\rightarrow G(S)$ is smooth.  In particular
$W_1\times _GH (S)$ is a smooth local complete intersection so a morphism from
$W_1(S)$ to the normal space $T(G)_e(S)/T(H)_e(S)$ of $W_1\times _GH (S)$,
with zero set contained in the complete intersection, must have
zero set which is a union of connected components of $W_1\times _GH (S)$.
Containing the identity, it contains the connected component of
the identity.

 In particular, our points $S_{mn_0}\rightarrow W_1\times _GH$ from before are
never in the connected component of $W_1\times _GH (S)$ which contains the
identity.  On the other hand, these  points all lift to $Z_2\rightarrow W_1$
(a scheme of finite type surjecting vertically to  $W_1\times _GH$).  Let
$Z_2(S)'$ denote the union of components of $Z_2(S)$ which contain liftings of
our points $S_{mn_0}\rightarrow W_1$.  We have a morphism $Z_2(S)'\rightarrow
W_1(Spec (k))$ whose image is a constructible set.  But the image contains all
of the points where the $S_{mn_0}$ are located, so the image must contain a
generic point of any irreducible component of the Zariski closure of the
$S_{mn_0}$. In particular, there is a component of $Z_2(S)'$ which maps to
something in $W_1(Spec (k))$ containing the identity in its closure.

Let $W_1(S)_e$ denote the inverse image of $e\in W_1(Spec (k))$ in $W_1(S)$.
Let $N\subset G(S)$ denote the image of $W_1(S)_e\rightarrow G(S)$.  We claim:
that $N$ is a unipotent subgroup of $G(S)$, and that the morphism
$W_1(S)_e\rightarrow N$ is a fibration with connected fibers.

Assume this claim for the moment.
The image of $W_1(S)\rightarrow W_1(Spec (k))$ is a closed subvariety
$R\subset W_1(Spec (k))$ (this can be seen since $W_1$ is etale over the vector
scheme $TX_1$).  We have a morphism $R\rightarrow G(S)/N$. On the other hand,
the above morphism $Z_1(S)' \rightarrow W_1(Spec (k))$ factors through a
morphism $Z_1(S)'\rightarrow R$, and the image of this map contains $1\in R$
in its closure.

The morphism $W_1(S)\rightarrow R$ is a fibration with fiber $W_1(S)_e$ in the
etale topology.
It suffices to prove that $TX_1 (S) \rightarrow TX_1(Spec (k))$ is a fibration
over its image, since locally in the etale topology $W_1$ is isomorphic to
$TX_1$.
In fact if $V$ is any vector scheme then $V(S)$ and $V(Spec (k))$ are vector
spaces so the morphism $V(S)\rightarrow V(Spec (k))$ is a fibration over its
image, with fiber the inverse image of the origin.

We now show that the morphism
$$
W_1\times _GH(S) \rightarrow R \times _{G(S)/N} H(S)
$$
is a fibration in the etale topology with fiber the kernel of
$W_1(S)_e\rightarrow N$.  Locally on $R$ we can choose a lifting $\lambda : R
\rightarrow W_1(S)$ and then we have a morphism
$$
R\times _{G(S)/N}H(S)\rightarrow N
$$
given by $(r,h)\mapsto h^{-1}im(\lambda (r))$.  We claim that (locally over
$R$)
$$
W_1(S)\times _{G(S)}H(S) = W_1(S)_e \times _N (R\times _{G(S)/N}H(S)).
$$
The morphism from right to left associates to the point
$(a,r,h)$ the point $(i(a)\ast \lambda (r) , h)$ where $i: W_1\rightarrow
W_1$ is
an etale-locally defined morphism covering the inverse.   This shows
that the morphism at the start of the paragraph is a fibration.

Suppose $A$ is an algebraic group with connected algebraic subgroups $B\subset
A$ and  $N\subset A$.  Then the morphism
$$
B / (B\cap N) \rightarrow A/N
$$
is proper over an open neighborhood of the class of the identity in $A/N$.
To prove this, proceed as follows.
Let $I\subset A/N$ denote the image.  Let $Z\subset A/N$ denote the subset of
points over which the map in question is not proper.  This can be constructed
as follows.  Let $X:= B/(B\cap N)$, and let $\overline{X}$ be a
relative completion with proper morphism $\overline{X}\rightarrow A/N$; and
suppose that $X\subset \overline{X}$ is open and dense. Then $Z$ is the image
of $\overline{X}-X$.  Since the map $X\rightarrow A/N$ is injective, we have
that the dimension of the image $Z$ is strictly less than the dimension of
the image $I$ of $X$.  In particular, there is a point $y\in I$ such that
the morphism in question is proper over a neighborhood $U$ of $y$.  But since
$B$ acts on $X$ and compatibly on $A$ (by left
multiplication) the morphism in question is proper over any translate of the
form $bU$.  Setting $b\in B$ equal to the inverse of a representative in for
$y$ we obtain a neighborhood $bU$ of the identity over which the map is proper.

Note that by the above claim that we are accepting for now, the fiber of
the fibration $W_1(S)_e\rightarrow N$ is connected.  On the other hand, by
the previous paragraph the map $H(S)\rightarrow G(S)/N$ induces a map
$H(S)/(H(S)\cap N)\rightarrow G(S)/N$ which is proper over a neighborhood of
the class of the identity. Let $I\subset G(S)/N$ be the image of this map. It is
a locally closed subset and the subset topology coincides with the topology of
the base of the fibration, at least near the identity. Note that $H(S)$ is
fibered over $I$ with fibers $H(S)\cap N$ which are connected because $N$ and
hence $H(S)\cap N$ are unipotent groups (unipotent groups are always connected).
Finally we have the following situation:
$$
W_1\times _GH(S)\rightarrow R\times _{G(S)/N} I
$$
is a fibration with connected fiber, whereas $I\subset G(S)/N$ is a locally
closed subset.  Since an etale fibration is an open map, the image of the
connected component of the identity in $W_1\times _GH(S)$ is an open
neighborhood of the identity in $R\times _{G(S)/N}I$.  Since the fibers of the
fibration are connected, the image of the complement of the identity component
is the complement of the image of the identity component.  In particular,
there is an open neighborhood of the identity in $R\times _{G(S)/N}I\subset R$
(and hence an open neighborhood of the identity in $R$) whose inverse image
doesn't meet any other connected component of $W_1\times _GH(S)$.  Finally,
since $R$ is a closed subset of $W_1(Spec(k))$, we obtain an open neighborhood
of the identity in $W_1(Spec (k))$ whose inverse image in $W_1\times _GH(S)$
is contained in the connected component of the identity.  This is a
contradiction to our earlier situation where $Z_2(S)' \rightarrow W_1(Spec
(k))$ has image a constructible set with the identity in its closure.  This
completes the proof modulo the following part.

We have to show that $N$ is a unipotent subgroup of $G(S)$ and that the
morphism $W_1(S)_e\rightarrow N$ is a fibration with connected fibers.

Write $S=Spec (A)$ with $A$ artinian, and choose a sequence of ideals
$I_j\subset A$ for example $I_j= {\bf m}^j$.  Let
$$
W_1(S)_j
$$
be the set of points of $W_1(S)$ which restrict to the identity on $Spec
(A/I_j)$.  In particular $W_1(S)_1=W_1(S)_e$.  Choose a good lift of the
multiplication
$$
W_1\times W_1 \rightarrow W_1
$$
in a formal neighborhood of our point $P$.  We obtain
$$
\ast : W_1(S)_e\times W_1 (S)_e\rightarrow W_1(S)_e.
$$
This operation is not a group, however we have the following property.
$$
\ast : W_1(S)_j\times W_1 (S)_j\rightarrow W_1(S)_j.
$$
Next, note that the fact that $W_1$ is etale over a vector scheme gives
another operation which we denote
$$
+ : W_1(S)_e\times W_1 (S)_e\rightarrow W_1(S)_e
$$
which is an abelian group structure.
We can write
$$
a*b = a+b + F(a,b)
$$
where
$$
F: W_1(S)_j\times W_1 (S)_j\rightarrow W_1(S)_{j+1}.
$$
This is because there is a unique good operation on the set of
elements of $W_1(Spec (A/I_{j+1})$ which restrict to the identity in $W_1(Spec
(A/I_j)$.  (Also note that the $+$-quotient $W_1(S)_j/W_1(S)_{j+1}$ injects into
this subset of $W_1(Spec(A/I_{j+1})$).

Because of this formula, we can define the quotient $W_1(S)_j/W_1(S)_{j+1}$
with respect to the operation $\ast$ and it is the same as the quotient with
respect to $+$. In particular note that the morphism
$$
W_1(S)_j\rightarrow W_1(S)_j/W_1(S)_{j+1}
$$
is a fibration with fiber a vector space.

The morphism $W_1(S)_e\rightarrow G(S)$ is compatible with the operation
$\ast$.  Let $N_j$ denote the image of $W_1(S)_j$ in $G(S)$ (in particular
$N_1=N$).  The $N_j$ are constructible sets and subgroups so they are
algebraic subgroups of $G(S)$.  We obtain a surjective morphism
$$
W_1(S)_j /W_1(S)_{j+1} \rightarrow N_j/N_{j+1},
$$
but from the previous formula the operation on the left is a unipotent
algebraic group, this shows that $N_j/N_{j+1}$ is a unipotent group, and since
extensions of unipotent groups are unipotent (and $N_j=\{ 1\}$ for $j$ large),
$N=N_1$ is unipotent.

We claim that
$W_1(S)_e\rightarrow N$ is smooth.  Suppose $R'\subset R$ is an inclusion
of artinian schemes over $Spec(k)$.  Look at the map
$$
W_1(S)_e(R)\rightarrow W_1(S)_e(R')\times _{N(R')}N(R).
$$
Suppose that we have a map $S\times R\rightarrow
G$ and a lifting over $S\times R'$ to $W_1$, sending $Spec(k)\times R$ to $e$.
We would like to find a lifting over $S\times R$ sending $Spec(k)\times R$ to
$e$.  We can do this whenever $R'$ is a union of $R_i$ and we have commuting
retracts from $R'$ to $R_i$, just apply the verticality property to $R\times
S$ with retracts to $R_i\times S$ and $R\times Spec( k)$.  This proves that
the morphism $W_1(S)_e\rightarrow G(S)$ is vertical with respect to
$Spec(k)$.  It is then immediate that the morphism $W_1(S)_e\rightarrow N$
is surjective (since $N$ is the image of the previous map).  Note that $N$ is
presentable, so by Lemma \ref{smooth}, the morphism
$W_1(S)_e\rightarrow N$ is smooth.

Let $K\subset W_1(S)_e$ be the inverse image of $1\in N$.
Then for any two points $a,b$ with $f(a)=f(b)$ there is a unique element
$k\in K$ such that $b=k\ast a$ (the existence and uniqueness of such an
element $k\in W_1(S)_e$ can be seen using the above grading and the expression
for $\ast$, and then it is immediate that $k\in K$ from compatibility of
$f$ with $\ast$). Any point $n\in N$ has an etale neighborhood $n\in
U\stackrel{p}{\rightarrow} N$ with a section $\sigma :U\rightarrow W_1(S)_e$.
Then we obtain a morphism
$$
K\times U \rightarrow W_1(S)_e\times _N U
$$
obtained by sending $(k,u)$ to $(k\ast \sigma (u), p(u))$.  This is an
isomorphism on the level of points by the above property for $K$, and both
sides are smooth, so it is an isomorphism.  This proves that
$W_1(S)_e\rightarrow N$ is a fibration in the etale topology.

Finally, to
show that the fibers are connected it suffices to show that $K$ is connected.
But since $W_1(S)_e$ and $N$ are vector spaces and the morphism $f$ is a
fibration in the etale topology, the associated analytic morphism is a
fibration in the usual topology, so the fiber is contractible.
\eop

\numero{The Lie algebra sheaf}

\begin{theorem}
\mylabel{lmn}
If $\Gg$ is a presentable group sheaf, and if
we set $Lie (\Gg) := T(\Gg )_1$, then there is a unique bilinear form
(Lie bracket)
$$
[\cdot , \cdot
] : Lie (\Gg )\times Lie (\Gg )\rightarrow Lie (\Gg )
$$
which, over artinian base schemes, reduces to the  usual Lie bracket.
\end{theorem}
{\em Proof:}
A section of $Lie (\Gg )$ over $S'= Spec (A)$ is a morphism $Spec (A[\epsilon
])\rightarrow  \Gg$ sending $Spec (A)$ to the identity section (in our notation
here $\epsilon$ denotes an element with $\epsilon ^2=0$).
Given two
such morphisms which we denote $\alpha$ and $\beta$ we obtain
$$
\alpha p_1, \, \beta p_2 : Spec (A[\epsilon , \epsilon '])\rightarrow \Gg
$$
(where also $(\epsilon ')^2=0$) and we can form the morphism
$$
\gamma := \alpha p_1 \cdot \beta p_2 \cdot (\alpha m p_1) \cdot (\beta m
p_2)
: Spec (A[\epsilon , \epsilon '])\rightarrow \Gg
$$
where the $g\cdot h$ denotes composition in $\Gg$, and where $m= A[\epsilon
]\rightarrow A[\epsilon ]$ is the involution sending $\epsilon$ to $-\epsilon$.
The morphism $\gamma$ restricts to the identity on $Spec (A[\epsilon , \epsilon
']/(\epsilon \epsilon '))$.  Let  $$
q: Spec (A[\epsilon , \epsilon '])\rightarrow Spec (A[\delta])
$$
denote the morphism sending $\delta$ to $\epsilon \epsilon '$ (here again
$\delta
^2=0$).

Our first claim is that if the morphism $\gamma$ factors as $\gamma = \varphi
\circ q$ then $\varphi $ is unique.  To see this suppose that $\phi$ and
$\varphi$ were two morphisms $Spec (A[\delta])\rightarrow \Gg$ with $\phi \circ
q = \varphi \circ q$. Let $X\rightarrow \Gg$ and $R\rightarrow X\times _{\Gg}X$
be the morphisms in a presentation of $\Gg$, with a chosen lift of the identity
section into $X$.   Choose liftings $\tilde{\varphi}$ and $\tilde{\phi}$
from $Spec
(A[\delta ])$ into $X$ sending $Spec (A)$ to the identity section of $X$
(here we
may have to localize on $S'=Spec (A)$ in the etale topology---but henceforth
ignore this point, much as we have already ignored it in lifting the identity
section into $X$\ldots ).  The fact that the compositions with $q$ are the same
means that the pair $(\tilde{\varphi}\circ q, \tilde{\phi}\circ q)$
defines a point which we denote
$$
\eta :Spec (A[\epsilon , \epsilon '])\rightarrow
X\times _{\Gg}X.
$$
 Note that $Spec (A[\epsilon , \epsilon ']/(\epsilon \epsilon
'))$ projects by $q$ to $Spec (A)\subset Spec (A[\delta ])$ and both
$\tilde{\varphi}$ and $\tilde{\phi}$ send $Spec (A)$ to the identity section
 (by
hypothesis on our liftings) so in particular $\eta$ sends $Spec (A[\epsilon ,
\epsilon ']/(\epsilon \epsilon '))$ to the identity pair $(e,e)$ in
$X\times _{\Gg
}X$. On the other hand we can take $Y=Spec (A[\epsilon , \epsilon '])$ and
$Y_1=Spec (A[\epsilon ])$ and $Y_2 = Spec (A[\epsilon '])$ and then apply the
lifting property $Lift _2(Y, Y_i)$ which holds for the morphism $R\rightarrow
X\times _{\Gg}X$ because (from the hypothesis in the property $P4$) this
morphism
is vertical.

Fix a lifting $e_R: S\rightarrow R$ of the identity pair in $X\times
_{\Gg}X$ and fix the values of the morphisms (denoted $\lambda _i$ in the
definition of the lifting property) as being $e_R$ on $Y_1$ and $Y_2$.
These are
indeed liftings of our given morphisms $Y_i \rightarrow X\times _{\Gg}X$
since, as
we have seen above, both $Y_1$ and $Y_2$ map to the identity pair (the subscheme
defined by $(\epsilon \epsilon ')$ is the union of $Y_1$ and $Y_2$).
We obtain by
the lifting property a lifting $Y\rightarrow R$ which agrees with $e_R$ on $Y_1$
and $Y_2$.  If we write (locally) $R=Spec (B)$ then this morphism corresponds to
a morphism $a:B\rightarrow A[\epsilon , \epsilon ']$ such that the projection of
$B$ modulo $\epsilon$ or modulo $\epsilon '$ is a constant morphism
$B\rightarrow
A$. It now follows that $a$ factors through $B\rightarrow A[\delta ]$.  We
obtain a morphism $Spec (A[\delta ])\rightarrow R$ whose projection  into
$X\times X$ is the pair $(\tilde{\varphi}, \tilde{\phi} )$ (that this is the
case is easy to check directly again by supposing that $X$ is affine).  This
implies that $(\tilde{\varphi}, \tilde{\phi} )$ has image in $X\times _{\Gg}X$,
in other words that the morphisms $\tilde{\varphi}$ and $\tilde{\phi}$ from
$Spec
(A[\delta ])$ into $X$ project to the same morphism into $\Gg$. Thus $\varphi =
\phi$, completing the proof of uniqueness.

Now we show existence of the factorization $\gamma = \varphi \circ q$.  The
preceding uniqueness result implies that it is sufficient to construct $\varphi$
after etale localization on $S'$.  Thus we may assume that $\alpha$ and $\beta$
lift to points $\tilde{\alpha}, \tilde{\beta}: Spec (A[\epsilon ])\rightarrow X$
sending $Spec (A)$ to the identity section.  There is a good lifting of the
multiplication in $\Gg$ to a multiplication $X\times X \rightarrow X$ which we
still denote $x\cdot y$, where goodness means  the property $x\cdot e = e
\cdot x =
x$.
We can
now put
$$
\tilde{\gamma }:= \tilde{\alpha }p_1 \cdot \tilde{\beta }p_2 \cdot
(\tilde{\alpha }m p_1)
\cdot (\tilde{\beta }m p_2)
: Spec (A[\epsilon , \epsilon '])\rightarrow X.
$$
We still have the formula that
$$
\alpha \cdot (\alpha m) = e
$$
(this is because the first order term of the composition is just addition of
vectors) and from this formula it follows that $\tilde{\gamma }$ sends the
subschemes $Spec (A[\epsilon ])$ and $Spec (A[\epsilon '])$ to $e$ (through
their
projections to $Spec (A)$). Since now $X$ is a scheme, this implies directly
the
existence of $\tilde{\varphi } : Spec (A[\delta ])\rightarrow X$ such that
$\tilde{\gamma } = \tilde{\varphi } \circ q$.  Projecting from $X$ to $\Gg$ we
get the factorization $\varphi$ desired.

Finally, we set $[\alpha , \beta ] := \varphi$ from the above construction.
It is of course completely clear from the construction that if $S'$ is artinian,
this gives the usual Lie bracket on the algebraic group $\Gg (S' )$.
it remains to be seen that this morphism is bilinear and satisfies the Jacobi
identity (i.e. that a certain deduced trilinear form vanishes).  But these
properties can be checked on values over artinian schemes $S'$, and there since
the bracket we have defined coincides with the usual one, we get bilinearity and
the Jacobi identity.
\eop

{\em Remark:} The subtlety in our whole situation
being essentially that the factorization, while immediate and obviously unique
in the case where the target of the map is a scheme, does not necessarily exist
and may not be unique even if it does exist, when the target of the
map is just a
sheaf.  One can give examples of $P2$ sheaves $\Hh$ on $\Xx /S$ together with
morphisms $Spec (A[\epsilon , \epsilon '])\rightarrow \Hh$ restricting to a
given
section $S\rightarrow \Hh$ over the subscheme defined by
$(\epsilon \epsilon ')$,
and where the morphism either doesn't factor through $Spec (A[\delta ])$ or else
such that the factorization isn't unique.

We indicate here a simpler example which shows the way toward the examples
refered to in the above paragraph.  Let $Y\rightarrow X$ be a degree $2$
morphism
of smooth curves completely ramified above a point $x\in X$.  Let $\Ff$ be the
image of this morphism (considered as a sheaf on $\Xx$).  Let $y$ be the point
lying over $x$ and suppose $f: Spec (k[\epsilon]/\epsilon ^3)\rightarrow
Y$ is a nonzero tangent vector located at $y$.  Then the associated element
of $\Ff(  Spec (k[\tau ]/\tau ^3))$ is constant (equal to the constant
point $y$) on the subscheme
$$
Spec (k[\tau ]/\tau ^2)\subset Spec (k[\tau ]/\tau ^3).
$$
Nevertheless there exists no factorization of the form
$$
Spec (k[\tau ]/\tau ^3)\subset Spec (k[\tau ]/\tau  ^2)\rightarrow
\Ff
$$
(this factorization would have existed had $\Ff$ been a scheme).

We can obtain
an example where a factorization of the type needed in the above theorem doesn't
exist, simply by composing this example with the morphism $Spec (k[\epsilon ,
\epsilon ']) \rightarrow Spec(k[\tau ]/\tau ^3$ sending $\tau$ to $\epsilon +
\epsilon '$.

The sheaf $\Ff$ in this example is not $P4$ with respect to $Spec (k)$.
Of course $\Ff$ is not a group sheaf. As stated elsewhere, I am not
sure about whether a $P2$ group sheaf might not automatically have to be $P4$
for example (or at least satisfy some of the properties we use here).
For example we have seen that an algebraic space (of finite type) which is a
group
is automatically a scheme.

\subnumero{The adjoint representation}
Suppose $\Gg$ is a presentable group sheaf. Then $\Gg$ acts on itself by
conjugation, by the formula
$$
Int (g)(h):= ghg^{-1}.
$$
More precisely this action is a morphism $\Gg \times \Gg \rightarrow \Gg$
and if we put in the identity map on the first projection we obtain a morphism
$\Gg \times\Gg \rightarrow \Gg \times \Gg$ which is a morphism of group objects
(the second variable) over the first variable $\Gg$. From this and from the
invariance of the above definition of the Lie algebra object, this action
induces
an action (the {\em adjoint action})
$$
\Gg \times Lie (\Gg )\rightarrow Lie (\Gg )
$$
which preserves the bracket.

If $(\Ll , [,])$ is a Lie algebra sheaf (that is to say, $\Ll$ is a vector sheaf
with bilinear morphism $[,]: \Ll \times \Ll \rightarrow \Ll$ satisfying the
Jacobi identity) then we obtain a group sheaf $Aut (\Ll , [,])$.

\begin{lemma}
\mylabel{AutLie}
If $(\Ll , [,])$ is a Lie algebra sheaf then $Aut (\Ll , [,])$
is a presentable group sheaf.
\end{lemma}
{\em Proof:}
The group sheaf $Aut (\Ll )$ of automorphisms of the vector sheaf $\Ll$ is
presentable and in particular $P4$ by Theorem \ref{I.1.g}.   The Lie bracket can
be considered as a morphism
$$
\Ll \otimes _{\Oo} \Ll \rightarrow \Ll .
$$
The subgroup $Aut (\Ll , [,])\subset Aut (\Ll )$ may thus be represented as the
equalizer of two morphisms
$$
Aut (\Ll ) \rightarrow Hom (\Ll \otimes _{\Oo} \Ll , \Ll ).
$$
Note that $Hom (\Ll \otimes _{\Oo}\Ll ,  \Ll )$ is a vector sheaf by
Lemma \ref{I.s} and the definition of tensor product following that lemma; and
presentable by Theorem \ref{I.1.g}. In particular $Aut (\Ll , [,])$ is $P4$ by
Lemma \ref{I.1.a} and presentable by Corollary \ref{I.z}. \eop

The adjoint action may be interpreted as a morphism of presentable group sheaves
$$
Ad : \Gg \rightarrow Aut (Lie (\Gg ), [,] ).
$$
We can of course forget about the bracket and compose this with the morphism
into
$Aut(Lie (\Gg ))$ which is just the automorphism sheaf of a vector sheaf.

\begin{proposition}
\mylabel{Adjoint}
Suppose $\Gg$ is a connected presentable group sheaf.  Then the kernel of the
morphism $Ad$ is the center $Z(\Gg )$ (that is to say the sheaf whose values
are the centers of the values of $\Gg$).
\end{proposition}
{\em Proof:}
The statement amounts to saying that a section $g$ of $\Gg$ acts trivially on
$\Gg$ if and only if it acts trivially on $Lie (\Gg )$.  This statement is true,
in fact, of any automorphism (defined over any base scheme $S'\rightarrow S$).
It
suffices to prove this last statement for the values over artinian base schemes
(if an automorphism agrees with the identity on the values over all artinian
base schemes then it must be equal to the identity).  In the case of values over
artinian base schemes it is just the statement that an automorphism  which acts
trivially on the Lie algebra of a connected algebraic group must act trivially
 on
the whole group. \eop

\begin{corollary}
\mylabel{centerPres}
If $\Gg$ is a connected presentable group sheaf, then the center
$Z(\Gg )$  is again presentable.
\end{corollary}
{\em Proof:} By Proposition \ref{Adjoint} the center is the kernel of a morphism
of presentable group sheaves.  By Theorem \ref{I.1.e}, this kernel is
presentable.
\eop

{\em Question}
Suppose $\Gg$ is a presentable group sheaf, not necessarily connected. Is the
center $Z(\Gg )$ presentable?

This is related to the following question.

{\em Question}
Suppose $H$ is a finite presentable group sheaf. Is $Aut (H)$ presentable?

A positive response here would allow us to prove that the center $Z(\Gg )$ is
connected, because it is the kernel of the action of
$Z(\Gg ^o)$ on the group of connected components $H=\Gg /\Gg ^o$.

\subnumero{Determination of presentable group sheaves by their Lie algebras}

The object of this section is to prove the following theorem, which is a
generalization of the well known principle that a Lie group is determined by its
Lie algebra, up to finite coverings, if the center is unipotent.

\begin{lemma}
\mylabel{123}
Suppose $F, G \subset H$ are two presentable group subsheaves of a presentable
group sheaf $H$, and suppose $F$ and $G$ are connected.  If $Lie (F)=Lie (G)$ as
subsheaves of $Lie (H)$ then $F=G$.
\end{lemma}
{\em Proof:}
By the properties B1 and B2 and artin approximation, it suffices to show that
for
any artinian $S'$ we have $F(S')=G(S')$.  But these two are connected Lie
subgroups of $H(S')$ which by hypothesis have the same Lie algebras; thus they
are equal.
\eop

\begin{corollary}
\mylabel{abc}
Suppose $F$ and $G$ are connected presentable group sheaves on $\Xx$.
Suppose $Lie (F)\rightarrow Lie(G)$ is an isomorphism of Lie algebras.  Then
this isomorphism lifts to a unique isomorphism
$F/Z(F)\cong G/Z(G)$ where $Z()$ denotes the center.
\end{corollary}
{\em Proof:}
Note that the center of a
connected presentable group sheaf is presentable by \ref{centerPres}, so
$F/Z(F)$ and $G/Z(G)$ are presentable.

Let $L=Lie (F)=Lie (G)$ and let $A=Aut (L)$ (automorphisms of the vector
sheaf or of the Lie algebra sheaf, we don't care).  We get maps $F\rightarrow A$
and $G\rightarrow A$.  Let $F_1$ and $G_1$ denote the images.  We have
$$
Lie (F_1) = im (L\rightarrow Lie (A)) = Lie (G_1)
$$
as subsheaves of $Lie (A)$, so by Lemma \ref{123} we have $F_1=G_1$. On the
other
hand, note that $Z(F)$ is the kernel of the map $F\rightarrow A$ because if an
element of $F$ acts trivially on $Lie (F)$ then by exponentiation and the fact
that $F$ is connected, it acts trivially on all $F(S')$ for $S'$ artinian
hence in
fact it acts trivially on $F$. Thus $F_1 = F/Z(F)$ and similarly $G_1=G/Z(G)$.
\eop

We have now finished verifying that the class of presentable group sheaves
satisfies the properties set out in the introduction. In effect:
\newline
Property 1 is Corollary \ref{uvw};
\newline
Property 2 is Theorem \ref{I.1.e};
\newline
Property 3 is Lemma \ref{I.1.h};
\newline
Property 4 is Lemma \ref{I.1.i};
\newline
Property 5 is Theorem \ref{I.1.m};
\newline
Property 6 is Theorem \ref{I.1.o} and Corollary \ref{connex};
\newline
Property 7 is Theorem \ref{lmn}; and
\newline
Property 8 is Theorem \ref{abc}.

\subnumero{Questions}

We present in further detail some other questions analogous to well known
properties of algebraic Lie groups, which seem to be more difficult here.

{\bf 1.} \,
 (Existence) {\em If $(\Ll , [,])$ is a
Lie algebra sheaf (i.e. a vector sheaf with bilinear operation satisfying the
Jacobi identity) then does there exist a presentable group sheaf $\Gg$ with $Lie
(\Gg )= (\Ll , [,])$?}  One has the following idea for a proof of existence in
a formal sense.  Take a resolution of $\Ll$ by vector schemes, and lift the
bracket to a bracket (not necessarily satisfying the Jacobi identity) on the
vector scheme $X$ surjecting to $\Ll$.  Then use an explicit version  of
Baker-Campbell-Hausdorff to define a composition law on the formal completion of
$X$ along the zero section.  This composition law will not be associative, but
one should be able to use the second part of the resolution of $\Ll$ to define a
relation scheme $R$ (formally), such that when we set $\Gg$ to be the quotient
of $X$ by $R$ we get a group sheaf.  One would have to check that the maps are
vertical.   Of course this idea for a proof skirts the main question of how to
integrate the formal structure out into an actual presentable group sheaf.

{\bf 2.} \, {\em Does every (connected, say) presentable group sheaf have a
faithful representation on a vector sheaf?}  I guess that the answer is probably
no, but I don't have a specific example in mind.

{\bf 3.} \, {\em Suppose $Lie (\Ff )\rightarrow Lie (\Gg )$ is a morphism of
vector Lie algebras.  Under what conditions does this lift to a morphism $\Ff
'\rightarrow \Gg$ where $\Ff' \rightarrow \Ff$ is a finite covering?}

{\bf 4.} \, {\em What happens in Theorem \ref{abc} if we don't divide out by the
centers?}

{\bf 5.} \, {\em Suppose $G\subset Aut (V)$ is a presentable subgroup of the
automorphisms of a vector sheaf.  Is there a vector subsheaf
$U\subset T^{a,b}(V)$ of a tensor power of $V$ (or possibly a cotensor power
or a
mixture\ldots ) such that $U$ is preserved by the action of $G$ and such that
$G$
is characterized as the subgroup of $Aut (V)$ preserving $U$?}  One of the main
problems in trying to
prove such a statement is that the vector sheaves (and similarly $P4$ or $P5$
sheaves) don't satisfy any nice chain condition.

Note that in the situation of question 4, for any sub-vector sheaf $U$ of a
tensor and cotensor combination of $V$, the subgroup of $Aut (V)$ of elements
preserving $U$ is a presentable subgroup, so at least we obtain a way of
constructing examples, even if we don't know whether we get everything this way.

\numero{Presentable $n$-stacks}

Recall that an $n$-groupoid in the sense of \cite{Tamsamani} is essentially the
same thing as an $n$-truncated homotopy type \cite{Tamsamani2}.  In
view of this,
we can approach the theory of $n$-stacks (we assume from here on that this
means $n$-stack of
$n$-groupoids and drop the word ``groupoid'' from the notation) via the
theory of
presheaves of topological spaces or equivalently simplicial
presheaves\cite{Jardine1}.  We adopt a working convention that by {\em
$n$-stack}
we mean the presheaf of $n$-groupoids associated to a fibrant presheaf of spaces
\cite{Jardine1} \cite{kobe} or, a bit more generally, any presheaf of
$n$-groupoids
such that the associated simplicial presheaf (taking the diagonal of
the nerve) is
fibrant in the sense of \cite{kobe} which means that it satisfies the
global part
of the fibrant condition of Jardine \cite{Jardine1}.

Some special cases are worth mentioning.  A $0$-stack is simply a sheaf of
sets.
A $1$-stack is what is usually called a stack---it is a sort of sheaf of
groupoids.  The notions of $2$-stack and $3$-stack were explored heuristically
from the category-theoretic point of view in \cite{Breen23}.

We suppose given an adequate theory of morphism $n$-stacks $Hom (R,T)$; and of
homotopy fiber products $T\times _RT'$ for $n$-stacks. These can be had, for
example, within the realm of presheaves of spaces \cite{Jardine1} \cite{kobe}
\cite{flexible}.

The path-stack
$P^{t_1,t_2}T$ on $\Xx /S$ between two basepoints (i.e. objects)
$t_1,t_2\in T(S)$
is then well defined.  We denote by $\pi _0(T)$ the truncation
down to a sheaf of
sets, and from this and the path space construction we obtain the homotopy group
sheaves $\pi _i(T,t)$ over $\Xx /S$ for an $n$-stack $T$ and object $t\in T(S)$.

In terms of the easier-to-understand version version of the theory involving
presheaves of spaces, the homotopy group sheaves are defined as follows.  If
$t\in T(S)$ then for any $Y\rightarrow S$ we get a basepoint $t|_Y\in T(Y)$.
The
functor $$
Y/S\mapsto \pi _i (T(Y), t|_Y)
$$
is a presheaf on $\Xx /S$ which we denote by $\pi _i^{\rm pre}(T,t)$. Then $\pi
_i(T,t)$ is sheaf associated to this presheaf.

There is probably a good extension of the theory to $\infty$-stacks which would
correspond to presheaves of spaces which are not necessarily truncated (and I
suppose that it again becomes equivalent to Jardine's theory but there may be a
few subtleties hidden here).  Generally below when we speak of $n$-stacks,
$n$ will
be indeterminate.  There is probably not too much difference between the
theory of
$\infty$-stacks and the projective limit of the theories of $n$-stacks, so we
will stick to the notation $n$-stack.

For $t_1, t_2\in T(S)$ use the notation $\varpi _1(T,t_1,t_2)$ for the sheaf on
$\Xx /S$ of paths in $T|_{\Xx /S}$ from $t_1$ to $t_2$ up to homotopy.  Thus
$$
\varpi _1(T,t_1,t_2)= \pi _0 (P^{t_1,t_2}T).
$$

We make the
following definition.
\newline
---We say that an $n$-stack $T$ on $\Xx$ is {\em presentable} if it satisfies
the
following conditions:

\begin{enumerate}

\item The sheaf $\pi _0(T)$ is P1 over $k$.

\item For any finite type morphism of schemes $Z\rightarrow Y$ and any
two sections
$\eta  : Y\rightarrow T$ and $\eta ': Z\rightarrow T$ the sheaf
$\varpi _1(T|_{\Zz
/Z}, \eta |_Z, \eta ' )$, when restricted down from $Z$ to $Y$, is $P4$
over $Y$.

\item For any scheme $Y$ and section $\eta : Y\rightarrow T$, the higher
homotopy  group sheaves $\pi _i( (T|_{\Zz /Y}), \eta )$, for $i\geq 1$, are
presentable group sheaves ($P5$) over $Y$.

\end{enumerate}

(Recall that if $\Hh$ is a sheaf on $\Xx / Z$ then it can also be considered
as a
sheaf on $\Xx$ with a map to $Z$; the restriction down to $Y$ is the same sheaf
taken
with the composed map to $Y$, then considered as a sheaf on $\Xx /Y$. This
shouldn't be confused with the direct image from $Z$ to $Y$. In
heuristic topological terms the fiber over $y\in Y$ of the restriction is
obtained by taking the direct union of the fibers of $\Hh$ over the points $z$
lying over $y$, whereas  the fiber of the direct image is obtained by taking the
direct product of the fibers of $\Hh$ over points $z$ lying over $y$.)

{\bf Caution:}  This definition of presentability is very slightly different
from
the definition given in \cite{kobe}.  The older version of presentability for
$T$ as defined in \cite{kobe}
corresponds to the property $P3$ for $\pi _0$ (see Theorem
\ref{I.1.q.1kobe} below); whereas the
present definition corresponds to the property $P3\frac{1}{2}$ (see
Theorem
\ref{I.1.q.1} below). I
hope that the present version corresponding to $P3\frac{1}{2}$ will be the most
useful. The reason for changing the definition was to be able to state
Theorem \ref{stability} in a nice way, i.e. to have a reasonable definition of
{\em presentable morphism} of $n$-stacks.

{\em Caution:}  If $T$ is $0$-truncated, that is a sheaf of sets, and happens
to have a group structure, then this notion is not the same as the  notion
that $T$ be a presentable group sheaf.  The presentability in $T$ as
defined here refers to the higher homotopy groups. In fact, presentability in
this case corresponds to the property $P3\frac{1}{2}$ rather than $P4$
(see below).

We can also reasonably use the notations {\em presentable homotopy sheaf};  {\em
presentable space over $\Xx$} or just {\em presentable space}; or {\em
presentable fibrant presheaf of spaces}, for the notion of presentable
$n$-stack.

Property $1$ implies the seemingly stronger statement that
there is a section $f: Z\rightarrow T$ over a
scheme $Z$ of finite type over $k$, such that the associated morphism
$Z\rightarrow \pi _0(T)$ is surjective.

The second condition reduces, in the case $\eta = \eta '$, to the statement
that for any scheme $Y$ and section $\eta : Y\rightarrow T$, the fundamental
group sheaf $\pi _1( (T|_{\Zz /Y}), \eta )$ is a presentable group sheaf over
$Y$.

We can give an alternative characterization, from which it follows
that any truncation $\tau _{\leq n}T$ of a presentable space is again
presentable.
Recall that we have defined a
condition $P3\frac{1}{2}$ which is intermediate between $P2$ and $P4$.

\begin{theorem}
\mylabel{I.1.q.1}
Suppose $T$ is an $n$-stack over $X$.  Then $T$ is presentable
if and only if the sheaf $\pi _0$ is $P3\frac{1}{2}$, and for any $Y\in \Xx$ and
$t\in T(Y)$, the sheaves $\pi _i (T|_{\Xx /Y}, t)$ are presentable group sheaves
($P5$) over $Y$.
\end{theorem}
{\em Proof:}
Suppose $T$ is presentable.  Then we just have to show that $\pi _0$ is
$P3\frac{1}{2}$.  We know that it is P1, so there is a surjection
$Y\rightarrow \pi _0$.  By replacing $Y$ by an etale cover, we may assume that
this comes from a point $t\in T(Y)$.  The path space $P^{p_1^{\ast}t,
p_2^{\ast}t}T$ maps to $Y\times Y$, and
$$
\varpi _1(T, p_1^{\ast}t,
p_2^{\ast}t )=\pi _0(P^{p_1^{\ast}t,
p_2^{\ast}t}T)\rightarrow Y\times _{\pi _0} Y
$$
is surjective. Let $G\rightarrow Y$ be the sheaf of groups $\pi _1(T|_Y,t)$.
It is presentable by hypothesis, and $G$ acts freely on
(the restriction from $Y\times Y$ down to $Y$ of) $\varpi _1(T, p_1^{\ast}t,
p_2^{\ast}t )$ with quotient  $Y\times _{\pi _0}Y$. Finally, we know
that (the restriction from $Y\times Y$ down to $Y$ of)
$\varpi _1(T, p_1^{\ast}t,
p_2^{\ast}t )$ is $P4$ over $Y$; thus the quotient
$Y\times _{\pi _0}Y$ is $P4$ over $Y$ by Theorem
\ref{I.1.d}. Now by definition there exists a surjective morphism
$Q\rightarrow Y\times _{\pi _0}Y$ which is $Y$-vertical. This is what
is required to show that $\pi _0$ is $P3\frac{1}{2}$.

Now suppose that $\pi _0$ is $P3\frac{1}{2}$ and that the other homotopy
group sheaves are presentable. We obtain immediately that $\pi _0$ is P1.   Let
$X\rightarrow \pi _0$ be the surjection given by the property $P3\frac{1}{2}$.
Then we have
an $X$-vertical surjection $Q\rightarrow X\times _{\pi _0}X$ (where the morphism
to $X$ is the first projection).  Suppose $X'\rightarrow X$ is an etale
surjection
chosen so that the map $X\rightarrow \pi _0$ lifts to $t\in T(X')$.  Let $Q'$ be
the pullback of $X' \times X'$ to $Q$.  Then $Q'= (X'\times _{\pi _0}X')\times
_{X\times _{\pi _0}X}Q$ so $Q'\rightarrow X'\times _{\pi _0}X'$ is
$X$-vertical, and hence $X'$-vertical. This implies that $X'\times _{\pi
_0}X'$ is $P4$ over $X'$, because  we can take as the relation scheme  $$
Q'\times
_{X'\times _{\pi _0}X'}Q'= Q'\times _{X'\times X'}Q'
$$
which is already a
scheme of finite type (and the identity is vertical). Now we have a sheaf of
groups $G= \pi _1(T|_{X'}, t)$ over $X'$ which is by hypothesis presentable, and
$G$ acts freely on $\varpi _1(T, p_1^{\ast}t,
p_2^{\ast}t )$ with
quotient $X'\times _{\pi _0}X'$.  By Theorem \ref{I.1.d},
$\varpi _1(T, p_1^{\ast}t,
p_2^{\ast}t )$ is $P4$ over $X'$.

Now suppose that we have a finite type morphism $q:Z\rightarrow Y$ and  two
points $\eta _1 \in T(Y)$ and $\eta _2 \in T(Z)$, and we show that the
restriction
from $Z$ to $Y$ of the path space $\varpi _1(T, \eta _1 |_Z,\eta _2)$
is $P4$ over
$Y$. There are etale surjections $ Y'\rightarrow Y$ and $Z'\rightarrow Z$ (of
finite type) with $Z'\rightarrow Y'$  and there are morphisms
$f_1:Y'\rightarrow X'$ and $f_2: Z'\rightarrow X'$ such that $f_1^{\ast} (t)$ is
homotopic to  $\eta _1|_{Y'}$ and
$f_2^{\ast} (t)$ is
homotopic to  $\eta _2|_{Z'}$. Let $(f_1|_{Z'},f_2): Z'\rightarrow X'\times
X'$ denote the resulting morphism (the first projection of which factors through
$Y'$).
Then
$$
\varpi _1(T,\eta _1|_{Z}, \eta _2)|_{Z'}=
\varpi _1(T,\eta _1|_{Z'}, \eta _2|_{Z'})= (f_1|_{Z'},f_2)^{\ast}
\varpi _1(T,p_1^{\ast}t, p_2^{\ast}t)
$$
$$
= (q, f_2)^{\ast}[\varpi _1(T,p_1^{\ast}t, p_2^{\ast}t)|_{Y'\times X'}].
$$
Note that $\varpi _1(T,p_1^{\ast}t, p_2^{\ast}t)|_{Y'\times X'}$ is
$P4$ with respect to $Y'$, so by the appendix to the proof below,
one gets that the restriction down to $Y'$ of $\varpi _1(T,\eta _1|_{Z}, \eta
_2)|_{Z'}$ is $P4$ with respect to $Y'$.    By Corollary \ref{I.1.j.1}, the
restriction down to $Y$ of $\varpi _1(T,\eta _1|_{Z}, \eta _2)$ is $P4$
over $Y$.
\eop

{\em Appendix to the proof:}
Suppose $Z\rightarrow Y$ is a finite type morphism, and suppose $\Ff$ is a
sheaf on
$Y$.  Then the restriction from $Z$ down to $Y$ of the pullback $\Ff |_Z$ is
equal to the fiber product $Z\times _Y\Ff$.  Note also that $Z$ is $P4$ over
$Y$.  Thus if $\Ff$ is $P4$ over $Y$ then the restriction of the pullback is
again $P4$.

\begin{corollary}
\mylabel{truncation}
If $T$ is a presentable $n$-stack and if $m<n$ then $\tau _{\leq m}T$ is a
presentable $m$-stack.
\end{corollary}
{\em Proof:}
Indeed the truncation operation preserves the homotopy group sheaves (and the
homotopy sheaf $\pi _0$). By the theorem, presentability is expressed solely in
terms of these sheaves so it is preserved by truncation.
\eop

We have a similar theorem for the old version of presentability of $T$ \cite
{kobe}.

\begin{theorem}
\mylabel{I.1.q.1kobe}
Suppose $T$ is an $n$-stack over $X$.  Then $T$ is presentable in the sense of
\cite{kobe} if and only if the sheaf $\pi _0$ is $P3$, and for any
$Y\in
\Xx$ and $t\in T(Y)$, the sheaves $\pi _i (T|_{\Xx /Y}, t)$ are presentable
group
sheaves ($P5$) over $Y$.
\end{theorem}
{\em Proof:} The proof is the same as above only very slightly easier.
The details
are left to the reader.
\eop

\subnumero{Very presentable $n$-stacks}

We make the following more restrictive definition.  Say that a presentable
group sheaf $G$ on $\Xx /S$ is {\em affine} if, for any artinian $S$-scheme
$S'$, the group scheme $G(S')$ over $Spec (k)$ is affine. A truncated homotopy
sheaf $T$ is {\em very presentable} if $T$ is presentable and if for any $\eta
\in T_Y$ we have that $\pi _1(T/Y,\eta )$ is affine, and $\pi _i(T/Y, \eta )$
are vector sheaves for $i\geq 2$.

The idea behind the definition of ``very presentable'' is that we want to
require
the higher homotopy groups to be unipotent.
Note that
if we don't require $\pi _1$ to be affine, or $\pi _i$ to be unipotent $(i\geq
2$), then the comparison between algebraic and analytic de Rham cohomology
(announced in \cite{kobe}) is no longer true, even over the base $S=Spec (k)$
when all of the groups are representable.   This is the reason for making the
definition of ``very presentable''.

I make the following conjecture:

\begin{conjecture}
\mylabel{I.1.r}
If $G$ is an abelian affine presentable group sheaf on $\Xx /S$ such that
for any
artinian $S'\rightarrow S$ the group scheme $G(S')$ over $k$ is a direct sum
of additive groups, then $G$ is a vector sheaf.
\end{conjecture}

If we knew this conjecture, we could replace the condition of being a vector
sheaf by the condition that the $G(S')$ are unipotent (hence additive) for
$G=\pi _i$, $i\geq 2$; this would then be along the same lines as the affineness
condition for $\pi _1$. As it is, we need to require the condition of $\pi
_i$ being vector schemes ($i\geq 2$) for many of the arguments concerning
de Rham
cohomology sketched in \cite{kobe} to work.

{\em Remark:} The categories of presentable and very presentable $n$-stacks are
closed under weak equivalences and fiber products but not under cofiber products
(push-outs); thus they are not closed model categories.

{\em Remark:} We have the same statement as Corollary \ref{truncation} for very
presentable stacks (if $T$ is very presentable then $\tau _{\leq m} T$ is very
presentable).

\subnumero{Other presentability conditions}

Recall from \cite{kobe} that we used the notation $P6$ for affine presentable
group sheaves and $P7$ for vector sheaves.  An $n$-stack $T$ on $\Xx$ is
{\em $(a_0,\ldots , a_n)$-presentable} (with $a_i \in \{ 0,1, 2 ,3,
3\frac{1}{2},
4, 5,6,  7\}$) if  $\pi _0(T)$ is $Pa_0$ and if for any scheme $Y$ and
$t\in T(Y)$,
$\pi _i (T, t)$ is $Pa_i$ over $Y$.  Here by convention $P0$ means no
condition at
all. Thus a presentable $n$-stack in our previous notation becomes a
$(3\frac{1}{2},5,5, \ldots )$-presentable $n$-stack in this notation.  A very
presentable $n$-stack is a $( 3\frac{1}{2},
6, 7, 7, \ldots )$-presentable $n$-stack.  The old notions of presentability
and very presentability as defined in \cite{kobe} are respectively
$(3,5,5,\ldots
)$-presentability and $(3,6,7,7, \ldots )$ presentability.  There may be some
interest in considering, for example, the $(2,2,2,\ldots )$-presentable
$n$-stacks, or the $(0,0, 7,7,7,\ldots )$-presentable $n$-stacks.

Some other useful versions might be $(4, 5, 5, \ldots
)$-presentable $n$-stacks, or $(4, 6, 7, 7, \ldots )$-presentable $n$-stacks
for example. Here the condition $P4$ on $\pi _0$ would be with
respect to $S=Spec
(k)$. For example an algebraic stack with smooth  morphisms from the morphism
scheme to the object scheme (or even more strongly a Deligne-Mumford
stack where
these morphisms are etale) would be a $(4, 5)$-presentable stack.  The converse
is not true since in the condition of $(4,5)$-presentability, the morphism
sheaves are not necessarily representable. In fact we will never see the
condition of representability of the morphism sheaves in our context, since this
is unnatural from the point of view of higher-order stacks (and even in the
context of algebraic stacks, one may wonder why the morphism object itself was
never allowed
to be an algebraic space?).

{\em Remark:}
Again we have the statement of Corollary \ref{truncation}: if $T$ is an
$(a_0,\ldots , a_n)$-presentable $n$-stack then $\tau _{\leq m}T$ is an
$(a_0,\ldots , a_m)$-presentable $m$-stack.

{\em Remark:}  A good convention for using all of these different notions would
be to chose some variables $A$, $B$, etc. and set them to be specific $(a_0,
a_1, \ldots )$ at the start of a discussion, then to use the notation
``$A$-presentable'' or ``$B$-presentable'' throughout the discussion.

\subnumero{A relative version of presentability}

We can make a relative definition.  In general, say that a morphism
$T\rightarrow
R$ of $n$-stacks is {\em $(a_0,\ldots , a_n)$-presentable} if for any scheme
$Y\in \Xx$ and any morphism $Y\rightarrow R$, the fiber $T\times _RY$ is
$(a_0,\ldots , a_n)$-presentable.  In particular we obtain the notions of
presentable and very presentable morphisms by taking
$(3\frac{1}{2},5,5, \ldots )$ and $(3\frac{1}{2},6,7,7, \ldots )$ respectively.

It is clear that if $T\rightarrow R$ is an $(a_0,\ldots , a_n)$-presentable
morphism and if $R'\rightarrow R$ is any morphism of $n$-stacks then the
morphism
$T':= T\times _RR'\rightarrow R'$ is $(a_0,\ldots , a_n)$-presentable.

\begin{lemma}
\mylabel{structural}
Suppose that $a_0 \leq 5$.
An $n$-stack $T$ on $\Xx$ is
$(a_0,\ldots , a_n)$-presentable if and only if the structural morphism
$T\rightarrow \ast$ is $(a_0,\ldots , a_n)$-presentable.
\end{lemma}
{\em Proof:}
Since $\ast$ is itself a scheme of finite type (it is $Spec (k)$) the structural
morphism being $(a_0,\ldots , a_n)$-presentable implies that $T$ is
$(a_0,\ldots , a_n)$-presentable.

For the other implication, suppose $T$ is
$(a_0,\ldots , a_n)$-presentable, then for any scheme of finite type $Y$ we have
that $T\times Y = T\times _{\ast}Y$ is $(a_0,\ldots , a_n)$-presentable (since
a scheme $Y$ is $a_0$-presentable for any $a_0 \leq 5$).
\eop

{\em Remark:} If $\Gg$ is a sheaf of groups on $\Xx /S$ then $\Gg$ is a
presentable group sheaf if and only if $K(\Gg , 1)\rightarrow S$ is a
presentable morphism of $1$-stacks. This is the correct point of view relating
our terminologies ``presentable group sheaf'' and ``presentable morphism'' or
``presentable $n$-stack'', i.e. the answer to the terminological problem posed
by the caution at the start of this section.

\begin{theorem}
\mylabel{stability}
Suppose $R$ is a presentable (resp. very presentable) $n$-stack.  Then
a morphism
$T\rightarrow R$ is presentable (resp. very presentable) if and only if $T$
itself is presentable (resp. very presentable).
\end{theorem}

The proof of this theorem will be given in the next subsection below.
We first state a few corollaries.

\begin{corollary}
\mylabel{fiberprod}
Suppose $T\rightarrow R$ and $S\rightarrow R$ are morphisms between presentable
(resp. very presentable) $n$-stacks.  Then the fiber product $T\times _RS$ is
presentable (resp. very presentable).
\end{corollary}
{\em  Proof:}
From the theorem, the morphism $T \rightarrow R$ is presentable, hence the
morphism $T\times _RS$ is presentable and since $S$ is presentable, again from
the theorem we conclude that $T\times_RS$ is presentable. The same goes for very
presentable. \eop

\begin{lemma}
\mylabel{basechange}
Suppose $R'\rightarrow R$ is a morphism inducing a surjection on $\pi _0$.
Then a morphism $T\rightarrow R$ is presentable (resp. very presentable) if and
only if the morphism $T':= T\times _RR'\rightarrow R'$ is presentable (resp.
very
presentable).  \end{lemma}
{\em  Proof:}
One direction follows directly from the first remark after the definition above.
For the other direction, suppose that $T'\rightarrow R'$ is presentable
(resp. very
presentable).  Then
for any scheme $Y\rightarrow R$ there is an etale covering $Y' \rightarrow Y$
and a lifting $Y'\rightarrow R'$, and we have
$$
(T\times _RY)\times _YY'=T\times _RY'=T' \times _{R'}Y',
$$
which is presentable (resp. very
presentable) by hypothesis.  The conditions on homotopy sheaves for
being presentable (resp. very
presentable) are etale-local, so $T\times _RY$ is presentable (resp. very
presentable).  \eop

\begin{corollary}
\mylabel{composition}
Suppose $R\rightarrow S$ and $S\rightarrow T$ are presentable (resp. very
presentable) morphisms of $n$-stacks. Then the composition $R\rightarrow T$ is a
presentable (resp. very presentable) morphism.
\end{corollary}
{\em Proof:}
Suppose $X$ is a scheme of finite type with a morphism $X\rightarrow T$.
Then
$$
X\times _TR = (X\times _TS) \times _SR.
$$
By hypothesis, $(X\times _TS)$ is presentable (resp. very presentable),
and by the other hypothesis and the base change property given at the start of
the subsection, the morphism $(X\times _TS) \times _SR\rightarrow
(X\times _TS)$ is presentable (resp. very presentable).  Theorem \ref{stability}
now implies that  $X\times _TR$ is presentable (resp. very presentable).  By
definition then, the morphism $R\rightarrow T$ is presentable (resp. very
presentable).
\eop

\begin{corollary}
\mylabel{check}
Suppose $f:T\rightarrow R$ is  a morphism such that $R$ is presentable (resp.
very presentable), and suppose $X\rightarrow
R$ is a morphism from a scheme of finite type $X$ which is surjective on
$\pi _0$.
Then $T$ and the morphism $f$ are  presentable (resp. very presentable)
if and only if $T\times _RX$ is presentable (resp. very presentable).
\end{corollary}
{\em  Proof:}
By Lemma \ref{basechange} the morphism $f$ is presentable if and only if the
morphism $p_2: T\times _RX\rightarrow X$ is presentable.  On the other
hand, $T$ is
presentable if and only if $f$ is presentable, from Theorem \ref{stability}.
Similarly  $T\times _RX$ is presentable if and only if $p_2$ is presentable
again by \ref{stability}. This gives the desired statement (the
same proof holds for
very presentable).
\eop

We now give some results that will be used in the proof of Theorem
\ref{stability}.

\begin{lemma}
\mylabel{vector?}
Suppose $V$ is a vector sheaf and $G$ is a presentable group sheaf on
$\Xx /S$. If
$f: V\rightarrow G$ is a morphism of group sheaves then the kernel of $f$ is a
vector sheaf.
\end{lemma}
{\em Proof:}
There is a natural isomorphism of vector sheaves $\varphi : V \cong Lie (V)$,
such that $\varphi$ reduces to the exponential on the values over artinian $S'$.
To construct $\varphi$ note that a section of $V$ may be interpreted as a map
$\Oo \rightarrow V$.  We have a tautological section of $Lie (\Oo )$ so for
every section of $V$ the image of this tautological section is a section of
$Lie (V)$.    This map is an isomorphism on values over artinian schemes,
so it is an isomorphism.

Let $U \subset Lie (V)$ be the kernel of
$$
Lie (f) : Lie
(V)\rightarrow Lie (G).
$$
Since $Lie (f)$ is a morphism of vector sheaves, its
kernel $U$ is a vector sheaf.  We claim that $\varphi ^{-1}(U)$ is the kernel of
$f$.  In order to prove this claim it suffices to prove it for the values over
artinian $S'$ (since both are presentable and contained in $V$, and using
\ref{Krull}). Here it reduces to the following statement about Lie groups: the
kernel of an algebraic morphism from a vector space to a Lie group is the
exponential of the kernel of the corresponding morphism of Lie algebras.  To
prove this notice first that this exponential is a subvector subspace; we can
take the quotient and then we are reduced to the case where the map is injective
on Lie algebras.  The kernel is thus a finite subgroup, but a vector space
contains no finite subgroups so we are done.
\eop

\begin{proposition}
\mylabel{I.1.s.3}
Suppose $R$, $S$ and $T$ are $n$-stacks over $\Xx$, with  morphisms
$R\rightarrow T$ and $S\rightarrow T$.  Suppose $Z\in \Xx$ and $(r,s)\in
R\times _TS(Z)$. Let $t\in T(Z)$ be the common image of $r$ and $s$. Then we
have the following long exact sequence of homotopy group sheaves on $\Xx /Z$:
$$
\ldots \rightarrow \pi _i (R\times _TS|_{\Xx /Z},(r,s))\rightarrow
\pi _i(R|_{\Xx /Z},r)\times \pi _i (S|_{\Xx /Z})\rightarrow
$$
$$
\pi _i (T|_{\Xx /Z},t)\rightarrow \pi _{i-1}(R\times _TS|_{\Xx /Z},(r,s))
\rightarrow \ldots ,
$$
terminating with the sequence
$$
\pi _2(R|_{\Xx /Z},r)\times \pi _2(S|_{\Xx /Z},s)\rightarrow
\pi _2(T|_{\Xx /Z},t)\rightarrow \pi _1(R\times
_TS|_{\Xx /Z}, (r,s))\rightarrow
$$
$$
\pi _1(R|_{\Xx /Z},r)\times
\pi _1(S|_{\Xx /Z},s) \stackrel{\displaystyle
\rightarrow }{\rightarrow } \pi _1(T|_{\Xx /Z},t)
$$
(the last part meaning that the image is equal to the equalizer of the two
arrows). Furthermore, we have a similar sequence for the path spaces.  Suppose
$(r_1,s_1)$ and $(r_2,s_2)$ are two points, with images $t_1$ and $t_2$.
We have the exact sequence
$$
\pi _2(R|_{\Xx /Z},r_1)\times \pi _2(S|_{\Xx /Z},s_1)\rightarrow
\pi _2(T|_{\Xx /Z},t_1)\stackrel{acts\; on}{\rightarrow} \varpi _1(R\times
_TS|_{\Xx /Z}, (r_1,s_1),(r_2,s_2))
$$
$$
\mbox{with quotient the equalizer of} \varpi
_1(R|_{\Xx /Z},r_1,r_2)\times  \varpi _1(S|_{\Xx /Z},s_1,s_2)
\stackrel{\displaystyle \rightarrow }{\rightarrow } \varpi _1(T|_{\Xx
/Z},t_1,t_2). $$
\end{proposition}
{\em Proof:}
We show that we have similar exact sequences at the homotopy presheaf level;
then the sequences for the homotopy sheaves follow by sheafification. To
define the exact sequences at the presheaf level, we can work object by
object, so it suffices to give functorial exact sequences for fibrations of
topological spaces $R\rightarrow T$ and $S\rightarrow T$ with basepoints
$(r,s)$ mapping to $t$.  The morphisms are defined as follows.  The morphism
from $\pi _i (R\times _TS,(r,s))$ to $\pi _i(R,r)\times \pi _i (S,s)= \pi _i
(R\times S, (r,s))$ comes from the inclusion $R\times _TS\rightarrow R\times
S$. The morphism from the product to $\pi _i (T,t)$ is the difference of the
two projection maps.  The morphism from $\pi _i (T,t)$ to $\pi _{i-1}(R\times
_TS,(r,s))$ is obtained as a composition
$$
\pi _i (T,t)\stackrel{\delta}{\rightarrow} \pi _{i-1} (R_t,r)
\stackrel{(1, 0_s)}{\rightarrow }\pi _{i-1}(R_t\times S_t,(r,s))
\stackrel{i}{\rightarrow}\pi _{i-1}(R\times _TS,(r,s))
$$
where $\delta$ is the connecting homomorphism for the fibration $R\rightarrow
T$, $0_s$ is the constant class at the basepoint $s$, and $i$ is the inclusion
of the fiber $i: R_t\times S_t\rightarrow R\times _TS$.  If we took $(1,1)$
instead of $(1,0_s)$ we would get the connecting morphism for the fibration
$R\times _TS\rightarrow T$, which goes to zero in the homotopy of the total
space $R\times _TS$.  Thus, our map is the same as the map which would be
obtained by putting in $-(0_r,1)$ instead.  From the equality of these two
maps, one obtains that the composition of this map with the difference of
projections, is equal to zero.  That the other compositions are zero is easy to
see.  Exactness follows by making a diagram with this sequence on one
horizontal row, with the sequence
$$
\pi _i(R_t\times S_t, (r,s))= \pi _i(R_t,r)\times \pi _i (S_t,s)\rightarrow
0\rightarrow \ldots
$$
on the row above, and the sequence
$$
\pi _i (T,t)\rightarrow \pi _i (T,t)\times \pi _i (T,t)\rightarrow \pi _i
(T,t) \stackrel{0}{\rightarrow} \pi _{i-1}(T,t)
$$
on the row below.  The vertical rows then have the exact fibration sequences
going downwards.  One obtains the exactness of the sequence of homotopy
groups in question (this works at the end by using the extension of the
homotopy sequence for a fibration, to the action of $\pi _1$ of the base on
$\pi _0$ of the fiber, with the $\pi _1$ of the total space being the
stabilizer of the component of $\pi _0$ of the fiber containing the
basepoint.

Finally, we treat the case of the path spaces. What is written on the left
means, more precisely, that the cokernel of the first map acts freely on the
middle sheaf, with quotient equal to the equalizer.  The action in question is
by the map to $\pi _1(R\times _TS, (r_1,s_1))$ which itself acts on the path
space. Now if $\varpi _1(R\times _TS, (r_1,s_1), (r_2,s_2))$ is empty then
the equalizer in question is also empty (any element of the equalizer can be
realized as a pair of paths mapping to exactly the same path in $T$, giving a
path in the fiber product).  Note that we count an action on the empty set as
free.  So we may assume that $\varpi _1(R\times _TS, (r_1,s_1), (r_2,s_2))$ is
nonempty, and choose an element.  This choice gives compatible choices in all
the other path spaces, so composing with the inverse of this path we reduce to
the exact sequence for fundamental groups. \eop

{\em Remark:}  We can extend this sequence to a statement involving $\pi _0$,
specially in the case of a fibration sequence. This will be done as we need it
below.

\begin{lemma}
\mylabel{kernel}
Suppose $S$ is a base scheme and suppose $\Ff$ is a sheaf on $\Xx /S$ whose
restriction down to $\Xx$ is $P3\frac{1}{2}$.  Suppose that $\Gg$ is a
$P4$  sheaf on $\Xx /S$ with morphism $\Gg \rightarrow \Ff$, and finally suppose
that $\eta : S\rightarrow \Ff$ is a section.  Then the inverse image $\Hh
\subset
\Gg$ of the section $\eta$ is a $P4$ sheaf.
\end{lemma}
{\em Proof:}
Let $X\rightarrow \Gg$ and $W\rightarrow X\times _{\Gg}X$ be the $S$-vertical
surjections for $\Gg$.  Fix a surjection $Z\rightarrow \Ff$ and a surjection
$W\rightarrow Z\times _{\Ff}Z$ which is vertical with respect to the first
projection to $Z$.  Fix a lifting $\eta '$ of the section to $Z$ (note that we
are allowed to etale-localize on the base $S$).  Let $U:= S\times _{Z}W$
where the
morphism in the fiber product is the first projection from $W$ to $Z$ (note that
$U$ is a scheme of finite type over $S$). The surjective morphism
$$
U\rightarrow S\times _{Z}(Z\times _{\Ff}Z) = S\times _{\Ff} Z
$$
is $S$-vertical since the morphism $W\rightarrow Z\times _{\Ff}Z$
was $Z$-vertical.
We can choose a lifting $X\rightarrow Z$ of the morphism $\Gg \rightarrow \Ff$.
Then
$$
S\times _{\Ff} X= (S\times _{\Ff}Z)\times _ZX
$$
so there is an $S$-vertical morphism
$$
U\times _Z X \rightarrow S\times _{\Ff} X.
$$
On the other hand the $S$-vertical morphism $X\rightarrow \Gg$ gives
an $S$-vertical morphism
$$
S\times _{\Ff} X\rightarrow S\times _{\Ff} \Gg = \Hh .
$$
Note that $Y:= U\times _Z X$ is a scheme of finite type with a surjective
vertical
morphism to $\Hh$.  Since $\Gg$ is $P4$ there exists a scheme of finite type $V$
and an $S$-vertical morphism
$$
V\rightarrow Y\times _{\Gg}Y = Y\times _{\Hh} Y.
$$
This gives the condition $P4$ for $\Hh$.
\eop

The following lemma is a consequence of Corollary
\ref{fiberprod}, but we need
it in the proof of Theorem \ref{stability}.

\begin{lemma}
\mylabel{I.1.s.?}
If $R$ and $S$ are presentable (resp. very presentable) $n$-stacks over
$\Xx$ and $X$ a scheme of finite type, with  morphisms $R\rightarrow S$ and
$X\rightarrow S$, then the homotopy fiber product $X\times _SR$ is presentable
(resp. very presentable).
\end{lemma}
{\em Proof:}
Suppose $f:Y\rightarrow X\times _SR$ is a morphism. Let $r: Y\rightarrow R$ and
$s: Y\rightarrow S$ be the composed morphisms. Then (since $X$ is
zero-truncated)
for $i\geq 1$ we have
$$
\pi _i(X\times _SR |_{\Xx /Y}, f)= \pi _i (Y\times _SR/Y, r).
$$
The latter fits into a homotopy exact sequence, which we can therefore write
$$
\ldots \pi _{i+1}(S|_{\Xx /Y}, s)\rightarrow
\pi _i(X\times _SR |_{\Xx /Y}, f)\rightarrow \pi _i(R|_{\Xx /Y}, r)\rightarrow
\ldots .
$$
In the presentable case we obtain immediately from Theorem \ref{I.1.e} that
$\pi _i(X\times _SR |_{\Xx /Y}, f)$ is a presentable group sheaf over $Y$.
In the very presentable case,
for $i\geq 3$ we obtain immediately (from Corollary \ref{I.j} and Theorem
\ref{I.k}) that $\pi _i(X\times _SR |_{\Xx /Y}, f)$ is a vector sheaf.  For
$i=2$ we obtain the same conclusion but must also use Lemma \ref{vector?}.
For $i=1$ we obtain that $\pi _1(X\times _SR |_{\Xx /Y}, f)$ is $P5$. In fact it
is an extension of the kernel of a morphism of $P6$ group sheaves, by a vector
sheaf.  Therefore it is also affine (since kernels and extensions by vector
sheaves at least, preserve the affineness property).  Thus it is $P6$.

We just have to prove (in both the presentable and very  presentable case)
that $\pi _0(X\times _SR)$ is $P3\frac{1}{2}$.  Let $a: X\rightarrow S$ denote
the given morphism. Recall that $\pi _0(X\times _SR)/X$ denotes this sheaf
considered as a sheaf on $\Xx /X$.
We  have an action of $\pi _0(S|_{\Xx /X}, a)$ (which is a $P5$ group sheaf over
$X$) on $\pi _0(X\times _SR)/X$, and the quotient is the fiber product $X\times
_{\pi _0(S)}\pi _0(R)/X$ (i.e. again considered as a sheaf over $\Xx /X$).  This
is the same thing as the inverse image of the given section $a$ via the map $\pi
_0(R|_{\Xx /X})\rightarrow \pi _0(S|_{\Xx /X})$.  By Corollary \ref{P3c} or
\ref{P3d} the quotient by the action is $P3\frac{1}{2}$.  Finally by Proposition
\ref{P3e}, the sheaf $\pi _0(X\times _SR)$ is $P3\frac{1}{2}$.
\eop

{\em Remark:}
A similar technique allows one to directly prove Corollary \ref{fiberprod}, that
if $R$, $S$ and $T$ are presentable (resp. very presentable) $n$-stacks with
morphisms $R\rightarrow S$ and $T\rightarrow S$ then the fiber product $R\times
_ST$ is presentable (resp. very presentable). This is left to the reader.  Our
technique is to use only the above special case to get Theorem \ref{stability},
and then to deduce Corollary \ref{fiberprod} as a consequence.

\subnumero{The proof of Theorem \ref{stability}}

Lemma \ref{I.1.s.?} immediately implies one direction in Theorem
\ref{stability},
namely that if $R$ and $S$ are presentable then the morphism $f$ is presentable.

We have to show the other direction: suppose $S$ is a presentable
$n$-stack, $R$ is an $n$-stack, and $f:R\rightarrow S$ is a presentable
morphism.  Choose a scheme of finite type $X$ with a morphism $X\rightarrow S$
inducing a surjection on $\pi _0$.  We will show that if $X\times _SR$ is
presentable then $R$ is presentable.

First of all the morphism $\pi _0(X\times _SR)\rightarrow \pi _0(R)$ is
surjective so if $\pi _0(X\times _SR)$ is $P1$ then so is $\pi _0(R)$.
For the higher homotopy groups, suppose that $s:Z\rightarrow R$ is a morphism.
Lift the projection into $S$ (denoted by $s$) to a morphism $Z\rightarrow X$.
This gives a point $f: Z\rightarrow X\times _SR$ and by composition
$f_Z: Z\rightarrow Z\times _SR= Z\times _X(X\times _SR)$.
Then we have the exact sequence
$$
\ldots \rightarrow \pi _i(Z\times _SR |_{\Xx /Z}, f_Z)
\rightarrow \pi _i(R|_{\Xx /Z}, r)\rightarrow \pi _i(S|_{\Xx /Z}, s)\rightarrow
\ldots .
$$
But since $Z$ and $X|_{\Xx /Z}$ are zero-truncated, and we have that
$Z\times _SR = Z\times _X(X\times _SR)$, the higher
homotopy groups $\pi _i(Z\times _SR |_{\Xx /Z}, f_Z)$ are the same as the
$\pi _i(X\times _SR |_{\Xx /Z}, f)$.
Thus we can write the exact sequence as
$$
\ldots \rightarrow \pi _i(X\times _SR |_{\Xx /Z}, f)
\rightarrow \pi _i(R|_{\Xx /Z}, r)\rightarrow \pi _i(S|_{\Xx /Z}, s)\rightarrow
\ldots .
$$
Note that (in the very
presentable case) the kernel of the morphism
$$
\pi _2(S|_{\Xx /Z}, s)\rightarrow \pi _1(X\times _SR |_{\Xx /Z}, f)
$$
is a vector sheaf by Lemma \ref{vector?}.  In the other cases the kernel (and
the
cokernel on the other end) are automatically vector sheaves by Corollary
\ref{I.j}.  Since the property of being a vector sheaf is preserved under
extension we get the condition that the $\pi _i(R|_{\Xx /Z}, r)$ are vector
sheaves ($i\geq 2$).  In the presentable case the exact sequence immediately
gives the property $P5$ for  the $\pi _i(R|_{\Xx /Z}, r)$ for  ($i\geq 2$).

We have to treat the case of $\varpi _1$.  Suppose $Z\rightarrow Y$ is a
morphism
of finite type and suppose $r, r': Z\rightarrow R$ are points such that $r$
factors through $Y$.  Let $s,s'$ denote the images in $S$ and assume that they
lift to points $f, f'$ and $f_Z, f'_Z$ as above (with $f$ or $f_Z$ factoring
through $Y$).

We first study everything on the level of sheaves on $\Xx /Z$.  Note first that
$$
Z \times _{S|_{\Xx /Z}}(R|_{\Xx /Z}) \rightarrow R|_{\Xx /Z}\rightarrow
S|_{\Xx /Z}
$$
is a fibration sequence (this should actually have been pointed out above in the
treatment of the $\pi _i$, $i\geq 2$), over the basepoint $s\in S(Z)$.
On the other hand note that $r': Z\rightarrow R$ is a point lying over $s'$.
Consider the map
$$
\varpi _1(S|_{\Xx /Z}, s, s')\rightarrow
\pi _0(Z \times _{S|_{\Xx /Z}}(R|_{\Xx /Z}))
$$
which sends a path to the point obtained by transporting $f'_Z$ along the path
from $s'$ back to $s$.
The fibration sequence gives the following statement:

{\em The group $\pi _1(Z \times _{S|_{\Xx /Z}}(R|_{\Xx /Z}, f_Z)$ acts on
$\varpi _1(R|_{\Xx /Z}, r, r')$ with quotient the inverse image in
$\varpi _1(S|_{\Xx /Z}, s, s')$ of the section $f_Z: Z \rightarrow \pi _0(
Z \times _{S|_{\Xx /Z}}(R|_{\Xx /Z})$.
}

Now we note that
$$
\pi _1(Z \times _{S|_{\Xx /Z}}(R|_{\Xx /Z}, f_Z) =
\pi _1(X\times _SR|_{\Xx /Z}, f),
$$
and
$$
\pi _0(Z \times _{S|_{\Xx /Z}}(R|_{\Xx /Z})\subset \pi _0(X\times _S
R|_{\Xx /Z}).
$$
The transport of $f'$ along the path from $s'$ to $s$ again gives a map
$$
\varpi _1(S|_{\Xx /Z}, s, s')\rightarrow
\pi _0(X\times _SR|_{\Xx /Z})
$$
and we obtain the following statement.

{\em The group $\pi _1(X\times _SR|_{\Xx /Z}, f)$ acts on
$\varpi _1(R|_{\Xx /Z}, r, r')$ with quotient the inverse image in
\linebreak
$\varpi _1(S|_{\Xx /Z}, s, s')$ of the section $f: Z \rightarrow \pi
_0( X\times _SR|_{\Xx /Z})$.
}

Now we look at everything in terms of sheaves on $\Xx /Y$.
Let $Res _{Z/Y}$ denote the restriction from $Z$ down to $Y$, and let
$\tilde{f}$ denote the $Y$-valued point corresponding to $f$.
Note that
$$
Res _{Z/Y} \pi _0(X\times _SR|_{\Xx /Z}) = \pi _0(X\times _SR|_{\Xx /Y})\times
_YZ.
$$
In general if $\Aa$ is a sheaf over $Z$ and $\Bb$ a sheaf over $Y$ with a
section $Y\rightarrow \Bb$ then
$$
Res _{Z/Y}(\Aa \times _{\Bb |_{\Xx /Z}}Z) = (Res _{Z/Y}\Aa )\times _{\Bb}Y.
$$
In particular the inverse image in $\varpi _1(S|_{\Xx /Z}, s, s')$ of the
section $f: Z \rightarrow \pi _0(
X\times _SR|_{\Xx /Z})$ restricts down to $Y$ to the inverse image in
$Res _{Z/Y}\varpi _1(S|_{\Xx /Z}, s, s')$ of the section
$\tilde{f}: Y \rightarrow \pi _0(X\times
_SR|_{\Xx /Y})$.

Another general principal is that if $\Gg$ is a group sheaf on $Y$ such that
$\Gg |_{\Xx /Z}$ acts on a sheaf $\Hh$ then $\Gg$ acts on $Res _{Z/Y}\Hh$ with
quotient equal to $Res _{Z/Y}(\Hh /(\Gg |_{\Xx /Z}))$.

With these things in mind,
our above statement becomes:

{\em The group $\pi _1((X\times _SR |_{\Xx /Y}, \tilde{f})$ acts on
$Res _{Z/Y}\varpi _1(R|_{\Xx /Z}, r, r')$ with quotient the inverse image in
$Res _{Z/Y}\varpi _1(S|_{\Xx /Z}, s, s')$ of the section
$\tilde{f}: Y \rightarrow \pi _0(X\times
_SR|_{\Xx /Y})$.
}

Now the facts that $\pi _0(X\times
_SR|_{\Xx /Y})$ is $P3\frac{1}{2}$ and that
$Res _{Z/Y}\varpi _1(S|_{\Xx /Z}, s, s')$ is $P4$  (which comes by
hypothesis)
imply that the inverse image
in question is $P4$ (Lemma \ref{kernel}); then the theorem on quotients (Theorem
\ref{I.1.d}) and the fact that the group  $\pi _1((X\times _SR |_{\Xx /Y},
\tilde{f})$ is $P5$ over $Y$ gives the condition that  $Res _{Z/Y}\varpi
_1(R|_{\Xx /Z}, r, r')$ is $P4$ over $Y$.

This is the condition on $\varpi _1$  needed to insure that $R$ is presentable.
This completes the proof of Theorem \ref{stability}.
\eop

We have the following characterization of presentable morphisms via the relative
homotopy group sheaves.

\begin{proposition}
\mylabel{characterization}
Suppose $f: R\rightarrow S$ is a morphism of $n$-stacks.  Then
$f$ is presentable (resp. very presentable) if and only if the following
conditions are satisfied for any scheme $X$ of finite type:
\newline
---for any  morphism $X\rightarrow S$,
the sheaf $\pi
_0(X\times _SR)$ is $P3\frac{1}{2}$; and
\newline
---for any morphism $r: X\rightarrow R$ the sheaves $\pi _i(X\times _SR/X, r)$
on $\Xx /X$ are presentable group sheaves over $X$ (resp. $\pi _1$ is affine
presentable and $\pi _i$ are vector sheaves for $i\geq 2$).
\end{proposition}
{\em Proof:}
This falls out of the  proof of \ref{stability}.
\eop

{\em Exercise:}  For which values of $(a_0,a_1,\ldots )$ does Theorem
\ref{I.1.s.?} hold for $(a_0,a_1,\ldots )$-presentable spaces?  Place these
conditions in Corollary \ref{I.1.u} below.

\subnumero{Going to the base of a fibration}

It is an interesting question to ask, if $R\rightarrow S$ is a morphism
of $n$-stacks such that $R$ is presentable and such that for every scheme-valued
point $X\rightarrow S$ the fiber product $X\times _SR$ is presentable, then
is $S$ presentable?  The answer is surely no in this generality. We need to make
additional hypotheses.  Directly from the fibration exact sequences, one can see
that if $\pi _0(S)$ is assumed to be $P3\frac{1}{2}$ (a hypothesis which seems
unavoidable) and if we suppose that for any point $a:X\rightarrow S$, the action
of $\pi _1(S|_{\Xx /X}, a)$ on $\pi _0(X\times _SR)$ factors through a
presentable group sheaf over $X$, then $S$ will be  presentable.

As a particular case, if the morphism $R\rightarrow S$ is relatively
$0$-connected (i.e. the fibers are connected) and surjective on $\pi _0$, then
presentability of $R$ implies presentability of $S$.

One might look for other weaker conditions, for example that the fibers satisfy
some sort of artinian condition (e.g. there is a surjection from a scheme finite
over $X$, to $\pi _0(X\times _SR)$).  I don't know if this can be made to work.

\subnumero{Presentable shapes}

We have a notion of internal $Hom$ for $n$-stacks.  In the topological presheaf
interpretation (\cite{kobe} \S 2), recall that $\underline{Hom}(R,T)$ is
defined to be the presheaf $X\mapsto Hom (R'_X,T|_{\Xx /X})$ where $R'_X$ is a
functorial replacement of $R|_{\Xx /X}$ by a cofibrant presheaf.

\begin{corollary}
\label{I.1.u}
Suppose $W$ is a finite CW complex, and let $W_{\Xx}$ denote the constant
$n$-stack with values $\Pi _n(W)$ (or in terms of presheaves of spaces,
it is the fibrant presheaf associated to the constant presheaf with values
$\tau _{\leq n}W$). If $T$ is a presentable (resp. very presentable)
$n$-stack over $X$
then the $n$-stack $\underline{Hom}(W_{\Xx}, T)$ is presentable (resp.
very presentable).
\end{corollary}
{\em Proof:}
We first show this for $W=S^m$, the $m$-sphere.  Do this by induction on $m$.
It is clear for $m=0$ because then $W$ consists of two points and
$\underline{Hom}(W_{\Xx}, T)=T\times T$. For any $m$, write $S^m$ as the
union of two copies of $B^m$ joined along $S^{m-1}$.  We get
$$
\underline{Hom}(S^m_{\Xx}, T)=T\times _{
\underline{Hom}(S^{m-1}_{\Xx}, T)}T,
$$
since $\underline{Hom}(B^m_{\Xx}, T)=T$.
By Theorem \ref{I.1.s.?}, $\underline{Hom}(S^m_{\Xx}, T)$ is presentable
(resp. very presentable).  This shows the corollary for the spheres.

We now treat the case of general $W$, by induction on the number of cells.  We
may thus write $W=W'\cup B^m$ with the cell $B^m$ attached over an attaching
map $S^{m-1}\rightarrow W'$, and where we know the result for $W'$. Then
$$
\underline{Hom}(W_{\Xx}, T)=\underline{Hom}(W'_{\Xx}, T)\times _{
\underline{Hom}(S^{m-1}_{\Xx}, T)}T.
$$
Again by Theorem \ref{I.1.s.?}, we obtain the result for $W$.
\eop

Let $Pres ^n/\Xx$ denote the $n+1$-category of presentable
$n$-stacks.  We define the {\em presentable shape} of $W$ to be the
$n+1$-functor
$$
Shape (W):T\mapsto \underline{Hom}(\underline{W}, T)
$$
from $Pres ^n/\Xx$ to $Pres ^n/\Xx$.

One can show (using the calculations of \cite{kobe} Corollary 3.9 over $S=Spec
(k)$) that if $W$ is connected and simply connected then this functor is
homotopy-representable by an object $Hull (W)\in Pres /\Xx $. On the other hand,
in most cases where $W$ is not simply connected, the presentable shape is not
representable.  We could try to interpret  the hull of $W$ as the inverse limit
of $Shape (W)$, but this is not a standard kind of inverse limit.  It is a
question for further study, just what information is contained in $Shape (W)$.

{\em Example:}  Take $G=GL(n)$ and $T= K(G, 1)$.
Fix a finite CW complex $U$.
Then $M:=\underline{Hom}(U, T)$
is the moduli stack for flat principal $G$-bundles (i.e. flat vector bundles of
rank $n$) on $U$.

More generally it should be interesting to look at presentable or
very presentable {\em connected} $T$, these are objects whose homotopy group
sheaves are algebraic Lie groups over $Spec (k)$. Note that if $k$ is
algebraically closed then there is an essentially unique choice of basepoint
$t\in T(Spec (k))$. If $G= \pi _1(T, t)$ then we have a fibration
$T\rightarrow K(G,1)$ and we get a morphism
$$
\underline{Hom} (U, T) \rightarrow \underline{Hom}(U, K(G, 1)).
$$
This expresses $\underline{Hom} (U, T)$ as a presentable $n$-stack over the
moduli stack $M$ of flat principal $G$-bundles over $U$.

One can see from this example that we should consider the notion of vector
sheaf as a candidate for the higher homotopy group sheaves.

\subnumero{Leray theory}

We develop here a nonabelian Leray theory and K\"unneth formula.  This is in
some sense one of the principal reasons for going to nonconnected $n$-stacks, as
they can intervene as intermediate steps even when the original coefficient
stacks
were connected.

We give some notation for the stack of sections.  If $T\rightarrow S$ is a
morphism of  $n$-stacks on $\Xx$ (or on any site) then we denote by
$\underline{\Gamma}(S, T)$ the $n$-stack of sections, i.e. of diagrams
$$
\begin{array}{ccc}
S&\rightarrow &T\\
& {\displaystyle =}\searrow&\downarrow \\
& & S
\end{array}
$$
(with homotopy making the diagram commutative).

We also have a notion of relative morphism stack.  Suppose that $T\rightarrow S$
and $T' \rightarrow S$ are two morphisms of $n$-stacks.  Then we obtain an
$n$-stack together with morphism to $S$
$$
\underline{Hom}(T/S, T'/S) \rightarrow S.
$$
In topological language this corresponds to the space whose fiber over $s$
is the
space of morphisms from $T_s$ to $T'_s$.  This should not be confused with
another useful construction in the same situation, the space
$$
\underline{Hom}_S(T, T')
$$
which is the $n$-stack of diagrams
$$
\begin{array}{ccc}
T&\rightarrow &T'\\
& \searrow&\downarrow \\
& & S
\end{array}
$$
(again with homotopy making the diagram commutative).

These things can be constructed using the point of view of simplicial
presheaves or presheaves of spaces---cf for example \cite{flexible}.
It remains to
be seen how to give constructions of these things purely within the
realm of stacks
(and consequently to extend the same constructions to stacks of $n$-categories
which are not necessarily $n$-groupoids).

We have the following relationships among the above constructions.
First of all,
$\underline{\Gamma}(S, T) = \underline{Hom}_S(S, T)$. Then,

\begin{lemma}
Suppose $T\rightarrow S$ and $T' \rightarrow S$ are morphisms of $n$-stacks.
There is a natural equivalence
$$
\underline{\Gamma} (S, \underline{Hom}(T/S, T'/S)) \cong
\underline{Hom}_S(T,T').
$$
\end{lemma}
{\em Proof:}  From the point of view of presheaves of spaces, see
\cite{flexible}.
\eop

Finally note that if $T$ is an $n$-stack and $R\rightarrow S$ is a morphism of
$n$-stacks then
$$
\underline{Hom}_S(R/S, T\times S/S) \cong \underline{Hom}(R, T).
$$
From the above lemma we obtain a method of ``devissage'':
\begin{corollary}
Suppose $T$ is an $n$-stack and $R\rightarrow S$ is a morphism of
$n$-stacks, then
$$
\underline{Hom}(R, T) \cong \underline{\Gamma} (S, \underline{Hom}(R/S, T\times
S/S)).
$$
\end{corollary}
\vspace*{-.5cm}
\eop

In words this says that to calculate the stack of morphisms from $R$ to $T$ we
first look at the fiberwise morphisms from $R/S$ to $T$, and then we take the
sections over $S$.
Rather than taking the internal morphism and section spaces we can take the
external ones, removing the underline in the notation which means taking the
sections over $\ast$ (which is $Spec (k)$ in our case). We get
the statement
$$
Hom(R, T) \cong \Gamma (S, \underline{Hom}(R/S, T\times
S/S)).
$$
Note that it is still essential to look at the internal $\underline{Hom}$ inside
the space of sections.

It might be worthwhile looking at how this works in the case
of usual cohomology.  Suppose $\Aa$ is a sheaf of abelian groups on $\Xx$.
Let $T= K(\Aa , n)$, so that $Hom(R, T)$ is an $n$-groupoid with
homotopy groups
$$
\pi _i = H^{n-i}(R, \Aa ).
$$
Similarly $\underline{Hom}(R/S, T)$ is an $n$-stack over $S$ whose relative
homotopy group sheaves over $S$ are the higher direct images
$$
\pi _i = R^{n-i}f_{\ast} (\Aa |_R).
$$
There is a spectral sequence for the $n$-stack of sections going from the
cohomology of $S$ with coefficients in the relative homotopy sheaves to the
homotopy groups  of the space of sections, which turns out to be the Leray
spectral
sequence in this case.

This version of Leray theory is due to
Thomason \cite{Thomason}, who developed it mostly
in the context of presheaves of spectra.

We finally introduce one more bit of notation combining the previous notations,
that is the {\em relative section stack}.  Suppose $R\rightarrow S\rightarrow T$
are morphisms of $n$-stacks.  Then we obtain the $n$-stack
$$
\underline{\Gamma}(S/T, R/T)\rightarrow T
$$
which is geometrically the ``fiberwise space of sections of the morphism
$R\rightarrow S$ along the fibers of $S\rightarrow T$''.
The above Leray theory can itself be presented in a relative context:

\begin{lemma}
\mylabel{RelativeLeray}
Suppose $R\rightarrow S\rightarrow T\rightarrow U$ are morphisms of $n$-stacks.
Then
$$
\underline{\Gamma}(T/U, \underline{\Gamma}(S/T, R/T)/U) \cong
\underline{\Gamma}(T/U, R/U).
$$
\end{lemma}
\eop

Of course, given four morphisms there should be a diagram expressing
compatibility
of these Leray equivalences (and  further diagrams of homotopy between
the homotopies).

\subnumero{Leray theory for presentable and very presentable $n$-stacks}
Now we get back to presentable and very presentable $n$-stacks.  Our goal  is to
show that in certain cases the Leray theory stays within the world of
presentable
$n$-stacks.

The first task is to generalize Corollary
\ref{I.1.u} to the case of a local coefficient system, i.e. a presentable
morphism of $n$-stacks to our given finite CW complex.

\begin{lemma}
\mylabel{Leray2}
Suppose $U$ is a constant $n$-stack associated to the $n$-groupoids
associated to a finite CW complex.
Suppose $T\rightarrow U$ is a presentable (resp. very presentable) morphism of
$n$-stacks. Then the $n$-groupoid of sections
$\underline{\Gamma}(U, T)$ is a presentable (resp. very presentable) $n$-stack.
\end{lemma}
{\em Proof:}
The proof is identical to that of Corollary \ref{I.1.u} but we repeat it here
for the reader's convenience.
As before, we first treat the case $U=S^m$ by induction on $m$. It is clear for
$m=0$ because then $W$ consists of two points $a,b$ and  $\underline{\Gamma
}(W_{\Xx}, T)=T_a\times T_b$, with the fibers $T_a$ and $T_b$ being presentable
(resp. very presentable). Now for any $m$, write  $S^m$ as the union of two
copies of $B^m$ joined along $S^{m-1}$ and let $T_a$ be the fiber of $T$ over a
basepoint. This fiber is presentable (resp. very presentable).  We get
$$
\underline{\Gamma }(S^m_{\Xx}, T)=T_a\times _{
\underline{\Gamma }(S^{m-1}_{\Xx}, T)}T_a,
$$
since $\underline{\Gamma }(B^m_{\Xx}, T)\cong T_a$.
By the induction hypothesis and Theorem \ref{I.1.s.?}, $\underline{\Gamma
}(S^m_{\Xx}, T)$ is presentable (resp. very presentable).  This shows the
lemma for the spheres.

We now treat the case of general $U$, by induction on the number of cells.  We
may thus write $U=U'\cup B^m$ with the cell $B^m$ attached over an attaching
map $S^{m-1}\rightarrow U'$, and where we know the result for $U'$. Again let
$T_a$ be the fiber over a basepoint in the attached cell. Then
$$
\underline{\Gamma }(U_{\Xx}, T)=\underline{\Gamma }(U'_{\Xx}, T)\times _{
\underline{\Gamma }(S^{m-1}_{\Xx}, T)}T_a.
$$
By Theorem \ref{I.1.s.?} and the above result for spheres, we obtain the result
for $U$.
\eop

Say that a morphism $U\rightarrow V$ of $n$-stacks is {\em of finite CW type}
if for any scheme of finite type $X$ with morphism $X\rightarrow V$ there
is a covering family $\{ Y_{\alpha} \rightarrow X\}$ and finite CW complexes
$W^{\alpha}$ such that $Y_{\alpha} \times _V U \cong Y_{\alpha} \times
W^{\alpha}_{\Xx}$ (with $W^{\alpha}_{\Xx}$ being the constant $n$-stack
associated to $\Pi _n(W^{\alpha})$ as defined previously).

\begin{theorem}
\mylabel{Leray}
Suppose $U\rightarrow V$ is a morphism of $n$-stacks of finite CW type,
and suppose $T\rightarrow U$ is a presentable (resp. very presentable) morphism
of $n$-stacks. Then $\underline{\Gamma} (U/V, T/V)\rightarrow V$ is a
presentable (resp. very presentable) morphism.
\end{theorem}
{\em Proof:}
Suppose $X$ is a scheme of finite type with a morphism $X\rightarrow V$.
Let $\{ Y^{\alpha} \rightarrow X\}$ be the covering family and $\{ W^{\alpha}\}$
the collection of finite CW complexes
with isomorphisms $U\times _VY^{\alpha}\cong W^{\alpha}_{\Xx}$
given by the fact
that $U\rightarrow V$ is a morphism of finite CW type.  It suffices to
prove that
$$
\underline{\Gamma} (U/V, T/V)\times _V Y^{\alpha}=
\underline{\Gamma} (U\times _VY^{\alpha}/Y^{\alpha}, T\times
_VY^{\alpha}/Y^{\alpha})
$$
is presentable (resp. very presentable).  Thus it suffices to prove the theorem
in the case where $V$ is a scheme of finite type and $U=V\times W_{\Xx}$ for
a finite CW complex $W$. With these hypotheses we return to the notations of the
theorem. If $W$ is a finite union of components then the section space in
question will be the product of the section spaces of each of the components.
Thus we may assume that $W$ is connected.  The $n$-stack of sections from
$W_{\Xx}$ to $V\times W_{\Xx}$ is isomorphic to $V$.  Thus the $n$-stack of
sections of the morphism $T\rightarrow W_{\Xx}$ maps to $V$, and this $n$-stack
of sections is the same as the relative section stack $\underline{\Gamma}(U/V,
T/V)$.
It suffices to prove that $\underline{\Gamma}(W_{\Xx}, T)$ is presentable (resp.
very presentable).
But the morphism $V\times W_{\Xx}\rightarrow W_{\Xx}$ is very
presentable,  so by Corollary \ref{composition} the morphism $T\rightarrow
W_{\Xx}$ is presentable (resp. very presentable), and Lemma \ref{Leray2}
applies to give that $\underline{\Gamma}(W_{\Xx}, T)$ is presentable (resp.
very presentable) as needed.
\eop

\begin{corollary}
\mylabel{Leray1}
Suppose $U\rightarrow V$ is a morphism of $n$-stacks of finite CW type
and
Suppose $T\rightarrow V$ is a presentable morphism of
$n$-stacks. Then the morphism
$$
\underline{Hom}(U/V, T/V)\rightarrow V
$$
is a presentable morphism.
\end{corollary}
{\em Proof:}
We have
$$
\underline{Hom}(U/V, T/V) = \underline{\Gamma}(U/V, T\times _VU/V)
$$
and $T\times _VU\rightarrow U$ is presentable by \ref{fiberprod}, so Theorem
\ref{Leray} applies.
\eop

\begin{corollary}
\mylabel{Leray1a}
Suppose $T$ is a presentable $n$-stack, and suppose $V\rightarrow U$ is a
morphism whose fibers are finite CW complexes in the sense of the above
theorem.
Then
$$
\underline{Hom}(V/U, T\times U/U)\rightarrow U
$$
is a presentable morphism.
\end{corollary}
{\em Proof:}
Indeed, the morphism $T\times U\rightarrow U$ is presentable.
\eop

We look at the case of a  morphism of $n$-groupoids
$f:U\rightarrow V$ such that $U$ and $V$ are the $n$-groupoids associated to
finite
CW complexes.  Suppose that the fibers of $f$ are
the $n$-groupoids associated to finite CW complexes. This is the case for
example
if $f$ comes from a smooth morphism of manifolds.
For a presentable
$n$-stack $T$ we can calculate
$$
\underline{Hom}(V, T) = \underline{\Gamma}(U, \underline{Hom}(V/U,
T\times U/U)).
$$
Corollary \ref{Leray1a} states that
$\underline{Hom}(V/U, T\times U/U)\rightarrow U$ is a presentable morphism, and
Lemma \ref{Leray2} (which is also a corollary of Theorem \ref{Leray})
states that
for any presentable morphism $R\rightarrow U$ the space of sections is
presentable.  We obtain in particular the presentability of  $\underline{Hom}(V,
T)$ (which we already knew beforehand).  The Leray devissage process thus stays
within the realm of presentable $n$-stacks.

{\em The K\"unneth formula:}  We can apply the above discussion to the
particular
case where $V=U\times U'$ is a product.  In this case the formula is
simplified:
$$
\underline{Hom}(U\times U', T) = \underline{Hom}(U, \underline{Hom}(U', T))
$$
and again (this time using only Corollary \ref{I.1.u})
this process of first taking $\underline{Hom}(U', T)$ and then
$\underline{Hom}(U, -)$ stays within the realm of presentable $n$-stacks.

Of course the entire discussion above works equally well if we replace
``presentable'' by ``very presentable''.

{\em Example:}  Take $G=GL(n)$ and $T= K(G, 1)$.
Then $M':=\underline{Hom}(U', T)$
is the moduli stack for flat principal $G$-bundles (i.e. flat vector bundles of
rank $n$) on $U'$.  After that, assuming that $U$ is connected,
$\underline{Hom}(U,
M')$ is the moduli stack of flat $G$-bundles on $U\times U'$.

More generally it should be interesting to look at presentable or
very presentable {\em connected} $T$, these are objects whose homotopy group
sheaves are algebraic Lie groups over $Spec (k)$. Note that if $k$ is
algebraically closed then there is an essentially unique choice of basepoint
$t\in T(Spec (k))$. If $G= \pi _1(T, t)$ then we have a fibration
$T\rightarrow K(G,1)$ and we get a morphism
$$
\underline{Hom} (U, T) \rightarrow \underline{Hom}(U, K(G, 1)).
$$
This expresses $\underline{Hom} (U, T)$ as a presentable $n$-stack over the
moduli stack $M$ of flat principal $G$-bundles over $U$.

\end{document}